\documentclass{article}
\usepackage{amsmath}
\usepackage{amsthm}
\usepackage{amsfonts}
\usepackage{amssymb}
\usepackage{geometry}
\usepackage{algorithm}
\usepackage{algpseudocode}
\usepackage{color}
\usepackage{enumitem}
\usepackage{graphicx}

\newtheorem{theorem}{Theorem}
\newtheorem{proposition}{Proposition}
\newtheorem{lemma}{Lemma}
\newtheorem{definition}{Definition}
\newtheorem{assumption}{Assumption}

\newtheorem{remark}{Remark}

\newtheorem*{standingassumption}{Standing Assumption}

\DeclareMathOperator{\dist}{dist}
\DeclareMathOperator{\range}{range}
\DeclareMathOperator{\Null}{Null}

\title{A Proximal Newton Framework with Local Superlinear Convergence for Simple Bilevel Optimization}
\author{
	Xiaoning Bai\thanks{Department of Mathematics, Southern University of Science and Technology, Shenzhen, China. Email: \texttt{12331002@mail.sustech.edu.cn}.}
	\and
	Boris S. Mordukhovich\thanks{Department of Mathematics \& Center for Artificial Intelligence and Data Sciences, Wayne State University, Detroit, MI 48202, USA. Email: \texttt{aa1086@wayne.edu}.} 
	\and
	Shangzhi Zeng\thanks{National Center for Applied Mathematics Shenzhen \& Department of Mathematics, Southern University of Science and Technology, Shenzhen, China. Email: \texttt{zengsz@sustech.edu.cn}.}
	\and
	Jin Zhang\thanks{Department of Mathematics \& National Center for Applied Mathematics Shenzhen, Southern University of Science and Technology, Shenzhen, China. Email: \texttt{zhangj9@sustech.edu.cn}.}
}
\date{}

\begin{document}
	\maketitle
	\begin{abstract}
		In this paper, we study a penalty method for simple bilevel optimization in which the lower-level objective has a convex composite structure. Our approach is based on a vanishing-penalty reformulation, which gives rise to a sequence of penalized problems as the penalty parameter decreases. A key difficulty is that metric subregularity of the lower-level optimality mapping does not, in general, yield a uniform error bound for the penalized problems as the penalty parameter vanishes, preventing a direct application of standard proximal Newton theory. To address this difficulty, we first analyze a regularized proximal Newton step by exploiting the specific structure of simple bilevel problems, and establish the local estimates that underpin superlinear convergence. Building on this local analysis, we propose SiB-Newton, which incorporates the proximal Newton step into a globalization strategy. We establish convergence guarantees for SiB-Newton and show that, under suitable local conditions, Newton steps are eventually accepted and the iterates converge at a superlinear rate. Numerical experiments on synthetic examples verify the predicted local superlinear convergence rate.
	\end{abstract}
	\section{Introduction}

	Bilevel optimization has recently attracted significant attention for modeling hierarchical decision-making structures in optimization and machine learning. Among the various bilevel optimization models, a special class is simple bilevel optimization, in which both the upper- and lower-level objectives are optimized over the same variable. This class of problems arises in a variety of applications, including sparse linear regression \cite{jiang2023conditional, samadi2024achieving}, ill-posed least-squares regression \cite{chen2024penalty, amini2019iterative}, linear inverse problems such as integral equation problems \cite{beck2014first, shehu2021inertial, jiang2025near}, and logistic regression \cite{wang2024near, latafat2025convergence}. In this paper, we consider the following simple bilevel convex optimization problem:
	\begin{equation}\label{simple}
		\min_{x \in \mathbb{R}^n} F(x) \quad \text{s.t.} \quad x \in \mathop{\arg\min}_{x \in \mathbb{R}^n} \varphi(x):=h(Ax) + g(x).
	\end{equation}
	Here, \(A\in\mathbb{R}^{m\times n}\), \(F:\mathbb{R}^n\to\mathbb{R}\) and \(h:\mathbb{R}^m\to\mathbb{R}\) are smooth convex functions, and \(g:\mathbb{R}^n\to\mathbb{R}\cup\{+\infty\}\) is convex but potentially nonsmooth.

	Numerous studies have investigated algorithms for simple bilevel optimization problems such as \eqref{simple}, with a particular focus on first-order methods. Among them, Beck and Sabach \cite{beck2014first} provided one of the earliest complexity results, introducing a projected-gradient method that integrates the upper-level gradient with the optimality condition of the lower level to achieve an $O(1/k)$ convergence rate under convexity. Related developments include fixed-point and inertial variants \cite{sabach2017first, shehu2021inertial}, as well as the Bi-Sub-Gradient method \cite{merchav2023convex}. Alternative first-order frameworks include alternating feasibility-descent steps \cite{helou2017}, projection-free conditional-gradient schemes that replace projections with linear optimization oracles \cite{jiang2023conditional, giang2025projection}, and online convex optimization formulations that utilize regret minimization to yield simultaneous guarantees for both objectives \cite{shen2023online}. More recently, acceleration via cutting-plane surrogates and Nesterov updates has been shown to achieve the optimal $O(1/k^2)$ rate in the smooth convex setting \cite{cao2024accelerated}. Another line of research focuses on bisection-based approaches. Wang et al. \cite{wang2024near} developed a bisection method with near-optimal complexity, matching the optimal first-order complexity for smooth or composite convex optimization up to logarithmic factors. Jiang et al. \cite{jiang2025near} further reformulated the problem as a one-dimensional root-finding task and developed bisection and Newton-type methods, both achieving a near-optimal operation complexity of $\widetilde{O}(\epsilon^{-1/2})$.

	Among the existing works, a major line of research in simple bilevel optimization is based on penalty techniques. The key idea is to approximate the original bilevel problem by penalizing the lower-level objective $\varphi(x)$. In this paper, we consider the following penalized subproblem:
	\begin{equation}\label{eq:penalty}
		\min_{x \in \mathbb{R}^n} \psi_{\sigma_k}(x) := \sigma_k F(x) + \varphi(x),
	\end{equation}
	where the sequence of parameters satisfies $\sigma_k > 0$ and $\sigma_k \to 0$ as $k \to \infty$.

	This perspective closely links the penalty method for simple bilevel optimization with classical Tikhonov regularization \cite{tikhonov1977solutions}. In standard Tikhonov regularization, an ill-posed minimization problem is stabilized by adding a strongly convex regularizer, scaled by a parameter that asymptotically approaches zero. Similarly, in formulation \eqref{eq:penalty}, the upper-level objective $F(x)$ acts as a generalized regularizer. This mechanism is fundamentally connected to the viscosity approximation methods pioneered in optimization and variational inequalities \cite{attouch1996viscosity, facchinei2003finite}.

	A number of first-order algorithms have been developed for simple bilevel optimization based on penalty formulations such as \eqref{eq:penalty}. Early contributions include the explicit penalty methods of Solodov \cite{solodov2007explicit, solodov2007bundle}, which established $O(1/k)$ convergence rates under convexity. Amini et al. \cite{amini2019iterative} proposed an iteratively regularized incremental subgradient method with sublinear complexity. Subsequent works investigated accelerated schemes \cite{samadi2024achieving}, achieving an $O(1/k^2)$ rate for the upper-level objective with a fixed regularization parameter, and the same rate for the lower level under weak sharpness. Chen et al. \cite{chen2024penalty} established a unified penalty framework under $q$-H\"olderian error bounds, deriving an $O(1/k^{1/q})$ convergence rate. More recently, \cite{samadi2025iteratively} proposed iteratively regularized schemes, proving linear convergence under strong convexity and weak sharpness, while guaranteeing sublinear rates in the general convex case. For convex composite simple bilevel problems, Latafat et al. \cite{latafat2025convergence} developed proximal-gradient schemes with adaptive and nonadaptive stepsize rules. Merchav et al. \cite{merchav2026dynamic} introduced a dynamic FISTA scheme, achieving an $O(1/k^2)$ rate for the lower-level objective and improved simultaneous convergence guarantees under an error-bound condition. Shtern and Taiwo \cite{shtern2026convergence} further established simultaneous convergence rates for the lower- and upper-level objectives and characterized the trade-off induced by the regularization schedule. In parallel, \cite{boct2026accelerating} adopted a continuous-time perspective to analyze diagonal (Tikhonov-type) methods, deriving unified convergence guarantees for both proximal-gradient and accelerated variants under H\"olderian growth and the weaker Attouch-Czarnecki condition, extending the analysis to general Hilbert space settings.

	Despite this substantial progress in first-order methods for simple bilevel optimization \eqref{simple}, to the best of our knowledge, no algorithm has been developed for problem~\eqref{simple} with a guarantee of local superlinear convergence. To bridge this gap, we investigate Newton-type methods, a classical framework renowned for achieving quadratic local convergence in smooth, unconstrained settings \cite{nocedal2006numerical, broyden1970convergence, goldfarb1970family, broyden1967quasi}. Specifically, when the objective features a convex composite structure,
	\begin{equation*}
		\min_{x\in\mathbb{R}^n} \phi(x)+g(x),
	\end{equation*}
	where the function decomposes into a smooth term $\phi(x)$, proximal Newton methods extend classical updates by solving a proximalized quadratic subproblem at each iteration:
	\begin{equation*}
		\min_{x \in \mathbb{R}^n}\, \langle\nabla \phi(x_k), x-x_k\rangle + \frac{1}{2}\langle x-x_k, \nabla^2 \phi(x_k)(x-x_k)\rangle + g(x).
	\end{equation*}
	These methods achieve superlinear local convergence in structured settings such as $\ell_1$-regularization \cite{lee2012proximal} and have been successfully extended to more general nonsmooth terms \cite{lee2014proximal, byrd2016inexact}. To secure global guarantees, strategies such as line search and trust-region frameworks have been widely deployed \cite{lee2019inexact, scheinberg2016practical, yue2019family}. More recent refinements leverage variational analytic techniques; for instance, \cite{mordukhovich2023globally} utilized metric $q$-subregularity to characterize local rates, while \cite{vom2024inexact} proposed an inexact regularized proximal Newton scheme without line search under H\"olderian error bounds, and \cite{garber2025weak} introduced weak proximal Newton oracles.

	Unlike standard proximal Newton methods for a fixed optimization problem, penalty-based approaches to simple bilevel optimization involve a sequence of penalized problems with \(\sigma_k\to0\). This varying problem structure makes the local Newton analysis more delicate. In this work, we develop a proximal Newton framework for this setting and show that superlinear local convergence can still be achieved. Our analysis exploits variational analysis techniques with the special structure of simple bilevel problems to directly address the resulting difficulties in the local analysis of proximal Newton methods.

	A central difficulty is that the penalized objective varies with the penalty parameter. Since metric subregularity is not stable under perturbations in general \cite{zheng2021perturbation,dontchev2004regularity,benko2018new,dontchev2020radius}, metric subregularity of the lower-level optimality mapping does not, in general, yield a uniform error bound for the family of penalized problems. This prevents a direct application of standard proximal Newton convergence analyses. 

	Our main contributions are summarized as follows. First, we exploit the structure of simple bilevel problems to derive local estimates for the Newton update directly relative to the solution set of problem~\eqref{simple}, without requiring uniform error bounds for the penalized problems. By suitably coupling the penalty, regularization, and inexactness parameters, we establish local superlinear convergence of the resulting proximal Newton method. Second, building on this local analysis, we propose SiB-Newton, which combines proximal Newton unit steps with a fallback line-search procedure. We establish global accumulation-point guarantees and show that, under suitable local assumptions, Newton unit steps are eventually accepted, leading to local superlinear convergence of the globalized method. Numerical experiments on a linear regression problem show that SiB-Newton reaches smaller errors more rapidly than the compared methods, consistent with the predicted local superlinear convergence behavior. Our analysis is carried out under the following standing assumptions on the problem data.
	\begin{standingassumption}\label{asup_upper_1}
		Impose the following structural properties on the simple bilevel problem~\eqref{simple}:
		\begin{enumerate}[label=(\alph*)]
			\item The functions \(F:\mathbb R^n\to\mathbb R\) and \(h:\mathbb R^m\to\mathbb R\) are convex and twice continuously differentiable, and \(g:\mathbb R^n\to\mathbb R\cup\{+\infty\}\) is proper, closed, and convex.

			\item The lower-level solution set \(X_\varphi^*:=\mathop{\arg\min}_{x\in\mathbb{R}^n}\varphi(x)\) is nonempty. 

			\item The gradients \(\nabla F\) and \(\nabla h\) are Lipschitz continuous on \({\rm dom}\,g\) and \(\mathbb{R}^m\), respectively, with Lipschitz constants \(L_F>0\) and \(L_h>0\).
			
			\item The upper-level objective \(F\) is bounded below on
			\(\mathbb R^n\), and the simple bilevel problem~\eqref{simple}
			admits at least one optimal solution.
		\end{enumerate}
	\end{standingassumption}
	The remainder of the paper is organized as follows. Section~\ref{sec2} introduces the notation and basic variational-analytic tools used in the subsequent analysis. Section~\ref{sec3} develops the local analysis of the proximal Newton step and establishes the estimates underlying superlinear convergence. Section~\ref{sec4} introduces the globalization strategy and establishes the corresponding global convergence and local superlinear convergence results. Section~\ref{sec5} presents numerical experiments illustrating the performance of SiB-Newton.

	\section{Preliminaries}\label{sec2}
	This section introduces the notation, key definitions, and results from variational analysis that are foundational to our subsequent development. We begin by introducing the notation used throughout the paper. Let $\mathbb{R}^n$ denote the $n$-dimensional Euclidean space, and let $\mathbb{B}_\epsilon(x)$ denote the Euclidean ball centered at $x$ of radius $\epsilon$, i.e., $\mathbb{B}_\epsilon(x) := \left\{ z \in \mathbb{R}^n : \|z-x\|_2 < \epsilon \right\}$. The standard inner product and Euclidean norm are denoted by $\langle \cdot, \cdot \rangle$ and $\|\cdot\|$, respectively. For any matrix \(A \in \mathbb{R}^{m \times n}\), we denote by \(\|A\|\) its spectral norm, i.e., $\|A\| := \max_{\|x\| \neq 0} \|Ax\|/\|x\|$. For a proper, closed, convex function $\Phi$, the proximal mapping of \( \Phi: \mathbb{R}^n \to \mathbb{R} \cup \{+\infty\} \) is defined by $\mathrm{Prox}_\Phi(x) := \arg\min_{y \in \mathbb{R}^n} \left\{ \Phi(y) + \|y - x\|^2 / 2\right\}.$ For a set-valued mapping \(\mathcal{F}: \mathbb{R}^n \rightrightarrows \mathbb{R}^m\), its graph is defined as $\mathrm{gph}\, \mathcal{F} := \left\{ (x, v) \in \mathbb{R}^n \times \mathbb{R}^m \;\middle|\; v \in \mathcal{F}(x) \right\}$, and its inverse mapping \(\mathcal{F}^{-1}: \mathbb{R}^m \rightrightarrows \mathbb{R}^n\) is given by $\mathcal{F}^{-1}(v) := \left\{ x \in \mathbb{R}^n \;\middle|\; v \in \mathcal{F}(x) \right\}$. The Cartesian product of two sets \( X \) and \( Y \) is denoted by \( X \times Y \). For a vector \( x \in \mathbb{R}^n \) and a closed convex set \( D \subset \mathbb{R}^n \), the distance from \( x \) to \( D \) is defined as \(\mathrm{dist}(x, D) = \min_{y \in D} \|x - y\|\). The Euclidean projection operator onto \( D \) is denoted by \(\mathrm{Proj}_D\). For \( x \in D \), the normal cone to \( D \) at \( x \) is defined by $N_D(x) := \left\{ v \in \mathbb{R}^n \;\middle|\; \langle v, z - x \rangle \le 0 \ \text{for all } z \in D \right\}$. The Fr\'{e}chet (regular) subdifferential \(\hat{\partial} \Phi\) and the Mordukhovich (limiting) subdifferential \(\partial \Phi\) at a point \( \bar{z} \in \mathrm{dom}\, \Phi := \{z \mid \Phi(z) < \infty\} \) are defined as:
	\[
	\begin{aligned}
		\hat{\partial} \Phi(\bar{z}) &= \left\{ v \in \mathbb{R}^n \;\middle|\; \Phi(z) \geq \Phi(\bar{z}) + \langle v, z - \bar{z} \rangle + o(\|z - \bar{z}\|), \; \forall z \in \mathbb{R}^n \right\}, \\
		\partial \Phi(\bar{z}) &= \left\{ v \in \mathbb{R}^n \;\middle|\; \exists z_k \xrightarrow{\Phi} \bar{z}, \; \exists v_k \in \hat{\partial} \Phi(z_k), \; v_k \to v \right\},
	\end{aligned}
	\]
	where \( o(z) \) denotes a function satisfying \( o(z)/z\to 0 \) as \( z \downarrow 0 \), and \( z_k \xrightarrow{\Phi} \bar{z} \) indicates \( z_k \to \bar{z} \) and \( \Phi(z_k) \to \Phi(\bar{z}) \) as \( k \to \infty \). If \( \bar{z} \notin \mathrm{dom}\, \Phi \), we define \(\partial \Phi(\bar{z}) = \emptyset\).
	
	We now recall the notion of metric subregularity, which plays a central role in our convergence analysis.
	
	\begin{definition}
		Let $\mathcal{F}: \mathbb{R}^n \rightrightarrows \mathbb{R}^m$ be a set-valued mapping, $(\bar{x}, \bar{v}) \in \mathrm{gph}~\mathcal{F}$. $\mathcal{F}$ is said to be metrically subregular at $(\bar{x}, \bar{v})$ with modulus $\kappa > 0$ if there exists $\epsilon > 0$ such that
			\begin{equation*}
				\mathrm{dist}(x, \mathcal{F}^{-1}(\bar{v})) \leq \kappa\mathrm{dist}(\bar{v}, \mathcal{F}(x)) \quad \text{for all } x \in \mathbb{B}_\epsilon(\bar{x}).
			\end{equation*}
	\end{definition}
	
	When applied to the subdifferential $\partial \Phi$ of a convex function $\Phi$, metric subregularity yields the following local quadratic-growth estimate \cite[Proposition 2]{ye2021variational}.
	
	\begin{proposition}\label{prop:subgradient_error}		
		Suppose that \(\Phi: \mathbb{R}^n \to \mathbb{R} \cup \{+\infty\}\) is a proper, lower semicontinuous, convex function, and that \(\partial \Phi\) is metrically subregular at \((\bar{x}, \bar{v}) \in \mathrm{gph}\,\partial \Phi\) with modulus \(\kappa > 0\). Then there exists \(\epsilon,c > 0\) such that for all \(x \in \mathbb{B}_{\epsilon}(\bar{x})\), the following inequality holds:
		\begin{equation*}
			\Phi(x) \geq \Phi(\bar{x}) + \langle \bar{v}, x - \bar{x} \rangle + c \mathrm{dist}\left(x, (\partial \Phi)^{-1}(\bar{v})\right)^2.
		\end{equation*}
	\end{proposition}

	Finally, we recall the standard notions of convergence rates used throughout the paper.
	
	\begin{definition}
		Let $\{x^k\} \subset \mathbb{R}^n$ be a sequence converging to $\bar{x} \in \mathbb{R}^n$ as $k \to \infty$ with $x^k \neq \bar{x}$ for all $k \in \mathbb{N}$.
		\begin{enumerate}
			\item[(i)] We say that $\{x^k\}$ converges $Q$-superlinearly to $\bar{x}$ if $\lim_{k \to \infty} \|x^{k+1} - \bar{x}\|/\|x^k - \bar{x}\| = 0$.
			\item[(ii)]	We say that $\{x^k\}$ converges $R$-superlinearly to $\bar{x}$ if there exists a sequence of nonnegative scalars $\{\nu_k\}$ such that
			\[
			\|x^k - \bar{x}\| \leq \nu_k \quad \text{for all } k\in \mathbb{N},
			\]
			and $\{\nu_k\}$ converges $Q$-superlinearly to zero.
			\item[(iii)] We say that $\{x^k\}$ converges $Q$-quadratically to $\bar{x}$ if $\limsup_{k\to\infty} \left\|x^{k+1}-\bar x\right\|/\left\|x^k-\bar x\right\|^2<+\infty$.
			\item[(iv)]	We say that $\{x^k\}$ converges $R$-quadratically to $\bar{x}$ if there exists a sequence of nonnegative scalars $\{\nu_k\}$ such that
			\[
			\|x^k - \bar{x}\| \leq \nu_k \quad \text{for all } k\in \mathbb{N},
			\]
			and $\{\nu_k\}$ converges $Q$-quadratically to zero.
		\end{enumerate}
	\end{definition}

\section{A Regularized Proximal Newton Method with Local Superlinear Convergence}\label{sec3}

In this section, we develop a regularized inexact proximal Newton method for problem \eqref{simple} and establish its local superlinear convergence.

\subsection{Regularized Proximal Newton Method}\label{sec:local-newton}
We now introduce the regularized inexact proximal Newton method. Starting from an initial point $x^0\in\operatorname{dom}g$, the method generates a sequence $\{x^k\}$ by approximately solving regularized proximal Newton subproblems associated with the penalized problem~\eqref{eq:penalty}. For notational convenience, let $f:=h\circ A$, so that $\varphi=f+g$. By the standing assumptions, \(\nabla f\) is Lipschitz continuous with constant \(L_{hA}:=L_h\|A\|^2\).

Given the current iterate \(x^k\), a penalty parameter \(\sigma_k>0\), a regularization parameter \(\alpha_k>0\), and a symmetric matrix \(B_k\), we consider the following proximal Newton subproblem:
\begin{equation}\label{eq:q_model_unit}
	\min_{x\in\mathbb{R}^n}\quad
	q_k(x):=
	\sigma_kF(x^k)+f(x^k)
	+\bigl(\sigma_k\nabla F(x^k)+\nabla f(x^k)\bigr)^\top(x-x^k)
	+\frac12(x-x^k)^\top H_k(x-x^k)
	+g(x),
\end{equation}
where $H_k:=B_k+\sigma_k\alpha_k I$. Here, \(I\) denotes the identity matrix, and \(B_k\) is an approximation to the Hessian of the smooth part of \(\psi_{\sigma_k}\) at \(x^k\), that is, $B_k\approx \nabla^2 f(x^k)+\sigma_k\nabla^2F(x^k)$. We impose the conditions \(B_k\succeq0\) and \(\left\|B_k\right\|\le M_B\) on the sequence of Hessian approximations, where \(M_B>0\) is a constant. The exact Hessian \(B_k=\nabla^2 f(x^k)+\sigma_k\nabla^2F(x^k)\) is a particular choice whenever the resulting sequence is uniformly bounded.

The regularization term \(\sigma_k\alpha_k I\) gives $H_k\succeq \sigma_k\alpha_k I\succ0$, so that \(q_k\) is proper, closed, and strongly convex and hence admits a unique minimizer. Moreover, when \(\sigma_k\) and \(\alpha_k\) are chosen such that \(\sigma_k\le\bar\sigma\) and \(\alpha_k\le\bar\alpha\) for some fixed constants \(\bar\sigma,\bar\alpha>0\), the boundedness of \(B_k\) yields $\left\|H_k\right\|\le M_H:=M_B+\bar\sigma\bar\alpha$.

Our method allows for inexact solutions of the proximal Newton subproblem \eqref{eq:q_model_unit}. To this end, define the proximal residual mapping by
\begin{equation}\label{eq:newton_residual}
	r_k(x):=x-\operatorname{Prox}_g\bigl(x-[\nabla f(x^k)+\sigma_k\nabla F(x^k)+H_k(x-x^k)]\bigr).
\end{equation}
We have \(r_k(x)=0\) if and only if \(x\) is the unique minimizer of \(q_k\).

Given an accuracy parameter \(\tau_k>0\), we find an approximate model point \(\widetilde x^{k+1}\in\operatorname{dom}g\) satisfying the inexactness criterion
\begin{equation}\label{eq:inexact_unit}
	\left\|r_k(\widetilde x^{k+1})\right\|\le c_r\sigma_k\alpha_k\tau_k,
\end{equation}
where \(c_r>0\) is fixed. We then construct the search direction \(d^k\) by
\begin{equation}\label{eq:unit_output}
	d^k:=\operatorname{Prox}_g\Bigl(\widetilde x^{k+1}-\bigl[\nabla f(x^k)+\sigma_k\nabla F(x^k)+H_k(\widetilde x^{k+1}-x^k)\bigr]\Bigr)-x^k.
\end{equation}
The next iterate is updated by
\begin{equation}\label{eq:unit_update}
	x^{k+1}:=x^k+d^k.
\end{equation}
By the definition of \(r_k\), the updated point satisfies $x^{k+1}=\widetilde x^{k+1}-r_k(\widetilde x^{k+1})$. The above update is well defined, since the exact minimizer of \(q_k\) has zero model residual, ensuring the existence of a point satisfying \eqref{eq:inexact_unit}.

\subsection{Structural Properties and Local Assumptions}\label{sec3.3}
Before establishing local superlinear convergence, we introduce additional assumptions and some preliminary properties needed for the analysis.

\begin{assumption}\label{asup:common-image}
	The function \(h\) is strictly convex.
\end{assumption}
\begin{remark}\label{rem:fixed-kernel}
	Assumption~\ref{asup:common-image} is a sufficient condition for the single-valuedness of the inverse mapping $(\nabla h)^{-1}$, which is the key property needed to establish the following results. Indeed, for any $z,s\in X_\varphi^*$, the convexity and lower-level optimality conditions imply that $\nabla h(Az)=\nabla h(As)$. Hence, the single-valuedness of $(\nabla h)^{-1}$ yields
	\[
	Az=As
	\qquad
	\forall\, z,s\in X_\varphi^*.
	\]
	Consequently,
	\[
	X_\varphi^*-X_\varphi^*\subseteq\Null A,
	\qquad
	\nabla^2 f(u)(z-s)=0
	\quad
	\forall\,u\in\mathbb{R}^n,\ z,s\in X_\varphi^*,
	\]
	where the second relation follows from $\nabla^2 f(u)=A^\top\nabla^2 h(Au)A$. These structural properties will be used in the subsequent local convergence analysis.
\end{remark}

Define the proximal residual associated with the lower-level problem by
\[
R_L(x):=\|x-\operatorname{Prox}_g\bigl(x-\nabla f(x)\bigr)\|.
\]
By the first-order optimality condition, we have $R_L(x)=0$ if and only if $x\in X_\varphi^*$. The following proposition gives the basic properties of the lower-level solution set and its residual, together with additional structural properties under Assumption~\ref{asup:common-image}.

\begin{proposition}\label{prop:simple-blo-geometry}
	Suppose that the standing assumptions hold. Then \(X_\varphi^*\) is nonempty, closed, and convex, \(R_L^{-1}(0)=X_\varphi^*\), the residual \(R_L\) is continuous on \(\mathbb{R}^n\), and
	\begin{equation}\label{eq:residual_by_subgradient}
		R_L(x)\le\dist(0,\partial\varphi(x)).
	\end{equation}
	If, in addition, Assumption~\ref{asup:common-image} holds, then
	\[
	\left\|A(x-s)\right\|\le\left\| A\right\|\,\dist(x,X_\varphi^*)
	\qquad\ \forall\ x\in\mathbb{R}^n,\ s\in X_\varphi^*,
	\]
	and the gradient \(\nabla f\) is constant on \(X_\varphi^*\).
\end{proposition}
\begin{proof}
	The lower-level solution set is nonempty, closed, and convex. Continuity of $R_L$ and its zero set follow from proximal optimality and nonexpansiveness. If $v\in\partial\varphi(x)$, then $x=\operatorname{Prox}_g(x-\nabla f(x)+v)$; nonexpansiveness gives $R_L(x)\le\|v\|$, proving \eqref{eq:residual_by_subgradient}; see also \cite[Theorem~3.5]{drusvyatskiy2018error}.

	Suppose now that Assumption~\ref{asup:common-image} also holds. The image estimate follows by taking $z=\operatorname{Proj}_{X_\varphi^*}(x)$ and using $A(x-s)=A(x-z)$, since $Az=As$ by the common-image property established in Remark~\ref{rem:fixed-kernel}; see also \cite[Proof of Proposition~2.6]{xiong2026computational} for the corresponding distance estimate under a linear mapping. For \(z,s\in X_\varphi^*\), the common image gives \(Az=As\), so \(f(z)=h(Az)=h(As)=f(s)\) and \(\nabla f(z)=A^\top\nabla h(Az)=A^\top\nabla h(As)=\nabla f(s)\); in particular, \(\nabla f\) is constant on \(X_\varphi^*\), and \(\left\langle \nabla f(s),z-s\right\rangle=\left\langle \nabla h(As),A(z-s)\right\rangle=0\).
\end{proof}

We next impose the following metric subregularity assumption at a given point \(\bar x\in X_\varphi^*\).

\begin{assumption}\label{asup:lower_regular}
	The mapping \(\partial\varphi\) is metrically subregular at \((\bar x,0)\): there exist \(\epsilon_m,\kappa_m>0\) such that
	\begin{equation*}
		\dist(x,X_\varphi^*)\le\kappa_m\dist(0,\partial\varphi(x))
		\qquad\ \forall\ x\in\mathbb B_{\epsilon_m}(\bar x).
	\end{equation*}
\end{assumption}

Denote the optimal value of the lower-level problem by \(\varphi^*:=\min_{x\in\mathbb{R}^n}\varphi(x)\). The following standard consequences of metric subregularity relate the lower-level residual to the distance from the solution set and the quadratic growth of the objective; see \cite[Proposition 2.2]{mordukhovich2023globally} for proximal residual estimates and \cite[Proposition 2]{ye2021variational} for local quadratic growth.

\begin{proposition}\label{prop:lower-error-bound}
	Suppose that the standing assumptions hold. Fix \(\bar x\in X_\varphi^*\) and suppose that Assumption~\ref{asup:lower_regular} holds at \(\bar x\). Then there exist \(\epsilon_L,c_m,c_L>0\) such that the following estimates hold:
	\begin{enumerate}[label=(\roman*)]
		\item The lower residual and the distance to the lower-level solution set satisfy
		\begin{equation*}
			\dist(x,X_\varphi^*)\le c_mR_L(x)
			\qquad\ \forall\ x\in\mathbb B_{\epsilon_L}(\bar x).
		\end{equation*}
		\item The lower-level objective satisfies the local quadratic-growth estimate
		\begin{equation*}
			\varphi(x)-\varphi^*\ge c_L\dist(x,X_\varphi^*)^2
			\qquad\ \forall\ x\in\mathbb B_{\epsilon_L}(\bar x).
		\end{equation*}
	\end{enumerate}
\end{proposition}

We now turn to the optimality properties of the simple bilevel problem \eqref{simple}. Denote the solution set of problem~\eqref{simple} by \(S^*:=\mathop{\arg\min}_{x\in X_\varphi^*}F(x)\). Since \(X_\varphi^*\) is closed and convex and \(F\) is continuous and convex, \(S^*\) is closed and convex. For the subsequent local convergence analysis, we restrict the reference point to \(\bar x\in S^*\subseteq X_\varphi^*\) and denote the optimal value of problem~\eqref{simple} by \(F^*:=\min_{x\in X_\varphi^*}F(x)\). We next impose a local quadratic-growth condition for the simple bilevel problem and local Lipschitz continuity of the Hessians.

\begin{assumption}\label{asup_quadratic_growth}
	There exist \(\epsilon_q,\kappa_q>0\) such that
	\[
	F(z)-F^*\ge\kappa_q\dist(z,S^*)^2
	\qquad\ \forall\ z\in X_\varphi^*\cap\mathbb B_{\epsilon_q}(\bar x).
	\]
\end{assumption}

\begin{assumption}\label{asup_hessian_lipschitz}
	There exist \(\epsilon_{F,2},\epsilon_{h,2}>0\) and \(L_{F,2},L_{h,2}>0\) such that the Hessians \(\nabla^2F\) and \(\nabla^2h\) are Lipschitz continuous on \(\mathbb B_{\epsilon_{F,2}}(\bar x)\) and \(\mathbb B_{\epsilon_{h,2}}(A\bar x)\), respectively, with Lipschitz constants \(L_{F,2}\) and \(L_{h,2}\).
\end{assumption}

\subsection{Local Superlinear Convergence}\label{sec3.2}

In this subsection, we establish local superlinear convergence of the regularized inexact proximal Newton method introduced in Section~\ref{sec:local-newton} under the parameter selection strategy specified below.

To this end, define the exact-Hessian proximal residual at \(x^{k+1}\) by
\begin{equation}\label{eq:exact_hessian_residual}
	\widetilde r^k:=x^{k+1}
	-\operatorname{Prox}_g\!\Bigl(x^{k+1}-\bigl[\nabla f(x^k)+\sigma_k\nabla F(x^k)+\bigl(\nabla^2 f(x^k)+\sigma_k\nabla^2F(x^k)
	+\sigma_k\alpha_k I\bigr)d^k\bigr]\Bigr).
\end{equation}
For each Newton output, define
\begin{equation}\label{eq:local_trial_measure}
	G_{k+1}:=\frac{\left\|\widetilde r^k\right\|}{\sigma_k}
	+\left\|d^k\right\|^2+\alpha_k\left\|d^k\right\|
	+\sqrt{R_L(x^{k+1})}.
\end{equation}

Given \(x^0\in\operatorname{dom}g\), choose \(G_0\ge\sqrt{R_L(x^0)}\) and \(\tau_0>0\), and set
\begin{equation}\label{eq:omega_def}
	\omega_k:=G_k+\tau_k.
\end{equation}
The penalty, regularization, and accuracy parameters are selected according to
\begin{equation}\label{eq:parameter-rule}
	\sigma_k:=\min\{\bar\sigma,c_\sigma\omega_k^{2(1+\rho)}\},\quad
	\alpha_k:=\min\{\bar\alpha,c_\alpha\omega_k^\rho\},\quad
	\tau_{k+1}:=\min\{\gamma\tau_k,c_\tau\omega_k^{1+\rho}\},
\end{equation}
where \(\bar\sigma,\bar\alpha>0\), \(0<\rho\le1\), \(0<\gamma<1\), and \(c_\sigma,c_\alpha,c_\tau>0\) are fixed constants.

The local convergence analysis proceeds in three stages. First, exploiting the lower-level structure and the local quadratic-growth conditions, we derive a distance estimate for the perturbed Newton model in Proposition~\ref{prop:measure-distance}. A key step is the lower-level feasibility estimate $\dist(x^k,X_\varphi^*)=O(\omega_k^2)$ in Lemma~\ref{lem:lower-distance}, which, together with the parameter rule \eqref{eq:parameter-rule}, allows us to control the terms involving the vanishing penalty parameter. This leads to the estimate $\dist(x^{k+1},S^*)\le c_S\omega_{k+1}$ in Lemma~\ref{lem:measure-distance-iteration}.

Next, we establish bounds on the Newton step and the lower-level residual in Lemmas~\ref{lem:step_estimate} and \ref{lem:lower-stationarity}, respectively. Combined with the exact-Hessian residual estimate, these bounds yield the superlinear recursion $\omega_{k+1}\le c_\omega\omega_k^{1+\rho}$. Finally, we show that the iterates remain in a neighborhood where the local estimates apply, and establish local superlinear convergence in Theorem~\ref{thm:local-recursion}.

We begin with a lower-level model estimate. Since \(f=h\circ A\), Assumption~\ref{asup_hessian_lipschitz} implies that \(\nabla^2 f\) is locally Lipschitz continuous around \(\bar x\). Let \(\epsilon_{f,2}>0\) and \(L_{f,2}>0\) denote a corresponding neighborhood radius and Lipschitz constant, respectively.

\begin{lemma}\label{lem:lower-model-growth}
	Suppose that the standing assumptions and Assumption~\ref{asup:common-image} hold, and that Assumptions~\ref{asup:lower_regular} and \ref{asup_hessian_lipschitz} hold at the reference solution \(\bar x\in S^*\). Then there exists \(\epsilon_f>0\) such that, for all \(x,\widehat x\in\mathbb B_{\epsilon_f}(\bar x)\) with \(\widehat x\in\operatorname{dom}g\), and every \(s\in X_\varphi^*\cap\mathbb B_{\epsilon_f}(\bar x)\),
	\begin{equation}\label{eq:lower_model_growth}
		\begin{aligned}
			&\quad\ g(\widehat x)-g(s)+\left\langle \nabla f(x),\widehat x-s\right\rangle
			+\frac12\Bigl[
			\left\langle \widehat x-x,\nabla^2 f(x)(\widehat x-x)\right\rangle
			-\left\langle s-x,\nabla^2 f(x)(s-x)\right\rangle
			\Bigr]
			\\
			&\ge
			\frac{c_L}{2}\dist(\widehat x,X_\varphi^*)^2
			-\frac{L_{f,2}}{2}
			\dist(x,X_\varphi^*)^2\dist(\widehat x,X_\varphi^*).
		\end{aligned}
	\end{equation}
	where \(c_L>0\) is the constant in Proposition~\ref{prop:lower-error-bound}(ii), and \(L_{f,2}>0\) is the local Lipschitz constant of \(\nabla^2 f\) specified above.
\end{lemma}
\begin{proof} 
	Let \(\epsilon_L,c_L>0\) be the fixed radius and constant of Proposition~\ref{prop:lower-error-bound}(ii), and let \(\epsilon_{f,2}>0\) and \(L_{f,2}\ge0\) be the fixed radius and modulus following Assumption~\ref{asup_hessian_lipschitz}. Set
	\[ 
	\epsilon_f:=
	\min\left\{1,\epsilon_L,\epsilon_{f,2},\frac{3c_L}{7\max\{L_{f,2},1\}}\right\}. 
	\] 
	Set \(\widehat x_\varphi:=\operatorname{Proj}_{X_\varphi^*}(\widehat x)\) and \(x_\varphi:=\operatorname{Proj}_{X_\varphi^*}(x)\). Since \(\bar x=\operatorname{Proj}_{X_\varphi^*}(\bar x)\), nonexpansiveness of the projection gives \(\left\|\widehat x_\varphi-\bar x\right\| \le\left\|\widehat x-\bar x\right\|\le\epsilon_f\) and \(\left\|x_\varphi-\bar x\right\| \le\left\|x-\bar x\right\|\le\epsilon_f\). Thus the local estimates of Proposition~\ref{prop:lower-error-bound}(ii) and Assumption~\ref{asup_hessian_lipschitz} apply to these points.
	
	We start from the left-hand side of \eqref{eq:lower_model_growth}. We obtain
	\begin{equation}\label{eq:lower-model-expansion}
		\begin{aligned} 
			&\quad\ g(\widehat x)-g(s)+\left\langle \nabla f(x),\widehat x-s\right\rangle 
			+\frac12\Bigl[ 
			\left\langle \widehat x-x,\nabla^2 f(x)(\widehat x-x)\right\rangle 
			-\left\langle s-x,\nabla^2 f(x)(s-x)\right\rangle 
			\Bigr] 
			\\ 
			&= 
			g(\widehat x)-g(s) 
			+\left\langle \nabla f(x)+\nabla^2 f(x)(s-x),\widehat x-s\right\rangle 
			+\frac12\left\langle \widehat x-s,\nabla^2 f(x)(\widehat x-s)\right\rangle 
			\\ 
			&= 
			g(\widehat x)-g(s) 
			+\left\langle \nabla f(s),\widehat x-s\right\rangle 
			+\left\langle \nabla f(x)+\nabla^2 f(x)(s-x)-\nabla f(s), 
			\widehat x-s\right\rangle+\frac12\left\langle \widehat x-s,\nabla^2 f(x)(\widehat x-s)\right\rangle. 
		\end{aligned} 
	\end{equation}
	
	We now simplify the three terms on the right-hand side. First split \(\widehat x-s=(\widehat x-\widehat x_\varphi) +(\widehat x_\varphi-s)\). Since \(s,\widehat x_\varphi\in X_\varphi^*\), Proposition~\ref{prop:simple-blo-geometry} and Remark~\ref{rem:fixed-kernel} give $\nabla f(\widehat x_\varphi)=\nabla f(s)$, $\left\langle \nabla f(s),\widehat x_\varphi-s\right\rangle=0$.	Therefore,
	\begin{equation}\label{eq:lower-model-value}
		g(\widehat x)-g(s)+\left\langle \nabla f(s),\widehat x-s\right\rangle= 
		\varphi(\widehat x)-\varphi^* 
		-\Bigl[ 
		f(\widehat x)-f(\widehat x_\varphi) 
		-\left\langle \nabla f(\widehat x_\varphi), 
		\widehat x-\widehat x_\varphi\right\rangle 
		\Bigr]. 
	\end{equation}
	
	Since \(\nabla^2 f(x)(\widehat x_\varphi-s)=0\) by Remark~\ref{rem:fixed-kernel}, we obtain
	\begin{equation}\label{eq:lower-model-hessian}
		\left\langle \widehat x-s,
		\nabla^2 f(x)(\widehat x-s)\right\rangle
		=
		\left\langle \widehat x-\widehat x_\varphi,
		\nabla^2 f(x)(\widehat x-\widehat x_\varphi)\right\rangle.
	\end{equation}
	
	Moreover, \(\nabla f(x)+\nabla^2 f(x)(s-x)-\nabla f(s) \in\operatorname{range}A^\top\), whereas \(\widehat x_\varphi-s\in\operatorname{Null}A\) by Remark~\ref{rem:fixed-kernel}. Their orthogonality yields
	\begin{equation}\label{eq:lower-model-gradient}
		\left\langle 
			\nabla f(x)+\nabla^2 f(x)(s-x)-\nabla f(s)
		,\widehat x-s\right\rangle=
		\left\langle 
			\nabla f(x)+\nabla^2 f(x)(s-x)-\nabla f(s)
		,\widehat x-\widehat x_\varphi\right\rangle.
	\end{equation}

	Substituting \eqref{eq:lower-model-value}, \eqref{eq:lower-model-hessian}, and \eqref{eq:lower-model-gradient} into \eqref{eq:lower-model-expansion} gives
	\begin{equation}\label{eq:lower-model-prebound}
		\begin{aligned} 
			&\quad\ g(\widehat x)-g(s)+\left\langle \nabla f(x),\widehat x-s\right\rangle 
			+\frac12\Bigl[ 
			\left\langle \widehat x-x,\nabla^2 f(x)(\widehat x-x)\right\rangle 
			-\left\langle s-x,\nabla^2 f(x)(s-x)\right\rangle 
			\Bigr] 
			\\ 
			&= 
			\varphi(\widehat x)-\varphi^* 
			-\Bigl[ 
			f(\widehat x)-f(\widehat x_\varphi) 
			-\left\langle \nabla f(\widehat x_\varphi), 
			\widehat x-\widehat x_\varphi\right\rangle 
			-\frac12 
			\left\langle \widehat x-\widehat x_\varphi, 
			\nabla^2 f(x)(\widehat x-\widehat x_\varphi)\right\rangle 
			\Bigr] 
			\\ 
			&\quad 
			+\left\langle  
				\nabla f(x)+\nabla^2 f(x)(s-x)-\nabla f(s) 
			,\widehat x-\widehat x_\varphi\right\rangle\\
			&\ge 
			c_L\dist(\widehat x,X_\varphi^*)^2 -\left\| 
				f(\widehat x)-f(\widehat x_\varphi) 
				-\left\langle \nabla f(\widehat x_\varphi), 
				\widehat x-\widehat x_\varphi\right\rangle 
				-\frac12 
				\left\langle \widehat x-\widehat x_\varphi, 
				\nabla^2 f(x)(\widehat x-\widehat x_\varphi)\right\rangle 
			\right\| 
			\\ 
			&\quad 
			-\left\| 
				\nabla f(x)+\nabla^2 f(x)(s-x)-\nabla f(s) 
			\right\| 
			\dist(\widehat x,X_\varphi^*),
		\end{aligned} 
	\end{equation}
	where the inequality follows from Proposition~\ref{prop:lower-error-bound}(ii) and Cauchy--Schwarz.
	
	For the first term, adding and subtracting the quadratic term with \(\nabla^2 f(\widehat x_\varphi)\) and applying the second-order Taylor remainder estimate, we obtain
	\begin{equation}\label{eq:lower-model-function-remainder}
		\begin{aligned} 
			&\quad\ \left\| 
				f(\widehat x)-f(\widehat x_\varphi) 
				-\left\langle \nabla f(\widehat x_\varphi), 
				\widehat x-\widehat x_\varphi\right\rangle 
				-\frac12 
				\left\langle \widehat x-\widehat x_\varphi, 
				\nabla^2 f(x)(\widehat x-\widehat x_\varphi)\right\rangle 
			\right\| 
			\\ 
			&\le 
			\frac{L_{f,2}}{6} 
			\left\|\widehat x-\widehat x_\varphi\right\|^3 
			+\frac12 
			\left\|\nabla^2 f(\widehat x_\varphi)-\nabla^2 f(x)\right\| 
			\left\|\widehat x-\widehat x_\varphi\right\|^2 
			\\ 
			&\le 
			\frac{L_{f,2}}{6} 
			\dist(\widehat x,X_\varphi^*)^3 
			+\frac{L_{f,2}}{2} 
			\left\|x-\widehat x_\varphi\right\| 
			\dist(\widehat x,X_\varphi^*)^2 
			\\ 
			&= 
			\left[ 
			\frac{L_{f,2}}{6}\dist(\widehat x,X_\varphi^*) 
			+\frac{L_{f,2}}{2}\left\|x-\widehat x_\varphi\right\| 
			\right] 
			\dist(\widehat x,X_\varphi^*)^2. 
		\end{aligned} 
	\end{equation}
	The second inequality uses the Lipschitz continuity of \(\nabla^2 f\) and \(\left\|\widehat x-\widehat x_\varphi\right\| =\dist(\widehat x,X_\varphi^*)\).
	
	For the remaining term, Proposition~\ref{prop:simple-blo-geometry} and Remark~\ref{rem:fixed-kernel} give $\nabla f(s)=\nabla f(x_\varphi)$, $\nabla^2 f(x)(s-x_\varphi)=0$. Hence
	\begin{equation}\label{eq:lower-model-gradient-remainder}
		\begin{aligned} 
			&\quad\ \|\nabla f(x)+\nabla^2 f(x)(s-x)-\nabla f(s)\|
			\\ 
			&= 
			\|\nabla f(x) 
			+\nabla^2 f(x)(x_\varphi-x) 
			-\nabla f(x_\varphi)\|\\ 
			&\le 
			\frac{L_{f,2}}{2}\left\|x-x_\varphi\right\|^2 
			= 
			\frac{L_{f,2}}{2}\dist(x,X_\varphi^*)^2,
		\end{aligned} 
	\end{equation}
	where the inequality follows from the gradient Taylor remainder between \(x\) and \(x_\varphi\), and the last equality follows from \(x_\varphi=\operatorname{Proj}_{X_\varphi^*}(x)\).
	
	Substituting \eqref{eq:lower-model-function-remainder} and \eqref{eq:lower-model-gradient-remainder} into \eqref{eq:lower-model-prebound} yields
	\[ 
	\begin{aligned} 
		&\quad\ g(\widehat x)-g(s)+\left\langle \nabla f(x),\widehat x-s\right\rangle 
		+\frac12\Bigl[ 
		\left\langle \widehat x-x,\nabla^2 f(x)(\widehat x-x)\right\rangle 
		-\left\langle s-x,\nabla^2 f(x)(s-x)\right\rangle 
		\Bigr] 
		\\ 
		&\ge 
		\left[ 
		c_L 
		-\frac{L_{f,2}}{6}\dist(\widehat x,X_\varphi^*) 
		-\frac{L_{f,2}}{2}\left\|x-\widehat x_\varphi\right\| 
		\right] 
		\dist(\widehat x,X_\varphi^*)^2-\frac{L_{f,2}}{2} 
		\dist(x,X_\varphi^*)^2 
		\dist(\widehat x,X_\varphi^*). 
	\end{aligned} 
	\] 
	
	Finally, since \(\dist(\widehat x,X_\varphi^*) \le\left\|\widehat x-\bar x\right\|\le\epsilon_f\) and $\left\|x-\widehat x_\varphi\right\| \le \left\|x-\bar x\right\|+\left\|\widehat x_\varphi-\bar x\right\| \le 2\epsilon_f$,
	\[ 
	c_L 
	-\frac{L_{f,2}}{6}\dist(\widehat x,X_\varphi^*) 
	-\frac{L_{f,2}}{2}\left\|x-\widehat x_\varphi\right\| 
	\ge 
	c_L-\frac{7L_{f,2}}{6}\epsilon_f 
	\ge 
	\frac{c_L}{2}, 
	\] 
	which proves \eqref{eq:lower_model_growth}.
\end{proof}

We next highlight a difficulty arising in the local convergence analysis. A direct error-bound analysis of the penalized optimality mapping
\[
T_\sigma(x):=\nabla f(x)+\sigma\nabla F(x)+\partial g(x)
\]
would require an error bound uniform as $\sigma\downarrow0$. However, metric subregularity of $T_0=\partial\varphi$ does not guarantee such a uniform bound, since metric subregularity is not generally stable under perturbations \cite{zheng2021perturbation,dontchev2004regularity,benko2018new,dontchev2020radius}. Consequently, a direct distance analysis based on the penalized solution sets cannot be justified under our assumptions. We therefore estimate the Newton output directly relative to the bilevel solution set $S^*$.

\begin{proposition}\label{prop:measure-distance}
	Suppose that the standing assumptions and Assumption~\ref{asup:common-image} hold. Fix \(\bar x\in S^*\) and suppose that Assumptions~\ref{asup:lower_regular}, \ref{asup_quadratic_growth}, and \ref{asup_hessian_lipschitz} hold at \(\bar x\). Then there exist \(\epsilon_G,c_G>0\) such that, for any \(x,\widehat x\in\mathbb B_{\epsilon_G}(\bar x)\cap\operatorname{dom}g\), \(\sigma,\alpha>0\), and \(v\in\mathbb R^n\) satisfying
	\begin{equation}\label{eq:exact-model-certificate}
		v-\nabla f(x)-\sigma\nabla F(x)
		-\bigl[\nabla^2 f(x)+\sigma\nabla^2F(x)+\sigma\alpha I\bigr]
		(\widehat x-x)
		\in\partial g(\widehat x),
	\end{equation}
	the following estimate holds
	\begin{equation}\label{eq:accepted-distance}
		\dist(\widehat x,S^*)
		\le
		c_G\left[
		\frac{\left\|v\right\|}{\sigma}
		+\left\|\widehat x-x\right\|^2
		+\alpha\left\|\widehat x-x\right\|
		+\left(1+\frac{\dist(x,X_\varphi^*)}{\sqrt{\sigma}}\right)
		\sqrt{R_L(\widehat x)}
		\right].
	\end{equation}
\end{proposition}

\begin{proof}
	Let \(\epsilon_L,\epsilon_q,\epsilon_f,\epsilon_{F,2}\) be the fixed radii of Proposition~\ref{prop:lower-error-bound}, Assumption~\ref{asup_quadratic_growth}, Lemma~\ref{lem:lower-model-growth}, and Assumption~\ref{asup_hessian_lipschitz}, respectively, and let \(c_m,L_{F,2}\) be the corresponding constants. Since \(R_L\) is continuous and \(R_L(\bar x)=0\), fix \(\epsilon_b>0\) such that \(R_L(u)\le1\) on \(\mathbb B_{\epsilon_b}(\bar x)\). Set
	\[
	\epsilon_G:=
	\min\{1,\epsilon_L,\epsilon_q,\epsilon_f,\epsilon_{F,2},\epsilon_b\}.
	\]
	By the continuity of \(\nabla F\), we may also fix \(M_F>0\) such that
	\[
	\|\nabla F(u)\|\le M_F
	\qquad \forall\,u\in\mathbb B_{\epsilon_G}(\bar x).
	\]
	Set \(p:=\operatorname{Proj}_{X_\varphi^*}(\widehat x)\) and \(s:=\operatorname{Proj}_{S^*}(p)\). Since the projections onto \(X_\varphi^*\) and \(S^*\) fix \(\bar x\), nonexpansiveness gives \(\left\|p-\bar x\right\|\le\left\|\widehat x-\bar x\right\|\) and \(\left\|s-\bar x\right\|\le\left\|p-\bar x\right\|\), so \(p,s\) remain in \(\mathbb B_{\epsilon_G}(\bar x)\).
		
	We first separate the distance to \(S^*\) into the lower-level feasibility error and the upper-level optimality error on \(X_\varphi^*\). By the triangle inequality and the projection identity \(\left\|\widehat x-p\right\|=\dist(\widehat x,X_\varphi^*)\), we obtain
	\begin{equation}\label{eq:prop3-distance-split}
		\dist(\widehat x,S^*)
		\le
		\left\|\widehat x-p\right\|+\left\|p-s\right\|
		=
		\dist(\widehat x,X_\varphi^*)+\left\|p-s\right\|.
	\end{equation}
	The first term is controlled by Proposition~\ref{prop:lower-error-bound}(i):
	\begin{equation}\label{eq:prop3-normal-bound}
		\dist(\widehat x,X_\varphi^*)\le c_mR_L(\widehat x)
		\le c_m\sqrt{R_L(\widehat x)},
	\end{equation}
	where the second inequality uses \(R_L(\widehat x)\le1\). We next estimate \(\left\|p-s\right\|\) using the upper-level quadratic-growth condition in Assumption~\ref{asup_quadratic_growth} and the lower-level model estimate in Lemma~\ref{lem:lower-model-growth}.

	Assumption~\ref{asup_quadratic_growth} at \(p\in X_\varphi^*\), together with \(\dist(p,S^*)=\left\|p-s\right\|\) and \(F(s)=F^*\), gives
	\begin{equation}\label{eq:prop3-tangential-start}
		\begin{aligned}
			\kappa_q\left\|p-s\right\|^2
			&\le F(p)-F(s) \\
			&= F(\widehat x)-F(s)+F(p)-F(\widehat x)\\
			&\le F(\widehat x)-F(s)
			+\left\langle \nabla F(p),p-\widehat x\right\rangle\\
			&\le F(\widehat x)-F(s)
			+M_F\dist(\widehat x,X_\varphi^*).
		\end{aligned}
	\end{equation}
	The second inequality follows from the convexity of \(F\), and the third from the Cauchy--Schwarz inequality, the local bound \(\left\|\nabla F(p)\right\|\le M_F\), and the projection identity \(\left\|p-\widehat x\right\|=\dist(\widehat x,X_\varphi^*)\). Hence it remains to bound \(F(\widehat x)-F(s)\) from above.
	
	By the convexity of \(F\), adding and subtracting the first-order quadratic-model terms at \(x\) gives
	\begin{equation}\label{eq:prop3-upper-gap-decomposition}
		\begin{aligned}
			&\quad\ \sigma\bigl[F(\widehat x)-F(s)\bigr]\\
			&\le
			\sigma\left\langle \nabla F(\widehat x),\widehat x-s\right\rangle\\
			&=
			\sigma\left\langle 
				\nabla F(\widehat x)-\nabla F(x)
				-\nabla^2F(x)(\widehat x-x),
			\widehat x-s\right\rangle+
			\sigma\left\langle 
				\nabla F(x)+\nabla^2F(x)(\widehat x-x),
			\widehat x-s\right\rangle.
		\end{aligned}
	\end{equation}
	The first inner product in \eqref{eq:prop3-upper-gap-decomposition} is a second-order gradient remainder of \(F\). It remains to control the second inner product: applying the subgradient inequality at \(s\) to \eqref{eq:exact-model-certificate} gives
	\[
	\begin{aligned}
		&\quad\ g(\widehat x)-g(s)\\
		&\le
		\left\langle v-\nabla f(x)-\sigma\nabla F(x)
			-\bigl[\nabla^2 f(x)+\sigma\nabla^2F(x)+\sigma\alpha I\bigr]
			(\widehat x-x),\widehat x-s\right\rangle\\
		&=
		\left\langle v,\widehat x-s\right\rangle
		-\left\langle \nabla f(x)+\nabla^2 f(x)(\widehat x-x),
		\widehat x-s\right\rangle-
		\sigma\left\langle \nabla F(x)+\nabla^2F(x)(\widehat x-x),
		\widehat x-s\right\rangle
		-\sigma\alpha\left\langle \widehat x-x,\widehat x-s\right\rangle.
	\end{aligned}
	\]
	Rearranging this inequality yields
	\[
	\begin{aligned}
		&\quad\ \sigma\left\langle \nabla F(x)+\nabla^2F(x)(\widehat x-x),
		\widehat x-s\right\rangle\\
		&\le
		\left\langle v,\widehat x-s\right\rangle
		-\sigma\alpha\left\langle \widehat x-x,\widehat x-s\right\rangle-
		\Bigl[
		g(\widehat x)-g(s)
		+\left\langle \nabla f(x)+\nabla^2 f(x)(\widehat x-x),
		\widehat x-s\right\rangle
		\Bigr].
	\end{aligned}
	\]
	Substituting this bound into \eqref{eq:prop3-upper-gap-decomposition} gives
	\begin{equation}\label{eq:prop3-proximal-combination}
		\begin{aligned}
			\sigma\bigl[F(\widehat x)-F(s)\bigr]
			\le{}&
			\left\langle v,\widehat x-s\right\rangle
			-\sigma\alpha\left\langle \widehat x-x,\widehat x-s\right\rangle+
			\sigma\left\langle 
				\nabla F(\widehat x)-\nabla F(x)
				-\nabla^2F(x)(\widehat x-x),
			\widehat x-s\right\rangle
			\\
			&-
			\Bigl[
			g(\widehat x)-g(s)
			+\left\langle \nabla f(x)+\nabla^2 f(x)(\widehat x-x),
			\widehat x-s\right\rangle
			\Bigr].
		\end{aligned}
	\end{equation}
	
	The lower-level term in \eqref{eq:prop3-proximal-combination} is controlled directly by Lemma~\ref{lem:lower-model-growth}. Indeed,
	\begin{equation}\label{eq:prop3-lower-model-bound}
		\begin{aligned}
			&\quad\ g(\widehat x)-g(s)
			+\left\langle \nabla f(x)+\nabla^2 f(x)(\widehat x-x),
			\widehat x-s\right\rangle
			\\
			&=
			g(\widehat x)-g(s)
			+\left\langle \nabla f(x),\widehat x-s\right\rangle+
			\frac12\Bigl[
			\left\langle \widehat x-x,\nabla^2 f(x)(\widehat x-x)\right\rangle
			-\left\langle s-x,\nabla^2 f(x)(s-x)\right\rangle
			\Bigr]
			+
			\frac12
			\left\langle \widehat x-s,\nabla^2 f(x)(\widehat x-s)\right\rangle
			\\
			&\ge
			g(\widehat x)-g(s)
			+\left\langle \nabla f(x),\widehat x-s\right\rangle+
			\frac12\Bigl[
			\left\langle \widehat x-x,\nabla^2 f(x)(\widehat x-x)\right\rangle
			-\left\langle s-x,\nabla^2 f(x)(s-x)\right\rangle
			\Bigr]
			\\
			&\ge
			\frac{c_L}{2}\dist(\widehat x,X_\varphi^*)^2
			-\frac{L_{f,2}}{2}
			\dist(x,X_\varphi^*)^2
			\dist(\widehat x,X_\varphi^*).
		\end{aligned}
	\end{equation}
	The first equality follows from the quadratic identity for \(\nabla^2 f(x)\). The first inequality uses \(\nabla^2 f(x)=A^\top\nabla^2h(Ax)A\succeq0\), while the second inequality is exactly Lemma~\ref{lem:lower-model-growth}.
	
	Substituting \eqref{eq:prop3-lower-model-bound} into \eqref{eq:prop3-proximal-combination} yields
	\begin{equation}\label{eq:prop3-combined-lower}
		\begin{aligned}
			\sigma\bigl[F(\widehat x)-F(s)\bigr]
			\le{}&
			\left\langle v,\widehat x-s\right\rangle
			-\sigma\alpha\left\langle \widehat x-x,\widehat x-s\right\rangle+
			\sigma\left\langle 
				\nabla F(\widehat x)-\nabla F(x)
				-\nabla^2F(x)(\widehat x-x),
			\widehat x-s\right\rangle
			\\
			&+
			\frac{L_{f,2}}{2}
			\dist(x,X_\varphi^*)^2
			\dist(\widehat x,X_\varphi^*)
			-\frac{c_L}{2}\dist(\widehat x,X_\varphi^*)^2.
		\end{aligned}
	\end{equation}
	
	By the Hessian Lipschitz continuity of \(F\),
	\[
	\left\|
		\nabla F(\widehat x)-\nabla F(x)
		-\nabla^2F(x)(\widehat x-x)\right\|
	\le
	\frac{L_{F,2}}{2}\left\|\widehat x-x\right\|^2.
	\]
	Applying the Cauchy--Schwarz inequality to the three inner products in \eqref{eq:prop3-combined-lower}, using this remainder bound, and dividing by \(\sigma>0\), we obtain
	\begin{equation}\label{eq:prop3-upper-model}
		F(\widehat x)-F(s)
		\le
		\left[
		\frac{\left\|v\right\|}{\sigma}
		+\alpha\left\|\widehat x-x\right\|
		+\frac{L_{F,2}}{2}\left\|\widehat x-x\right\|^2
		\right]\left\|\widehat x-s\right\|+
		\frac{L_{f,2}}{2\sigma}
		\dist(x,X_\varphi^*)^2
		\dist(\widehat x,X_\varphi^*)
		-\frac{c_L}{2\sigma}
		\dist(\widehat x,X_\varphi^*)^2.
	\end{equation}
	
	Substituting \eqref{eq:prop3-upper-model} into \eqref{eq:prop3-tangential-start}, using $\left\|\widehat x-s\right\| \le \dist(\widehat x,X_\varphi^*)+\left\|p-s\right\|$, and discarding the nonpositive quadratic term gives
	\begin{equation}\label{eq:prop3-young-start}
		\begin{aligned}
			\kappa_q\left\|p-s\right\|^2
			\le{}&
			\left[
			\frac{\left\|v\right\|}{\sigma}
			+\alpha\left\|\widehat x-x\right\|
			+\frac{L_{F,2}}{2}\left\|\widehat x-x\right\|^2
			\right]\left\|p-s\right\|+
			\left[
			\frac{\left\|v\right\|}{\sigma}
			+\alpha\left\|\widehat x-x\right\|
			+\frac{L_{F,2}}{2}\left\|\widehat x-x\right\|^2
			\right]
			\dist(\widehat x,X_\varphi^*)\\
			&+
			\left[
			M_F+
			\frac{L_{f,2}}{2\sigma}
			\dist(x,X_\varphi^*)^2
			\right]
			\dist(\widehat x,X_\varphi^*).
		\end{aligned}
	\end{equation}
	
	We now apply Young's inequality to the first two products on the right-hand side. Specifically,
	\[
	\left[
	\frac{\left\|v\right\|}{\sigma}
	+\alpha\left\|\widehat x-x\right\|
	+\frac{L_{F,2}}{2}\left\|\widehat x-x\right\|^2
	\right]\left\|p-s\right\|\le
	\frac{\kappa_q}{2}\left\|p-s\right\|^2
	+
	\frac{1}{2\kappa_q}
	\left[
	\frac{\left\|v\right\|}{\sigma}
	+\alpha\left\|\widehat x-x\right\|
	+\frac{L_{F,2}}{2}\left\|\widehat x-x\right\|^2
	\right]^2,
	\]
	and, for the second product,
	\[
	\begin{aligned}
		\left[
		\frac{\left\|v\right\|}{\sigma}
		+\alpha\left\|\widehat x-x\right\|
		+\frac{L_{F,2}}{2}\left\|\widehat x-x\right\|^2
		\right]
		\dist(\widehat x,X_\varphi^*)
		&\le
		\frac12
		\left[
		\frac{\left\|v\right\|}{\sigma}
		+\alpha\left\|\widehat x-x\right\|
		+\frac{L_{F,2}}{2}\left\|\widehat x-x\right\|^2
		\right]^2
		+
		\frac12\dist(\widehat x,X_\varphi^*)^2\\
		&\le
		\frac12
		\left[
		\frac{\left\|v\right\|}{\sigma}
		+\alpha\left\|\widehat x-x\right\|
		+\frac{L_{F,2}}{2}\left\|\widehat x-x\right\|^2
		\right]^2
		+
		\frac12\dist(\widehat x,X_\varphi^*),
	\end{aligned}
	\]
	where the second inequality uses \(\dist(\widehat x,X_\varphi^*)\le1\). Substituting these two bounds into \eqref{eq:prop3-young-start} and moving \({\kappa_q}\left\|p-s\right\|^2/2\) to the left-hand side yields
	\[
	\left\|p-s\right\|^2
	\le
	\left(
	\frac{1}{\kappa_q^2}
	+\frac{1}{\kappa_q}
	\right)
	\left[
	\frac{\left\|v\right\|}{\sigma}
	+\alpha\left\|\widehat x-x\right\|
	+\frac{L_{F,2}}{2}\left\|\widehat x-x\right\|^2
	\right]^2+
	\frac{2M_F+1}{\kappa_q}
	\dist(\widehat x,X_\varphi^*)
	+
	\frac{L_{f,2}}{\kappa_q\sigma}
	\dist(x,X_\varphi^*)^2
	\dist(\widehat x,X_\varphi^*).
	\]
	Taking square roots and using \(\sqrt{a+b+c}\le\sqrt a+\sqrt b+\sqrt c\) gives
	\begin{equation}\label{eq:prop3-tangential-bound}
		\begin{aligned}
			\left\|p-s\right\|
			\le{}&
			\sqrt{\frac{1}{\kappa_q^2}+\frac{1}{\kappa_q}}\,
			\left[
			\frac{\left\|v\right\|}{\sigma}
			+\alpha\left\|\widehat x-x\right\|
			+\frac{L_{F,2}}{2}\left\|\widehat x-x\right\|^2
			\right]
			\\
			&+
			\max\left\{
			\sqrt{\frac{2M_F+1}{\kappa_q}},
			\sqrt{\frac{L_{f,2}}{\kappa_q}}
			\right\}
			\left(
			1+\frac{\dist(x,X_\varphi^*)}{\sqrt{\sigma}}
			\right)
			\sqrt{\dist(\widehat x,X_\varphi^*)}.
		\end{aligned}
	\end{equation}
	
	It remains to substitute this estimate into \eqref{eq:prop3-distance-split}. By Proposition~\ref{prop:lower-error-bound}(i),
	\begin{equation}\label{eq:prop3-normal-sqrt}
		\sqrt{\dist(\widehat x,X_\varphi^*)}
		\le
		\sqrt{c_m}\sqrt{R_L(\widehat x)}.
	\end{equation}
	Substituting \eqref{eq:prop3-normal-sqrt} into \eqref{eq:prop3-tangential-bound}, and \eqref{eq:prop3-normal-bound} into \eqref{eq:prop3-distance-split}, proves \eqref{eq:accepted-distance} with the fixed constant
	\[
	c_G
	:=
	\max\left\{
	\max\Bigl\{1,\frac{L_{F,2}}{2}\Bigr\}
	\sqrt{\frac{1}{\kappa_q^2}+\frac{1}{\kappa_q}},\ 
	\sqrt{c_m}\,
	\max\Bigl\{
	\sqrt{\frac{2M_F+1}{\kappa_q}},
	\sqrt{\frac{L_{f,2}}{\kappa_q}}
	\Bigr\}
	+c_m
	\right\}.
	\]
	The proof is completed.
\end{proof}

\begin{remark}\label{rem:regularization-scaling}
	The estimate in Proposition~\ref{prop:measure-distance} also motivates the choice of the regularization term \(\sigma\alpha I\). While regularization ensures the strong convexity of the Newton subproblem, its magnitude must be balanced with the upper-level contribution \(\sigma F\). In particular, the scaling \(\sigma\alpha I\) yields the term \(\alpha\|\widehat x-x\|\) in \eqref{eq:accepted-distance}, whereas an unscaled regularization \(\alpha I\) would introduce the potentially unfavorable factor \(\alpha/\sigma\) as \(\sigma\downarrow0\).
\end{remark}

To apply Proposition~\ref{prop:measure-distance} to the inexact Newton update~\eqref{eq:unit_update}, the following lemma constructs a perturbation \(\eta^k\) satisfying the optimality condition~\eqref{eq:exact-model-certificate} and provides bounds on both \(\eta^k\) and the exact-Hessian residual \(\widetilde r^k\).

\begin{lemma}\label{lem:exact-residual-bridge}
	Let \(x^k\in\operatorname{dom}g\) be given, and suppose that \(\widetilde x^{k+1}\in\operatorname{dom}g\) satisfies \eqref{eq:inexact_unit}, with \(d^k\) and \(x^{k+1}\) defined by \eqref{eq:unit_output} and \eqref{eq:unit_update}. Assume further that
	\begin{equation}\label{eq:hessian-approx}
		\left\|B_k-[\nabla^2 f(x^k)+\sigma_k\nabla^2F(x^k)]\right\|
		\le\kappa_H\sigma_k\alpha_k,
		\qquad 0\le\kappa_H<1.
	\end{equation}
	Define
	\[
	\eta^k
	:=
	(I-H_k)r_k(\widetilde x^{k+1})
	-
	\bigl[
	B_k-\nabla^2 f(x^k)-\sigma_k\nabla^2F(x^k)
	\bigr]d^k.
	\]
	Then \(\eta^k\) satisfies
	\begin{equation}\label{eq:exact-Hess-certificate}
		\eta^k
		-\nabla f(x^k)
		-\sigma_k\nabla F(x^k)
		-
		\bigl[
		\nabla^2 f(x^k)
		+\sigma_k\nabla^2F(x^k)
		+\sigma_k\alpha_k I
		\bigr]d^k
		\in
		\partial g(x^{k+1}),
	\end{equation}
	and the residual $\widetilde r^k$ and perturbation $\eta^k$ bounds
	\[
	\frac{\left\|\widetilde r^k\right\|}{\sigma_k}
	\le\frac{\left\|\eta^k\right\|}{\sigma_k}
	\le
	\alpha_k\bigl[(1+M_H)c_r\tau_k+\kappa_H\left\|d^k\right\|\bigr].
	\]
\end{lemma}
\begin{proof}
	By \eqref{eq:newton_residual}, \eqref{eq:unit_output}, and \eqref{eq:unit_update},
	\[
	x^{k+1}=\widetilde x^{k+1}-r_k(\widetilde x^{k+1})
	=\operatorname{Prox}_g\!\left(
	\widetilde x^{k+1}-\nabla f(x^k)-\sigma_k\nabla F(x^k)
	-H_k(\widetilde x^{k+1}-x^k)\right).
	\]
	Proximal optimality, together with \(\widetilde x^{k+1}-x^{k+1}=r_k(\widetilde x^{k+1})\) and \(\widetilde x^{k+1}-x^k=d^k+r_k(\widetilde x^{k+1})\), yields
	\[
	(I-H_k)r_k(\widetilde x^{k+1})
	-\nabla f(x^k)
	-\sigma_k\nabla F(x^k)
	-H_kd^k
	\in\partial g(x^{k+1}).
	\]
	Substituting \((I-H_k)r_k(\widetilde x^{k+1})=\eta^k+\bigl[B_k-\nabla^2 f(x^k)-\sigma_k\nabla^2F(x^k)\bigr]d^k\) and \(H_k=B_k+\sigma_k\alpha_kI\) rewrites this inclusion exactly as \eqref{eq:exact-Hess-certificate}.
	
	We next prove the bounds. By the definition of \(r_k(\widetilde x^{k+1})\) in \eqref{eq:newton_residual} and the update rules \eqref{eq:unit_output}--\eqref{eq:unit_update},
	\[
	x^{k+1}=
	\widetilde x^{k+1}-r_k(\widetilde x^{k+1})=
	\operatorname{Prox}_g\!\left(
	\widetilde x^{k+1}
	-\nabla f(x^k)
	-\sigma_k\nabla F(x^k)
	-H_k(\widetilde x^{k+1}-x^k)
	\right).
	\]
	Hence, by \eqref{eq:exact_hessian_residual} and the nonexpansiveness of \(\operatorname{Prox}_g\),
	\[
	\begin{aligned}
		\left\|\widetilde r^k\right\|
		&=
		\Bigl\|
		\operatorname{Prox}_g\!\left(
		\widetilde x^{k+1}
		-\nabla f(x^k)
		-\sigma_k\nabla F(x^k)
		-H_k(\widetilde x^{k+1}-x^k)
		\right)\\
		&\qquad\quad
		-
		\operatorname{Prox}_g\!\left(
		x^{k+1}
		-\nabla f(x^k)
		-\sigma_k\nabla F(x^k)
		-\bigl[
		\nabla^2 f(x^k)
		+\sigma_k\nabla^2F(x^k)
		+\sigma_k\alpha_k I
		\bigr]d^k
		\right)
		\Bigr\|\\
		&\le
		\Bigl\|
		\widetilde x^{k+1}-x^{k+1}
		-H_k(\widetilde x^{k+1}-x^k)+
		\bigl[
		\nabla^2 f(x^k)
		+\sigma_k\nabla^2F(x^k)
		+\sigma_k\alpha_k I
		\bigr]d^k
		\Bigr\|\\
		&=
		\left\|\eta^k\right\|
		\\
		&\le
		\left\|I-H_k\right\|\,
		\left\|r_k(\widetilde x^{k+1})\right\|
		+
		\left\|
			B_k-\nabla^2 f(x^k)-\sigma_k\nabla^2F(x^k)
		\right\|\,
		\left\|d^k\right\|
		\\
		&\le
		(1+M_H)\left\|r_k(\widetilde x^{k+1})\right\|
		+
		\kappa_H\sigma_k\alpha_k\left\|d^k\right\|,
	\end{aligned}
	\]
	where the second equality uses $\widetilde x^{k+1}-x^{k+1} = r_k(\widetilde x^{k+1})$, $\widetilde x^{k+1}-x^k = r_k(\widetilde x^{k+1})+d^k$, and \(H_k=B_k+\sigma_k\alpha_k I\), so that the vector in that line is precisely \(\eta^k\). The second inequality follows from the triangle inequality, and the last inequality uses \(\left\|I-H_k\right\|\le1+\left\|H_k\right\|\le1+M_H\) and \eqref{eq:hessian-approx}.
	
	Finally, \eqref{eq:inexact_unit} gives
	\[
	\frac{\left\|\widetilde r^k\right\|}{\sigma_k}
	\le
	\frac{\left\|\eta^k\right\|}{\sigma_k}
	\le
	(1+M_H)
	\frac{\left\|r_k(\widetilde x^{k+1})\right\|}{\sigma_k}
	+
	\kappa_H\alpha_k\left\|d^k\right\|\le
	\alpha_k
	\bigl[
	(1+M_H)c_r\tau_k
	+\kappa_H\left\|d^k\right\|
	\bigr].
	\]
	The proof is completed.
\end{proof}

We next express the lower-level feasibility error directly in terms of the current measure. The resulting \(O(\omega_k^2)\) estimate controls the scaled term \(\dist(x^k,X_\varphi^*)/\sqrt{\sigma_k}\) when Proposition~\ref{prop:measure-distance} is applied to the Newton update.

\begin{lemma}\label{lem:lower-distance}
	Suppose that the standing assumptions hold and that Assumption~\ref{asup:lower_regular} holds at \(\bar x\in X_\varphi^*\). Let \(\{x^k\}\) be the iterates generated by the inexact proximal Newton method with the parameters selected according to \eqref{eq:parameter-rule}. Then, for any \(k\ge0\) such that \(x^k\in\mathbb B_{\epsilon_L}(\bar x)\),
	\[
	\dist(x^k,X_\varphi^*)\le c_mG_k^2,
	\]
	where \(G_k\) is given by \eqref{eq:local_trial_measure} for \(k\ge1\), with \(G_0\ge\sqrt{R_L(x^0)}\), and \(\epsilon_L\) and \(c_m\) are as in Proposition~\ref{prop:lower-error-bound}(i).
\end{lemma}
\begin{proof}
	Let \(\epsilon_L,c_m>0\) be the fixed radius and modulus of Proposition~\ref{prop:lower-error-bound}(i), so that $\dist(x,X_\varphi^*)\le c_mR_L(x)$ for all $x\in\mathbb B_{\epsilon_L}(\bar x)$.
	
	For \(k\ge1\), the update \eqref{eq:local_trial_measure} gives \(G_k\ge\sqrt{R_L(x^k)}\); the initialization gives the same inequality at \(k=0\). Hence, for every \(k\ge0\), $R_L(x^k)\le G_k^2$. Therefore, whenever \(x^k\in\mathbb B_{\epsilon_L}(\bar x)\),
	\[
	\dist(x^k,X_\varphi^*)
	\le
	c_mR_L(x^k)
	\le
	c_mG_k^2.
	\]
	The proof is completed.
\end{proof}

Combining Proposition~\ref{prop:measure-distance} with the perturbation and feasibility bounds of Lemmas~\ref{lem:exact-residual-bridge} and \ref{lem:lower-distance} gives local control of the actual Newton output relative to \(S^*\).

\begin{lemma}\label{lem:measure-distance-iteration}
	Suppose that the standing assumptions and Assumption~\ref{asup:common-image} hold. Fix \(\bar x\in S^*\) and suppose that Assumptions~\ref{asup:lower_regular}, \ref{asup_quadratic_growth}, and \ref{asup_hessian_lipschitz} hold at \(\bar x\). Let \(\{x^k\}\) be the iterates generated by the inexact proximal Newton method with the parameters selected according to \eqref{eq:parameter-rule}. Then there exist \(\epsilon_S,\delta_S,c_S>0\) such that, for any \(k\ge0\) satisfying \eqref{eq:inexact_unit} and \eqref{eq:hessian-approx} at iteration \(k\), \(x^k,x^{k+1}\in\mathbb B_{\epsilon_S}(\bar x)\), and \(\omega_k\le\delta_S\), we have
	\begin{equation}\label{eq:accepted-distance-iteration}
		\dist(x^{k+1},S^*)\le
		c_S \omega_{k+1}.
	\end{equation}
	where \(\omega_j=G_j+\tau_j\) for \(j=k,k+1\).
\end{lemma}

\begin{proof}
	By Lemma~\ref{lem:exact-residual-bridge}, the model perturbation \(\eta^k\) satisfies \eqref{eq:exact-Hess-certificate}, so Proposition~\ref{prop:measure-distance} applies with \(x=x^k\), \(\widehat x=x^{k+1}\), \(\sigma=\sigma_k\), \(\alpha=\alpha_k\), and \(v=\eta^k\). The same lemma also gives
	\begin{equation}\label{eq:certificate-bound}
		\frac{\left\|\eta^k\right\|}{\sigma_k}
		\le
		(1+M_H)c_r\alpha_k\tau_k
		+\kappa_H\alpha_k\left\|d^k\right\|.
	\end{equation}
	
	Let \(\epsilon_G,c_G>0\) be the fixed radius and constant of Proposition~\ref{prop:measure-distance}, and let \(\epsilon_L,c_m>0\) be the radius and modulus used in Lemma~\ref{lem:lower-distance}. Set
	\[
	\epsilon_S:=\min\{\epsilon_G,\epsilon_L\},
	\qquad
	\delta_S
	:=
	\min\left\{
	1,
	\left(\frac{\bar\sigma}{c_\sigma}\right)^{1/(2(1+\rho))}
	\right\}.
	\]
	Then \(\omega_k\le\delta_S\) implies that the penalty cap in \eqref{eq:parameter-rule} is inactive and hence $\sigma_k=c_\sigma\omega_k^{2(1+\rho)}$. Proposition~\ref{prop:measure-distance} therefore gives
	\begin{equation}\label{eq:iteration-distance-from-prop}
		\dist(x^{k+1},S^*)
		\le
		c_G\Biggl[
		\frac{\left\|\eta^k\right\|}{\sigma_k}
		+\left\|d^k\right\|^2
		+\alpha_k\left\|d^k\right\|+
		\left(
		1+
		\frac{\dist(x^k,X_\varphi^*)}{\sqrt{\sigma_k}}
		\right)
		\sqrt{R_L(x^{k+1})}
		\Biggr].
	\end{equation}
	
	By Lemma~\ref{lem:lower-distance}, \(G_k\le\omega_k\), and the expression for \(\sigma_k\),
	\begin{equation}\label{eq:iteration-penalty-ratio}
		\frac{\dist(x^k,X_\varphi^*)}{\sqrt{\sigma_k}}\le\frac{c_mG_k^2}{\sqrt{c_\sigma}\,\omega_k^{1+\rho}}\le\frac{c_m}{\sqrt{c_\sigma}}\omega_k^{1-\rho}\le\frac{c_m}{\sqrt{c_\sigma}},
	\end{equation}
	where the last inequality follows from \(\omega_k\le1\) and \(0<\rho\le1\).
	
	Substituting \eqref{eq:certificate-bound} and \eqref{eq:iteration-penalty-ratio} into \eqref{eq:iteration-distance-from-prop} yields
	\[
	\dist(x^{k+1},S^*)
	\le
	c_G\Biggl[
	(1+M_H)c_r\alpha_k\tau_k
	+(\kappa_H+1)\alpha_k\left\|d^k\right\|
	+\left\|d^k\right\|^2+
	\left(
	1+\frac{c_m}{\sqrt{c_\sigma}}
	\right)
	\sqrt{R_L(x^{k+1})}
	\Biggr].
	\]
	Since
	\[
	G_{k+1}
	=
	\frac{\left\|\widetilde r^k\right\|}{\sigma_k}
	+\left\|d^k\right\|^2
	+\alpha_k\left\|d^k\right\|
	+\sqrt{R_L(x^{k+1})},
	\]
	and all four terms are nonnegative,
	\[
	\dist(x^{k+1},S^*)
	\le
	c_0\bigl(G_{k+1}+\alpha_k\tau_k\bigr),
	\]
	where $c_0:=c_G\max\left\{(1+M_H)c_r,\,\kappa_H+3+c_m/\sqrt{c_\sigma}\right\}$.
	
	By the parameter rule $\tau_{k+1} =\min\{\gamma\tau_k,\,c_\tau\omega_k^{1+\rho}\}$, together with \(\alpha_k\le\bar\alpha\), \(\alpha_k\le c_\alpha\omega_k^\rho\), and \(\tau_k\le\omega_k\), we have
	\[
	\alpha_k\tau_k\le
	\min\left\{
	\frac{\bar\alpha}{\gamma}\,
	\gamma\tau_k,\,
	\frac{c_\alpha}{c_\tau}\,
	c_\tau\omega_k^{1+\rho}
	\right\}\le
	\max\left\{
	\frac{\bar\alpha}{\gamma},
	\frac{c_\alpha}{c_\tau}
	\right\}
	\tau_{k+1}.
	\]
	Therefore,
	\begin{align*}
		\dist(x^{k+1},S^*)
		\le
		c_0\bigl(G_{k+1}+\alpha_k\tau_k\bigr) \le
		c_0
		\max\left\{
		1,\frac{\bar\alpha}{\gamma},
		\frac{c_\alpha}{c_\tau}
		\right\}
		\bigl(G_{k+1}+\tau_{k+1}\bigr)=
		c_S\omega_{k+1},
	\end{align*}
	where $c_S:=c_0	\max\left\{	1,\bar\alpha/\gamma,c_\alpha/c_\tau	\right\}$. This proves \eqref{eq:accepted-distance-iteration}.
\end{proof}

To estimate \(G_{k+1}\), we consider the individual terms in its definition. The exact-Hessian residual term has already been controlled in Lemma~\ref{lem:exact-residual-bridge}. We next estimate the Newton displacement \(d^k\), which will also be used in the subsequent analysis of \(R_L(x^{k+1})\).

\begin{lemma}\label{lem:step_estimate}
	Suppose that the standing assumptions and Assumption~\ref{asup:common-image} hold, and that Assumptions~\ref{asup:lower_regular} and \ref{asup_hessian_lipschitz} hold at \(\bar x\in S^*\). Let \(\{x^k\}\) be the iterates generated by the inexact proximal Newton method with the parameters selected according to \eqref{eq:parameter-rule}. Given \(c_b>0\), there exist \(\epsilon_d,\delta_d,c_d>0\) such that, whenever \eqref{eq:inexact_unit} and \eqref{eq:hessian-approx} hold at iteration \(k\), \(x^k\in\mathbb B_{\epsilon_d}(\bar x)\), \(\omega_k\le\delta_d\), and \(\dist(x^k,S^*)\le c_b\omega_k\), we have
	\begin{equation}\label{eq:step_estimate}
		\left\|d^k\right\|\le c_d\omega_k,
	\end{equation}
	where \(\omega_k=G_k+\tau_k\).
\end{lemma}
\begin{proof}
	Fix \(c_b>0\); the radius, threshold, and constant constructed below may depend on this fixed value of \(c_b\), but are independent of \(k\). Let \(q_k^0\) be the model \eqref{eq:q_model_unit} with \(B_k\) replaced by \(\nabla^2 f(x^k)+\sigma_k\nabla^2F(x^k)\), retaining the regularization \(\sigma_k\alpha_kI\), and let \(y^0\) be its unique minimizer. Existence and uniqueness follow from \(\sigma_k\alpha_k\)-strong convexity. The target step satisfies
	\[
	\left\|d^k\right\|=\left\|x^{k+1}-x^k\right\|
	\le\left\|x^{k+1}-y^0\right\|+\left\|y^0-x^k\right\|.
	\]
	We first bound \(\left\|y^0-x^k\right\|\) by comparing exact-model values near \(S^*\), and then bound \(\left\|x^{k+1}-y^0\right\|\) using the model correction and the Hessian-approximation condition.
	
	Let \(\epsilon_L,\epsilon_f\) denote the fixed radii of Proposition~\ref{prop:lower-error-bound} and Lemma~\ref{lem:lower-model-growth}, respectively. Let \(\epsilon_{F,2},L_{F,2},M_F>0\) be the fixed radius and bounds following Assumption~\ref{asup_hessian_lipschitz}, so that
	\[
	\left\|\nabla F(u)\right\|\le M_F,\qquad
	\left\|\nabla^2F(u)-\nabla^2F(v)\right\|\le L_{F,2}\left\|u-v\right\|
	\qquad\ \forall\ u,v\in\mathbb B_{\epsilon_{F,2}}(\bar x).
	\]
	Set the common radius
	\[
	\epsilon_d:=\tfrac12\min\{\epsilon_L,\epsilon_f,\epsilon_{F,2}\}.
	\]
	Thus all these local estimates hold on \(\mathbb B_{2\epsilon_d}(\bar x)\). The common threshold \(\delta_d\), specified below, makes the parameter caps inactive and ensures \(L_0\omega_k\le\epsilon_d\).
	
	Set \(s:=\operatorname{Proj}_{S^*}(x^k)\). Since the projection onto \(S^*\) fixes \(\bar x\), nonexpansiveness gives $\left\|s-\bar x\right\|\le\left\|x^k-\bar x\right\|\le\epsilon_d$. The current-state inputs of the proof are the lower-feasibility estimate of Lemma~\ref{lem:lower-distance} at \(k\) and the assumed distance estimate of the statement:
	\begin{equation}\label{eq:step-input-distances}
		\dist(x^k,X_\varphi^*)\le c_mG_k^2\le c_m\omega_k^2,\qquad
		\left\|x^k-s\right\|=\dist(x^k,S^*)\le c_b\omega_k.
	\end{equation}
	
	For \(w\in\mathbb B_{2\epsilon_d}(\bar x) \cap\operatorname{dom}g\), set \(z:=\operatorname{Proj}_{X_\varphi^*}(w)\). Likewise, \(\left\|z-\bar x\right\|\le\left\|w-\bar x\right\|\le2\epsilon_d\), and \(\left\|w-z\right\|=\dist(w,X_\varphi^*)\). Expanding the model difference separates the lower-level, upper-level, and regularization contributions:
	\begin{equation}\label{eq:exact-model-expansion}
		\begin{aligned}
			q_k^0(w)-q_k^0(s)
			={}&g(w)-g(s)+\left\langle \nabla f(x^k),w-s\right\rangle+\frac12\Bigl[
			\left\langle w-x^k,\nabla^2 f(x^k)(w-x^k)\right\rangle
			-\left\langle s-x^k,\nabla^2 f(x^k)(s-x^k)\right\rangle\Bigr]\\
			&+\sigma_k\Bigl[
			\left\langle \nabla F(x^k)+\nabla^2F(x^k)(s-x^k),w-s\right\rangle
			+\frac12\left\langle w-s,\nabla^2F(x^k)(w-s)\right\rangle\Bigr]\\
			&+\frac{\sigma_k\alpha_k}{2}
			\bigl[\left\|w-x^k\right\|^2-\left\|s-x^k\right\|^2\bigr].
		\end{aligned}
	\end{equation}
	By Lemma~\ref{lem:lower-model-growth}, the lower-level contribution satisfies
	\[
	\begin{aligned}
		&\quad\ g(w)-g(s)+\left\langle \nabla f(x^k),w-s\right\rangle+\frac12\Bigl[
		\left\langle w-x^k,\nabla^2 f(x^k)(w-x^k)\right\rangle
		-\left\langle s-x^k,\nabla^2 f(x^k)(s-x^k)\right\rangle\Bigr]\\
		&\ge\frac{c_L}{2}\left\|w-z\right\|^2
		-\frac{L_{f,2}}{2}\dist(x^k,X_\varphi^*)^2\left\|w-z\right\|.
	\end{aligned}
	\]
	For the upper-level contribution, constrained optimality at \(s\in S^*\) gives \(\left\langle \nabla F(s),z-s\right\rangle\ge0\). Splitting \(w-s=(w-z)+(z-s)\), we obtain
	\[
	\left\langle \nabla F(s),w-s\right\rangle=\left\langle \nabla F(s),w-z\right\rangle+\left\langle \nabla F(s),z-s\right\rangle\ge-\left\|\nabla F(s)\right\|\left\|w-z\right\|
	\ge-M_F\left\|w-z\right\|,
	\]
	where the inequalities use the Cauchy--Schwarz inequality and the local gradient bound. The Hessian Lipschitz estimate on \([x^k,s]\) gives
	\[
	\left\|\nabla F(x^k)+\nabla^2F(x^k)(s-x^k)-\nabla F(s)\right\|
	\le\frac{L_{F,2}}{2}\left\|x^k-s\right\|^2.
	\]
	Adding and subtracting \(\nabla F(s)\) in the upper-model inner product, and using \(\nabla^2F(x^k)\succeq0\), therefore yields
	\[
	\begin{aligned}
		&\quad\ \left\langle \nabla F(x^k)+\nabla^2F(x^k)(s-x^k),w-s\right\rangle
		+\frac12\left\langle w-s,\nabla^2F(x^k)(w-s)\right\rangle\\
		&\ge\left\langle \nabla F(s),w-s\right\rangle
		+\left\langle \nabla F(x^k)+\nabla^2F(x^k)(s-x^k)-\nabla F(s),w-s\right\rangle\\
		&\ge-M_F\left\|w-z\right\|
		-\frac{L_{F,2}}{2}\left\|x^k-s\right\|^2\left\|w-s\right\|.
	\end{aligned}
	\]
	Finally, the regularization contribution is controlled by the identity
	\[
	\left\|w-x^k\right\|^2-\left\|s-x^k\right\|^2
	=2\left\langle s-x^k,w-s\right\rangle+\left\|w-s\right\|^2
	\ge\left\|w-s\right\|^2-2\left\|x^k-s\right\|\left\|w-s\right\|.
	\]
	
	Substituting the preceding three bounds into \eqref{eq:exact-model-expansion} gives
	\begin{equation}\label{eq:exact-model-lower-bound}
		\begin{aligned}
			q_k^0(w)-q_k^0(s)
			\ge{}&
			\frac{c_L}{2}\left\|w-z\right\|^2
			-\left[
			\frac{L_{f,2}}{2}\dist(x^k,X_\varphi^*)^2
			+\sigma_kM_F
			\right]\left\|w-z\right\|
			\\
			&+\frac{\sigma_k\alpha_k}{2}\left\|w-s\right\|^2
			-\sigma_k\left[
			\alpha_k\left\|x^k-s\right\|
			+\frac{L_{F,2}}{2}\left\|x^k-s\right\|^2
			\right]\left\|w-s\right\|
			\\
			\ge{}&
			\frac{\sigma_k\alpha_k}{2}\left\|w-s\right\|^2
			-\sigma_k\left[
			\alpha_k\left\|x^k-s\right\|
			+\frac{L_{F,2}}{2}\left\|x^k-s\right\|^2
			\right]\left\|w-s\right\|
			\\
			&-
			\frac{
				\left[
				\frac{L_{f,2}}{2}\dist(x^k,X_\varphi^*)^2
				+\sigma_kM_F
				\right]^2
			}{2c_L}.
		\end{aligned}
	\end{equation}
	For \(w\notin\operatorname{dom}g\), the bound \eqref{eq:exact-model-lower-bound} holds trivially because \(q_k^0(w)=+\infty\).
	
	To make the model lower bound positive on a sphere around \(s\), define the fixed coefficient
	\[
	\Lambda:=\max\left\{
	c_b+\frac{L_{F,2}c_b^2}{2c_\alpha},\,
	\frac{\bigl(L_{f,2}c_m^2/2+c_\sigma M_F\bigr)^2}
	{2c_Lc_\sigma c_\alpha}\right\},
	\qquad L_0:=2(\Lambda+1).
	\]
	We now fix the common threshold
	\[
	\delta_d:=\min\left\{
	1,\,
	\left(\frac{\bar\sigma}{c_\sigma}\right)^{1/[2(1+\rho)]},\,
	\left(\frac{\bar\alpha}{c_\alpha}\right)^{1/\rho},\,
	\frac{\epsilon_d}{L_0}\right\}.
	\]
	Then \(0<\omega_k\le\delta_d\) makes both caps in \eqref{eq:parameter-rule} inactive, so $\sigma_k=c_\sigma\omega_k^{2(1+\rho)}$, $\alpha_k=c_\alpha\omega_k^\rho$, $\omega_k\le1$. Using \eqref{eq:step-input-distances}, we obtain
	\[
	\frac{\left\|x^k-s\right\|}{\omega_k}
	+\frac{L_{F,2}\left\|x^k-s\right\|^2}{2\alpha_k\omega_k}
	\le c_b+\frac{L_{F,2}c_b^2}{2c_\alpha}\omega_k^{1-\rho}
	\le\Lambda.
	\]
	Moreover, \(4\ge2(1+\rho)\) and \(\omega_k\le1\) give
	\[
	\begin{aligned}
		\frac{L_{f,2}}2\dist(x^k,X_\varphi^*)^2+\sigma_kM_F
		\le\frac{L_{f,2}c_m^2}{2}\omega_k^4
		+c_\sigma M_F\omega_k^{2(1+\rho)}\le\left(\frac{L_{f,2}c_m^2}{2}+c_\sigma M_F\right)
		\omega_k^{2(1+\rho)},\\
	\end{aligned}
	\]
	and
	\[
	\frac{\left[\frac{L_{f,2}}2\dist(x^k,X_\varphi^*)^2
		+\sigma_kM_F\right]^2}
	{2c_L\sigma_k\alpha_k\omega_k^2}
	\le\Lambda\omega_k^\rho\le\Lambda.
	\]
	Since \(L_0\omega_k\le L_0\delta_d\le\epsilon_d\), every point in the ball \(\left\|w-s\right\|\le L_0\omega_k\) satisfies
	\[
	\left\|w-\bar x\right\|\le\left\|w-s\right\|+\left\|s-\bar x\right\|
	\le L_0\omega_k+\left\|x^k-\bar x\right\|\le2\epsilon_d.
	\]
	Thus \eqref{eq:exact-model-lower-bound} applies throughout this ball. Dividing it by \(\sigma_k\alpha_k\omega_k^2\) and using the two normalized error bounds gives, on \(\left\|w-s\right\|=L_0\omega_k\),
	\begin{equation}\label{eq:exact-model-sphere-bound}
		\frac{q_k^0(w)-q_k^0(s)}{\sigma_k\alpha_k\omega_k^2}
		\ge\frac{L_0^2}{2}-\Lambda L_0-\Lambda
		=\Lambda+2>0.
	\end{equation}
	If \(\left\|y^0-s\right\|\ge L_0\omega_k\), set \(\theta:=L_0\omega_k/\left\|y^0-s\right\|\in(0,1]\) and \(w:=(1-\theta)s+\theta y^0\). Then \(w\) lies on this sphere, whereas convexity and minimality of \(y^0\) imply
	\[
	q_k^0(w)\le(1-\theta)q_k^0(s)+\theta q_k^0(y^0)\le q_k^0(s),
	\]
	contradicting \eqref{eq:exact-model-sphere-bound}. Hence \(\left\|y^0-s\right\|<L_0\omega_k\), and \eqref{eq:step-input-distances} yields
	\begin{equation}\label{eq:exact-model-localization}
		\left\|y^0-x^k\right\|\le\left\|y^0-s\right\|+\left\|s-x^k\right\|
		\le(L_0+c_b)\omega_k.
	\end{equation}
	
	It remains to compare \(x^{k+1}\) with \(y^0\). The model correction \eqref{eq:unit_output} and exact-model optimality give, respectively,
	\[
	u^{k+1}:=(I-H_k)r_k(\widetilde x^{k+1})
	-\nabla f(x^k)-\sigma_k\nabla F(x^k)-H_kd^k
	\in\partial g(x^{k+1}),
	\]
	and
	\[
	u^0:=-\nabla f(x^k)-\sigma_k\nabla F(x^k)-\bigl[\nabla^2 f(x^k)+\sigma_k\nabla^2F(x^k)
	+\sigma_k\alpha_kI\bigr](y^0-x^k)
	\in\partial g(y^0).
	\]
	Monotonicity of \(\partial g\) gives \(\left\langle u^{k+1}-u^0,x^{k+1}-y^0\right\rangle\ge0\). Substituting these two subgradients and then using \(d^k=(x^{k+1}-y^0)+(y^0-x^k)\) and \(H_k=B_k+\sigma_k\alpha_kI\) gives
	\[
	\begin{aligned}
		0\le{}&\left\langle u^{k+1}-u^0,x^{k+1}-y^0\right\rangle\\
		={}&\left\langle (I-H_k)r_k(\widetilde x^{k+1}),x^{k+1}-y^0\right\rangle
		-\left\langle H_k(x^{k+1}-y^0),x^{k+1}-y^0\right\rangle\\
		&-\left\langle \bigl[B_k-\nabla^2 f(x^k)-\sigma_k\nabla^2F(x^k)\bigr]
			(y^0-x^k),x^{k+1}-y^0\right\rangle.
	\end{aligned}
	\]
	Since \(H_k\succeq\sigma_k\alpha_kI\), rearranging gives
	\begin{equation}\label{eq:approx-exact-model-comparison}
		\sigma_k\alpha_k\left\|x^{k+1}-y^0\right\|^2
		\le\bigl\langle (I-H_k)r_k(\widetilde x^{k+1})-\bigl[B_k-\nabla^2 f(x^k)-\sigma_k\nabla^2F(x^k)\bigr]
		(y^0-x^k),\ x^{k+1}-y^0\bigr\rangle.
	\end{equation}
	For \(x^{k+1}\ne y^0\), applying the Cauchy--Schwarz inequality to \eqref{eq:approx-exact-model-comparison}, dividing by \(\sigma_k\alpha_k\left\|x^{k+1}-y^0\right\|\), and using \(\left\|I-H_k\right\|\le1+M_H\) and \eqref{eq:hessian-approx} yields
	\[
	\begin{aligned}
		\left\|x^{k+1}-y^0\right\|
		&\le(1+M_H)\frac{\left\|r_k(\widetilde x^{k+1})\right\|}{\sigma_k\alpha_k}
		+\kappa_H\left\|y^0-x^k\right\|\\
		&\le(1+M_H)c_r\tau_k+\kappa_H\left\|y^0-x^k\right\|\\
		&\le\bigl[(1+M_H)c_r+\kappa_H(L_0+c_b)\bigr]\omega_k.
	\end{aligned}
	\]
	The second inequality uses \eqref{eq:inexact_unit}, the third uses \eqref{eq:exact-model-localization} and \(\tau_k\le\omega_k\). The same bound holds immediately if \(x^{k+1}=y^0\). Returning to the target step and using \eqref{eq:exact-model-localization}, we obtain
	\[
	\left\|d^k\right\|\le\left\|x^{k+1}-y^0\right\|+\left\|y^0-x^k\right\|\le\bigl[(1+M_H)c_r+(1+\kappa_H)(L_0+c_b)\bigr]\omega_k.
	\]
	Thus \eqref{eq:step_estimate} holds with the fixed constant $c_d:=(1+M_H)c_r+(1+\kappa_H)(L_0+c_b)$.
\end{proof}

Having established the estimate for \(d^k\), we next control the remaining lower-level residual term \(R_L(x^{k+1})\) in the definition of \(G_{k+1}\).

\begin{lemma}\label{lem:lower-stationarity}
	Suppose that the standing assumptions and Assumption~\ref{asup:common-image} hold, and that Assumptions~\ref{asup:lower_regular} and \ref{asup_hessian_lipschitz} hold at \(\bar x\in S^*\). Let \(\{x^k\}\) be the iterates generated by the inexact proximal Newton method with the parameters selected according to \eqref{eq:parameter-rule}. Given \(c_b>0\), there exist \(\epsilon_R,\delta_R,c_R>0\) such that, whenever \eqref{eq:inexact_unit} and \eqref{eq:hessian-approx} hold at iteration \(k\), \(x^k\in\mathbb B_{\epsilon_R}(\bar x)\), \(\omega_k\le\delta_R\), and \(\dist(x^k,S^*)\le c_b\omega_k\), we have
	\begin{equation}\label{eq:lower_stationarity_order}
		R_L(x^{k+1})\le c_R\omega_k^{2(1+\rho)},
	\end{equation}
	where \(\omega_k=G_k+\tau_k\).
\end{lemma}
\begin{proof}
	Fix \(c_b>0\), and let \(\epsilon_d,\delta_d,c_d\) be the constants furnished by Lemma~\ref{lem:step_estimate} for this value of \(c_b\). Set
	\[
	\epsilon_R:=\tfrac12\min\left\{\epsilon_L,\epsilon_d,\epsilon_m,
	\epsilon_{F,2},\epsilon_{f,2},\frac{1}{2\kappa_m\max\{L_{f,2},1\}}\right\},
	\qquad
	\delta_R:=\min\left\{\delta_d,\frac{\epsilon_R}{c_d}\right\},
	\]
	where \(\epsilon_L\) is the radius of Proposition~\ref{prop:lower-error-bound}, \(\epsilon_m,\kappa_m\) are those of Assumption~\ref{asup:lower_regular}, and \(\epsilon_{F,2},\epsilon_{f,2},L_{f,2}\) are those fixed after Assumption~\ref{asup_hessian_lipschitz}. Whenever \(x^k\in\mathbb B_{\epsilon_R}(\bar x)\), \(\omega_k\le\delta_R\), and \(\dist(x^k,S^*)\le c_b\omega_k\), Lemma~\ref{lem:step_estimate} applies at index \(k\) and gives
	\[
	\left\|x^{k+1}-\bar x\right\|
	\le\left\|x^k-\bar x\right\|+\left\|d^k\right\|
	\le\left\|x^k-\bar x\right\|+c_d\omega_k\le\epsilon_R+c_d\delta_R\le2\epsilon_R.
	\]
	Hence \(x^k,x^{k+1}\in\mathbb B_{2\epsilon_R}(\bar x)\), where the metric subregularity of Assumption~\ref{asup:lower_regular}, the Hessian Lipschitz estimates of Assumption~\ref{asup_hessian_lipschitz} (applied at the image points \(Ax^k,Ax^{k+1}\in\mathbb B_{\epsilon_{h,2}}(A\bar x)\)), and the bound \(\left\|\nabla F\right\|\le M_F\) all hold.
	
	By \eqref{eq:residual_by_subgradient}, it suffices to construct \(v^k\in\partial\varphi(x^{k+1})\) and bound \(\left\|v^k\right\|\). Adding \(\nabla f(x^{k+1})\) to the perturbed optimality condition \eqref{eq:exact-Hess-certificate} of Lemma~\ref{lem:exact-residual-bridge} and using \(\partial\varphi(x^{k+1})=\nabla f(x^{k+1})+\partial g(x^{k+1})\), we obtain \(v^k\in\partial\varphi(x^{k+1})\), where
	\[
	\begin{aligned}
		v^k
		:={}&
		\bigl[
		\nabla f(x^{k+1})
		-\nabla f(x^k)
		-\nabla^2 f(x^k)d^k
		\bigr]-
		\sigma_k
		\bigl[
		\nabla F(x^k)
		+\nabla^2F(x^k)d^k
		+\alpha_kd^k
		\bigr]
		+\eta^k.
	\end{aligned}
	\]
	Fix any \(s\in X_\varphi^*\). By Proposition~\ref{prop:simple-blo-geometry},
	\[
	\left\|Ad^k\right\|\le
	\left\|A(x^{k+1}-s)\right\|
	+
	\left\|A(x^k-s)\right\|\le
	\left\| A\right\|
	\bigl[
	\dist(x^{k+1},X_\varphi^*)
	+
	\dist(x^k,X_\varphi^*)
	\bigr].
	\]
	Therefore, by the Hessian Lipschitz continuity of \(h\),
	\begin{equation}\label{eq:lower-stationarity-taylor}
		\left\|
			\nabla f(x^{k+1})
			-\nabla f(x^k)
			-\nabla^2 f(x^k)d^k
		\right\|
		\le
		\frac{L_{f,2}}{2}
		\bigl[
		\dist(x^k,X_\varphi^*)
		+\dist(x^{k+1},X_\varphi^*)
		\bigr]^2\le
		L_{f,2}
		\bigl[
		\dist(x^k,X_\varphi^*)^2
		+\dist(x^{k+1},X_\varphi^*)^2
		\bigr].
	\end{equation}
	
	Moreover, since \(\omega_k\le1\) and \(\left\|d^k\right\|\le c_d\omega_k\), the standing Lipschitz bound on \(\nabla F\) and the local bound \(\left\|\nabla F(x^k)\right\|\le M_F\) give
	\begin{equation}\label{eq:lower-stationarity-upper}
		\sigma_k
		\left\|
			\nabla F(x^k)
			+\nabla^2F(x^k)d^k
			+\alpha_kd^k
		\right\|\le
		\sigma_k
		\left[
		M_F
		+(L_F+\bar\alpha)\left\|d^k\right\|
		\right]\le
		\sigma_k
		\left[
		M_F
		+(L_F+\bar\alpha)c_d
		\right].
	\end{equation}
	By Lemma~\ref{lem:exact-residual-bridge},
	\begin{equation}\label{eq:lower-stationarity-model-error}
		\left\|\eta^k\right\|
		\le
		\sigma_k\alpha_k
		\bigl[
		(1+M_H)c_r\tau_k
		+\kappa_H\left\|d^k\right\|
		\bigr]
		\le
		\sigma_k\bar\alpha
		\bigl[
		(1+M_H)c_r+\kappa_Hc_d
		\bigr],
	\end{equation}
	where the last inequality uses \(\tau_k\le\omega_k\le1\), \(\left\|d^k\right\|\le c_d\omega_k\le c_d\), and \(\alpha_k\le\bar\alpha\).
	
	Combining \eqref{eq:lower-stationarity-taylor}, \eqref{eq:lower-stationarity-upper}, and \eqref{eq:lower-stationarity-model-error} gives
	\begin{equation}\label{eq:lower-stationarity-v-bound}
		\left\|v^k\right\|
		\le
		L_{f,2}\dist(x^k,X_\varphi^*)^2
		+
		L_{f,2}\dist(x^{k+1},X_\varphi^*)^2+
		\Bigl[
		M_F
		+(L_F+\bar\alpha)c_d
		+\bar\alpha\bigl((1+M_H)c_r+\kappa_Hc_d\bigr)
		\Bigr]\sigma_k.
	\end{equation}
	
	Metric subregularity at \(x^{k+1}\), together with \eqref{eq:lower-stationarity-v-bound}, therefore yields
	\[
	\begin{aligned}
		&\quad\ \dist(x^{k+1},X_\varphi^*)\\
		&\le
		\kappa_m\dist(0,\partial\varphi(x^{k+1}))\\
		&\le
		\kappa_mL_{f,2}\dist(x^k,X_\varphi^*)^2
		+
		\kappa_mL_{f,2}\dist(x^{k+1},X_\varphi^*)^2+
		\kappa_m
		\Bigl[
		M_F
		+(L_F+\bar\alpha)c_d
		+\bar\alpha\bigl((1+M_H)c_r+\kappa_Hc_d\bigr)
		\Bigr]\sigma_k.
	\end{aligned}
	\]
	
	Since \(\left\|x^{k+1}-\bar x\right\|\le2\epsilon_R\), \(2\kappa_mL_{f,2}\epsilon_R\le1/2\), and \(\bar x\in X_\varphi^*\),
	\[
	\kappa_mL_{f,2}
	\dist(x^{k+1},X_\varphi^*)^2
	\le \kappa_mL_{f,2} \left\|x^{k+1}-\bar x\right\| \dist(x^{k+1},X_\varphi^*) \le
	\frac12
	\dist(x^{k+1},X_\varphi^*).
	\]
	Moving this term to the left gives
	\begin{equation}\label{eq:lower-stationarity-next-distance}
		\begin{aligned}
			\dist(x^{k+1},X_\varphi^*)
			&\le
			2\kappa_mL_{f,2}
			\dist(x^k,X_\varphi^*)^2+
			2\kappa_m
			\Bigl[
			M_F
			+(L_F+\bar\alpha)c_d
			+\bar\alpha\bigl((1+M_H)c_r+\kappa_Hc_d\bigr)
			\Bigr]\sigma_k\\
			&\le 2\kappa_m
			L_{f,2}c_m^2\omega_k^4+
			2\kappa_m c_\sigma
			\Bigl(
			M_F
			+(L_F+\bar\alpha)c_d
			+\bar\alpha\bigl((1+M_H)c_r+\kappa_Hc_d\bigr)
			\Bigr)
			\omega_k^{2(1+\rho)}\\
			&\le 2\kappa_m
			\Bigl[
			L_{f,2}c_m^2+
			c_\sigma
			\Bigl(
			M_F
			+(L_F+\bar\alpha)c_d
			+\bar\alpha\bigl((1+M_H)c_r+\kappa_Hc_d\bigr)
			\Bigr)
			\Bigr]
			\omega_k^{2(1+\rho)}.
		\end{aligned}
	\end{equation}
	where the second inequality follows from Lemma~\ref{lem:lower-distance}, that is $\dist(x^k,X_\varphi^*)^2\le c_m^2G_k^4\le c_m^2\omega_k^4$ and the rule gives $\sigma_k\le c_\sigma\omega_k^{2(1+\rho)}$. The last inequality uses $\omega_k^4	\le	\omega_k^{2(1+\rho)}$ with \(\omega_k\le1\) and \(2(1+\rho)\le4\).
	
	Substituting \eqref{eq:lower-stationarity-next-distance}, $\dist(x^k,X_\varphi^*)^2 \le c_m^2\omega_k^{2(1+\rho)}$, and $\sigma_k\le c_\sigma\omega_k^{2(1+\rho)}$ into \eqref{eq:lower-stationarity-v-bound}, and using \(\omega_k^{4(1+\rho)}\le\omega_k^{2(1+\rho)}\), gives
	\[
	\left\|v^k\right\|
	\le
	c_R\omega_k^{2(1+\rho)},
	\]
	with the fixed constant
	\[
	c_R:=L_{f,2}c_m^2
	+4\kappa_m^2L_{f,2}\bigl(L_{f,2}c_m^2+c_\sigma c_v\bigr)^2
	+c_\sigma c_v,
	\qquad
	c_v:=M_F+(L_F+\bar\alpha)c_d
	+\bar\alpha\bigl((1+M_H)c_r+\kappa_Hc_d\bigr).
	\]
	Finally, since \(v^k\in\partial\varphi(x^{k+1})\),
	\[
	R_L(x^{k+1})
	\le
	\dist(0,\partial\varphi(x^{k+1}))
	\le
	\left\|v^k\right\|
	\le
	c_R\omega_k^{2(1+\rho)}.
	\]
	This proves \eqref{eq:lower_stationarity_order}.
\end{proof}

We now combine the preceding estimates to establish local invariance and superlinear convergence of the proximal Newton sequence.

\begin{theorem}\label{thm:local-recursion}
	Suppose that the standing assumptions and Assumption~\ref{asup:common-image} hold, and that Assumptions~\ref{asup:lower_regular}--\ref{asup_hessian_lipschitz} hold at \(\bar x\in S^*\). Let \(c_S>0\) be the constant from Lemma~\ref{lem:measure-distance-iteration}. Then, for any \(c_b\ge c_S\), there exist \(\epsilon_+,\delta_+,c_\omega>0\) such that the following holds. Given \(x^0\in\operatorname{dom}g\cap\mathbb B_{\epsilon_+}(\bar x)\), \(G_0\ge\sqrt{R_L(x^0)}\), and \(\tau_0>0\) satisfying $\omega_0=G_0+\tau_0\le\delta_+$, $\dist(x^0,S^*)\le c_b\omega_0$, let \(\{x^k\}_{k\ge0}\) be generated by the inexact proximal Newton method, with \eqref{eq:inexact_unit} and \eqref{eq:hessian-approx} satisfied at each iteration and the parameters selected according to \eqref{eq:parameter-rule}. Then
	\[
	x^k\in\mathbb B_{2\epsilon_+}(\bar x),
	\qquad
	\omega_{k+1}\le c_\omega\omega_k^{1+\rho},
	\qquad \forall\,k\ge0.
	\]
	Consequently, \(\{\omega_k\}\) converges \(Q\)-superlinearly to zero, while \(\dist(x^k,S^*)\) converges \(R\)-superlinearly to zero. If \(\rho=1\), the corresponding convergence rates are \(Q\)-quadratic and \(R\)-quadratic, respectively.
\end{theorem}

\begin{proof}
	Let $\epsilon_S,\delta_S,c_S$ be the constants of Lemma~\ref{lem:measure-distance-iteration} at $\bar x$, and fix $c_b\ge c_S$. The proximal correction gives $x^{k+1}\in\operatorname{dom}g$. Moreover, \eqref{eq:local_trial_measure} and \eqref{eq:parameter-rule}, together with the initialization, imply $G_k\ge\sqrt{R_L(x^k)}$ and $0<\tau_k\le\omega_k$ for every $k\ge0$. Let $\epsilon_d,\delta_d,c_d$ and $\epsilon_R,\delta_R,c_R$ be the constants of Lemmas~\ref{lem:step_estimate} and \ref{lem:lower-stationarity}, respectively, for this value of $c_b$, all at $\bar x$. Consider an index $k$ for which $x^k\in\mathbb B_{\epsilon_R}(\bar x)$, $\omega_k\le\delta_R$, and $\dist(x^k,S^*)\le c_b\omega_k$. By the constructions in Lemmas~\ref{lem:step_estimate} and \ref{lem:lower-stationarity}, $\epsilon_R\le\epsilon_d/2$ and $\delta_R\le\delta_d\le(\bar\alpha/c_\alpha)^{1/\rho}$. Hence both lemmas apply, $\left\|d^k\right\|\le c_d\omega_k$, $\omega_k\le1$, and $\alpha_k=c_\alpha\omega_k^\rho$. Also $\tau_k\le\omega_k$.
	
	First, Lemma~\ref{lem:exact-residual-bridge} gives
	\[
	\frac{\left\|\widetilde r^k\right\|}{\sigma_k}
	\le
	\alpha_k
	\bigl[
	(1+M_H)c_r\tau_k+\kappa_H\left\|d^k\right\|
	\bigr]\le
	c_\alpha\omega_k^\rho
	\bigl[
	(1+M_H)c_r\omega_k+\kappa_Hc_d\omega_k
	\bigr]=
	c_\alpha
	\bigl[
	(1+M_H)c_r+\kappa_Hc_d
	\bigr]
	\omega_k^{1+\rho}.
	\]
	Similarly,
	\[
	\alpha_k\left\|d^k\right\|
	\le
	c_\alpha c_d\omega_k^{1+\rho},\quad \left\|d^k\right\|^2
	\le
	c_d^2\omega_k^2
	\le
	c_d^2\omega_k^{1+\rho},
	\]
	where the last inequality follows from \(\omega_k\le1\) and \(0<\rho\le1\).
	
	By Lemma~\ref{lem:lower-stationarity}, $\sqrt{R_L(x^{k+1})}
	\le
	\sqrt{c_R}\,\omega_k^{1+\rho}$. Substituting these four estimates into \eqref{eq:local_trial_measure} gives
	\begin{equation}\label{eq:residual-estimate}
		\begin{aligned}
			G_{k+1}
			&=
			\frac{\left\|\widetilde r^k\right\|}{\sigma_k}
			+\left\|d^k\right\|^2
			+\alpha_k\left\|d^k\right\|
			+\sqrt{R_L(x^{k+1})}
			\\
			&\le
			\Bigl[
			c_\alpha\bigl((1+M_H)c_r+\kappa_Hc_d\bigr)
			+c_d^2
			+c_\alpha c_d
			+\sqrt{c_R}
			\Bigr]
			\omega_k^{1+\rho}
			\\
			&=
			c_+\omega_k^{1+\rho},
		\end{aligned}
	\end{equation}
	where $c_+ := c_\alpha\bigl((1+M_H)c_r+\kappa_Hc_d\bigr) +c_d^2 +c_\alpha c_d +\sqrt{c_R}$. Set $c_\omega:=c_++c_\tau$. By \eqref{eq:omega_def} and the accuracy update in \eqref{eq:parameter-rule},
	\[
	\omega_{k+1}=G_{k+1}+\tau_{k+1}
	\le(c_++c_\tau)\omega_k^{1+\rho}
	=c_\omega\omega_k^{1+\rho}.
	\]
	Choose
	\[
	\epsilon_+:=\min\{\epsilon_R/2,\epsilon_S/4\},\qquad
	\delta_+:=\min\left\{\delta_R,\delta_S,
	\frac{\epsilon_+}{2c_d},
	\left(\frac{1}{2c_\omega}\right)^{1/\rho}\right\}.
	\]
	Then $2\epsilon_+\le\epsilon_R$, $2\epsilon_+<\epsilon_S$, and $c_\omega\delta_+^\rho\le1/2$. For any local Newton update with $x^k\in\operatorname{dom}g\cap\mathbb B_{\epsilon_+}(\bar x)$, $G_k\ge\sqrt{R_L(x^k)}$, $\tau_k>0$, $\omega_k\le\delta_+$, and $\dist(x^k,S^*)\le c_b\omega_k$, the preceding bounds apply and give $\left\|x^{k+1}-\bar x\right\|<\epsilon_++c_d\delta_+\le2\epsilon_+$. Lemma~\ref{lem:measure-distance-iteration} therefore gives $\dist(x^{k+1},S^*)\le c_S\omega_{k+1}$.
	
	To ensure that the preceding local estimates remain applicable at every iteration, we prove by induction that, for all $k\ge0$, $x^k\in\mathbb B_{2\epsilon_+}(\bar x)$, $\omega_k\le2^{-k}\omega_0$, $\dist(x^k,S^*)\le c_b\omega_k$. The initial conditions give the assertion for $k=0$. Suppose that it holds through index $k$. Since $x^k\in\mathbb B_{\epsilon_R}(\bar x)$ and $\omega_k\le\delta_+\le\delta_R$, the preceding estimates yield $\left\|d^k\right\|\le c_d\omega_k$, \eqref{eq:residual-estimate}, and
	\[
	\omega_{k+1}\le c_\omega\omega_k^{1+\rho}
	\le\frac12\omega_k
	\le2^{-(k+1)}\omega_0.
	\]
	Consequently,
	\[
	\left\|x^{k+1}-\bar x\right\|
	\le\left\|x^0-\bar x\right\|+\sum_{j=0}^k\left\|d^j\right\|
	<\epsilon_++c_d\sum_{j=0}^k\omega_j
	\le\epsilon_++2c_d\omega_0
	\le2\epsilon_+.
	\]
	Thus $x^k,x^{k+1}\in\mathbb B_{\epsilon_S}(\bar x)$ and $\omega_k\le\delta_S$. Lemma~\ref{lem:measure-distance-iteration} gives
	\[
	\dist(x^{k+1},S^*)\le c_S\omega_{k+1}\le c_b\omega_{k+1},
	\]
	which completes the induction.
	
	Consequently, $\omega_k\to0$ and $\omega_{k+1}/\omega_k\le c_\omega\omega_k^\rho\to0$, showing that $\{\omega_k\}$ converges $Q$-superlinearly to zero. Since $\dist(x^k,S^*)\le c_b\omega_k$, the sequence $\dist(x^k,S^*)$ converges $R$-superlinearly to zero. When $\rho=1$, $\omega_{k+1}\le c_\omega\omega_k^2$ yields the corresponding quadratic rates.
\end{proof}
\section{Globalization and Convergence Analysis}\label{sec4}

Section~\ref{sec3} establishes the local superlinear convergence of the proximal Newton iteration when initialized sufficiently close to a solution. 
We now embed this iteration into a globalization strategy based on a trial-and-fallback mechanism, allowing the method to start from an arbitrary point in \(\operatorname{dom}g\). At each outer iteration, the algorithm first attempts a full proximal Newton step (the \emph{unit step}); if this trial fails to achieve a sufficient contraction in an optimality measure, the algorithm adjusts the penalty parameter and falls back to finding an approximate stationary point of a new penalized objective via an inner line-search procedure. Under our local assumptions at an accumulation point, we show that the fallback branch eventually becomes inactive, so that the method automatically reduces to the fast local iteration analyzed in Section~\ref{sec3}.

\paragraph{Unit-step trial.}
Specifically, at outer iteration $k$, given the current iterate $x^k$ and the parameters $(\omega_k,\sigma_k,\alpha_k,\tau_{k+1})$ updated according to \eqref{eq:local_trial_measure}--\eqref{eq:parameter-rule} in Section~\ref{sec3.2}, we first generate a candidate point following the unit-step construction of Section~\ref{sec:local-newton}. That is, we form the regularized proximal Newton subproblem \eqref{eq:q_model_unit} and its residual mapping $r_k$ in \eqref{eq:newton_residual} at $x^k$ using $(\sigma_k,\alpha_k)$, and compute an approximate model solution $\widetilde x^{k+1}\in\operatorname{dom}g$ satisfying the inexactness criterion \eqref{eq:inexact_unit}. As in \eqref{eq:unit_output}, we construct the search direction \(d^k\) using \(\widetilde x^{k+1}\). Prior to the acceptance test, we denote the Newton output in \eqref{eq:unit_update} by
\[
\widehat x^{k+1}:=x^k+d^k=\widetilde x^{k+1}-r_k(\widetilde x^{k+1}),
\]
to distinguish it from the final outer iterate \(x^{k+1}\). Evaluating the exact-Hessian proximal residual $\widetilde r^k$ in \eqref{eq:exact_hessian_residual} at the trial point $\widehat x^{k+1}$ (in place of $x^{k+1}$), we define the trial progress measure by
\begin{equation}\label{eq:trial_measure}
	G_k^{\mathrm t}
	:=
	\frac{\left\|\widetilde r^k\right\|}{\sigma_k}
	+\left\|d^k\right\|^2
	+\alpha_k\left\|d^k\right\|
	+\sqrt{R_L(\widehat x^{k+1})}.
\end{equation}
The unit step is accepted if it satisfies the contraction condition
\begin{equation}\label{eq:unit_acceptance}
	G_k^{\mathrm t}+\tau_{k+1}\le\theta(G_k+\tau_k),
\end{equation}
where $0<\theta<1$ is a fixed parameter. In this case, we set $x^{k+1}:=\widehat x^{k+1}$ and $G_{k+1}:=G_k^{\mathrm t}$, which completes outer iteration $k$ and advances the algorithm directly to iteration $k+1$.

\paragraph{Fallback strategy and outer update.}
If \eqref{eq:unit_acceptance} fails, the unit step is rejected. We then decrease the penalty parameter to the safeguarded value
\begin{equation}\label{eq:line_search_penalty}
	\widehat{\sigma}_k
	:=
	\min\{\sigma_k,c_\sigma\tau_k^2\},
\end{equation}
and switch to approximately minimizing the new penalized objective $\psi_{\widehat{\sigma}_k}$ starting from $x^{k,0}:=x^k$. Specifically, using the proximal residual of $\psi_\sigma$ defined by
\begin{equation*}
	R_\sigma(x):=\left\|x-\operatorname{Prox}_g
	\bigl(x-\nabla f(x)-\sigma\nabla F(x)\bigr)\right\|,
\end{equation*}
we generate a sequence of inner iterates $\{x^{k,j}\}_{j\ge 0}$ that monotonically decreases $\psi_{\widehat{\sigma}_k}$ until we obtain an iterate $x^{k,j}$ satisfying the approximate stationarity condition
\begin{equation}\label{eq:line_search_stop}
	\frac{R_{\widehat{\sigma}_k}(x^{k,j})}{\widehat{\sigma}_k}\le\tau_k.
\end{equation}
Once \eqref{eq:line_search_stop} is met, the auxiliary problem for $\psi_{\widehat{\sigma}_k}$ is deemed sufficiently solved. We then terminate the inner loop and set the new outer iterate $x^{k+1}$ to the proximal point already evaluated in the stopping test \eqref{eq:line_search_stop}:
\begin{equation}\label{eq:fallback_output}
	x^{k+1}:=\operatorname{Prox}_g\!
	\left(x^{k,j}-\nabla f(x^{k,j})
	-\widehat{\sigma}_k\nabla F(x^{k,j})\right).
\end{equation}
Combining both the unit-step and fallback branches, the outer progress measure is updated at the end of iteration $k$ according to
\begin{equation}\label{eq:G_update}
	G_{k+1}:=
	\begin{cases}
		G_k^{\mathrm t},&\text{if the unit step is accepted},\\[1mm]
		\displaystyle\frac{R_{\widehat{\sigma}_k}(x^{k,j})}{\widehat{\sigma}_k}
		+\sqrt{R_L(x^{k+1})},&\text{otherwise}.
	\end{cases}
\end{equation}
In particular, \(G_{k+1}=0\) implies \(x^{k+1}\in S^*\) in either branch. For an accepted unit step, this follows from \(d^k=0\) and \(\widetilde r^k=0\). For the fallback branch, \(R_{\widehat\sigma_k}(x^{k,j})=0\) yields \(x^{k+1}=x^{k,j}\). Together with \(R_L(x^{k+1})=0\), these conditions imply that \(x^{k+1}\) minimizes both \(\varphi\) and \(\psi_\sigma\) for some \(\sigma>0\), and hence belongs to \(S^*\).

\paragraph{Inner proximal Newton line-search procedure.}
It remains to specify how the inner procedure generates a point $x^{k,j}$ satisfying \eqref{eq:line_search_stop}. To preserve the second-order structure of the algorithm, we still apply an inexact regularized proximal Newton method, now equipped with a backtracking line search, to minimize $\psi_{\widehat{\sigma}_k}$. Starting from $x^{k,0}:=x^k$, the parameters $\widehat{\sigma}_k$, $\alpha_k$, and $\tau_k$ are kept fixed throughout the inner iterations $j=0,1,2,\ldots$. At each inner iterate $x^{k,j}$ violating \eqref{eq:line_search_stop}, we construct the Hessian approximation, regularization matrix, quadratic model, and residual mapping in the same manner as in Section~\ref{sec:local-newton}, with $(x^k,\sigma_k)$ replaced by $(x^{k,j},\widehat{\sigma}_k)$; with a slight abuse of notation, we retain the symbols $B_k,H_k,q_k$, and $r_k$ for these inner-loop objects. We then compute an approximate model point $\widetilde x^{k,j}\in\operatorname{dom}g$ satisfying the residual condition \eqref{eq:inexact_unit} (with $\widehat{\sigma}_k$ in place of $\sigma_k$), and additionally require that the corrected point $\widehat x^{k,j}:=\widetilde x^{k,j}-r_k(\widetilde x^{k,j})$ satisfy the model decrease condition
\begin{equation}\label{eq:inner_model_decrease}
	q_k(\widehat x^{k,j})\le q_k(x^{k,j}).
\end{equation}
To ensure that the resulting search direction $d^{k,j}:=\widehat x^{k,j}-x^{k,j}$ is a sufficient descent direction for $\psi_{\widehat{\sigma}_k}$, we choose the inexactness constant $c_r>0$ to satisfy
\begin{equation}\label{eq:line_search_forcing}
	c_r\bar\alpha
	\bigl(1+M_B+\bar\sigma\bar\alpha\bigr)<1,
\end{equation}
which guarantees $c_r\alpha_k\left\|I-H_k\right\|<1$ whenever the fallback branch is invoked. Finally, we perform a backtracking line search along $d^{k,j}$ to find the smallest nonnegative integer $m_{k,j}$ such that the step size $t_{k,j}:=\bar t\beta^{m_{k,j}}$ satisfies the Armijo condition
\begin{equation}\label{eq:line_search_test}
	\psi_{\widehat{\sigma}_k}\!\left(x^{k,j}+t_{k,j} d^{k,j}\right)
	\le\psi_{\widehat{\sigma}_k}(x^{k,j})-c_A\widehat{\sigma}_k\alpha_k t_{k,j} \left\|d^{k,j}\right\|^2,
\end{equation}
and set $x^{k,j+1}:=x^{k,j}+t_{k,j}d^{k,j}$. 

The complete method is summarized in Algorithm~\ref{alg:global-newton}.

	\algdef{SE}[SUBBLOCK]{BeginBlock}{EndBlock}[1]{\textbf{#1}}{}
	\algtext*{EndBlock}
	
	\begin{algorithm}[!ht]
		\caption{Simple Bilevel Newton Algorithm (SiB-Newton)}
		\label{alg:global-newton}
		\begin{algorithmic}[1]
			\Statex \textbf{Input:}
			Initial point \(x^0\in\operatorname{dom}g\), \(\tau_0>0\), \(G_0:=1+\sqrt{R_L(x^0)}\);
			\(\rho\in(0,1]\), \(\gamma,\theta,\beta\in(0,1)\), \(c_A\in(0,1/2)\), \(0<\bar t\le1\), and \(\bar\sigma,\bar\alpha,c_\sigma,c_\alpha,c_\tau,c_r>0\) satisfying \eqref{eq:line_search_forcing}.
			
			\For{\(k=0,1,2,\ldots\)}
			\If{\(G_k=0\)}
			\State \textbf{terminate.}
			\EndIf
			\State \textbf{Parameter update:} Set \((\omega_k,\sigma_k,\alpha_k,\tau_{k+1})\) by \eqref{eq:omega_def} and \eqref{eq:parameter-rule}.
			
			\BeginBlock{Unit-step trial:}
			\State Form \(q_k,r_k\) at \((x^k,\sigma_k,\alpha_k)\) via \eqref{eq:q_model_unit}--\eqref{eq:newton_residual}.
			\State Find \(\widetilde x^{k+1}\in\operatorname{dom}g\) satisfying \eqref{eq:inexact_unit}.
			\State Set \(\widehat x^{k+1}:=\widetilde x^{k+1}-r_k(\widetilde x^{k+1})\) and \(d^k:=\widehat x^{k+1}-x^k\).
			\State Evaluate \(\widetilde r^k\) at \(\widehat x^{k+1}\) by \eqref{eq:exact_hessian_residual} and compute \(G_k^{\mathrm t}\) by \eqref{eq:trial_measure}.
			\EndBlock
			
			\If{\(G_k^{\mathrm t}+\tau_{k+1}\le\theta(G_k+\tau_k)\)} \Comment{Unit step accepted}
			\State Set \(x^{k+1}:=\widehat x^{k+1}\) and \(G_{k+1}:=G_k^{\mathrm t}\).
			\Else \Comment{Fallback to inner line search}
			\State Set \(\widehat{\sigma}_k:=\min\{\sigma_k,c_\sigma\tau_k^2\}\) and obtain \(x^{k,j}\) via the \textbf{Inner Proximal Newton Line-Search Procedure} below.
			\State Set \(x^{k+1}\) by \eqref{eq:fallback_output} and \(G_{k+1}:=R_{\widehat{\sigma}_k}(x^{k,j})/\widehat{\sigma}_k+\sqrt{R_L(x^{k+1})}^{}\).
			\EndIf
			\EndFor
			
			\vspace{2mm}\hrule\vspace{2mm}
			\Statex \textbf{Inner Proximal Newton Line-Search Procedure} (minimizing \(\psi_{\widehat{\sigma}_k}\) with fixed \(\widehat{\sigma}_k,\alpha_k,\tau_k\)):
			\State Initialize \(x^{k,0}:=x^k\) and \(j:=0\).
			\While{\(R_{\widehat{\sigma}_k}(x^{k,j})/\widehat{\sigma}_k>\tau_k\)}
			\State Form \(B_k,H_k,q_k,r_k\) at \((x^{k,j},\widehat{\sigma}_k)\) with fixed \(\alpha_k,\tau_k\).
			\State Find \(\widetilde x^{k,j}\in\operatorname{dom}g\) satisfying \eqref{eq:inexact_unit} (with \(\widehat{\sigma}_k\)).
			\State Set \(\widehat x^{k,j}:=\widetilde x^{k,j}-r_k(\widetilde x^{k,j})\) satisfying \eqref{eq:inner_model_decrease} and \(d^{k,j}:=\widehat x^{k,j}-x^{k,j}\).
			\State Find the smallest integer \(m_{k,j}\ge 0\) such that \(t_{k,j}:=\bar t\beta^{m_{k,j}}\) satisfies \eqref{eq:line_search_test}.
			\State Set \(x^{k,j+1}:=x^{k,j}+t_{k,j}d^{k,j}\) and \(j:=j+1\).
			\EndWhile
			\State \textbf{Return} \(x^{k,j}\).
		\end{algorithmic}
	\end{algorithm}

\subsection{Auxiliary Results for Global Convergence}\label{sec:conv-prelim}

We now relate small algorithmic residuals and measures to feasibility and optimality for the original bilevel problem.

Since \(x^0\in\operatorname{dom}g\), the required model candidates and all proximal outputs lie in \(\operatorname{dom}g\). Accepted inner updates have finite penalized objective, and both outer updates select proximal outputs. Thus all algorithmic iterates and proximal trial points lie in \(\operatorname{dom}g\). The segments joining them also lie in \(\operatorname{dom}g\) by convexity.

For the subsequent analysis, set
\begin{equation}\label{eq:global-gradient-bound}
	L_s:=L_{hA}+\bar\sigma L_F.
\end{equation}
Then \(\nabla(f+\sigma F)\) is \(L_s\)-Lipschitz continuous on \(\operatorname{dom}g\), uniformly for \(0<\sigma\le\bar\sigma\).

The next lemma bounds the upper-level objective gap at a proximal point by the scaled penalized residual.

\begin{lemma}\label{lem1}
	Suppose that the standing assumptions hold. For \(u\in\operatorname{dom}g\) and \(0<\sigma\le\bar\sigma\), set $y:=\operatorname{Prox}_g\bigl(u-\nabla f(u)-\sigma\nabla F(u)\bigr)$. Then
	\begin{equation}\label{eq:proximal_certificate}
		w:=(u-y)+\nabla f(y)-\nabla f(u)
		+\sigma[\nabla F(y)-\nabla F(u)]\in\partial\psi_\sigma(y),
		\qquad
		\left\| w\right\|\le(1+L_s)R_\sigma(u),
	\end{equation}
	and, for every \(s\in X_\varphi^*\),
	\begin{equation}\label{eq:proximal-certificate-gap}
		F(y)-F(s)
		\le(1+L_s)\frac{R_\sigma(u)}{\sigma}\left\|y-s\right\|.
	\end{equation}
\end{lemma}

\begin{proof}
	Proximal optimality in \(y=\operatorname{Prox}_g(u-\nabla f(u)-\sigma\nabla
	F(u))\) gives \(u-\nabla f(u)-\sigma\nabla F(u)-y\in\partial g(y)\). Adding the smooth gradient of \(\psi_\sigma\) at \(y\) yields a penalized subgradient, $w\in\partial\psi_\sigma(y)$, and
	\[
	\left\| w\right\|\le(1+L_{hA}+\sigma L_F)\left\|u-y\right\|\le(1+L_s)R_\sigma(u),
	\]
	where the first inequality uses the standing Lipschitz bounds on the segment \([u,y]\subseteq\operatorname{dom}g\), and the second uses \(\sigma\le\bar\sigma\) and \(R_\sigma(u)=\left\|u-y\right\|\). The subgradient inequality at \(y\) and \(s\) gives
	\[
	\varphi(y)-\varphi^*+\sigma[F(y)-F(s)]
	=\psi_\sigma(y)-\psi_\sigma(s)\le\left\langle w,y-s\right\rangle.
	\]
	Since \(\varphi(y)\ge\varphi^*\), the Cauchy--Schwarz inequality and \eqref{eq:proximal_certificate} yield
	\[
	F(y)-F(s)
	\le\frac{\left\| w\right\|}{\sigma}\left\|y-s\right\|
	\le(1+L_s)\frac{R_\sigma(u)}{\sigma}\left\|y-s\right\|,
	\]
	which proves \eqref{eq:proximal-certificate-gap}.
\end{proof}

The next proposition bounds the distance of a returned line-search proximal point to $S^*$. This estimate supplies the compatibility condition needed to initialize the local iteration from a line-search output.
\begin{proposition}\label{prop:line-search-distance}
	Suppose that the standing assumptions hold. Fix \(\bar x\in S^*\) and suppose that Assumptions~\ref{asup:lower_regular} and \ref{asup_quadratic_growth} hold at \(\bar x\). Then there exist \(\epsilon_P,c_P>0\) such that, for every \(u\in\operatorname{dom}g\), \(0<\sigma\le\bar\sigma\), and
	\[
	y=\operatorname{Prox}_g(u-\nabla f(u)-\sigma\nabla F(u))
	\in\mathbb B_{\epsilon_P}(\bar x),
	\]
	one has
	\begin{equation}\label{eq:line-search-distance}
		\dist(y,S^*)\le c_P\left(\frac{R_\sigma(u)}{\sigma}+\sqrt{R_L(y)}\right).
	\end{equation}
\end{proposition}
\begin{proof}
	Let \(\epsilon_L,\epsilon_q\) be the fixed radii of Proposition~\ref{prop:lower-error-bound} and Assumption~\ref{asup_quadratic_growth}, and let \(c_m\) be the modulus of Proposition~\ref{prop:lower-error-bound}(i). Set the common radius and gradient bound
	\[
	\epsilon_P:=\min\{1,\epsilon_L,\epsilon_q\},
	\qquad
	M_P:=\left\|\nabla F(\bar x)\right\|+L_F.
	\]
	Since \(\epsilon_P\le1\) and \(\nabla F\) is \(L_F\)-Lipschitz on \(\operatorname{dom}g\), every \(v\in\mathbb B_{\epsilon_P}(\bar x)\cap\operatorname{dom}g\) satisfies \(\left\|\nabla F(v)\right\|\le\left\|\nabla F(\bar x)\right\|+L_F\left\|v-\bar x\right\|\le M_P\). Set \(p:=\operatorname{Proj}_{X_\varphi^*}(y)\) and \(s:=\operatorname{Proj}_{S^*}(p)\). Since the projections onto \(X_\varphi^*\) and \(S^*\) fix \(\bar x\), nonexpansiveness gives \(\left\|p-\bar x\right\|\le\left\|y-\bar x\right\|\le\epsilon_P\) and \(\left\|s-\bar x\right\|\le\left\|p-\bar x\right\|\), so \(p,s\) remain in \(\mathbb B_{\epsilon_P}(\bar x)\). The target distance splits into a normal and a tangential term:
	\begin{equation}\label{eq:prop4-distance-split}
		\dist(y,S^*)\le\left\|y-p\right\|+\left\|p-s\right\|
		=:h_y+\left\|p-s\right\|,
		\qquad h_y:=\dist(y,X_\varphi^*).
	\end{equation}
	
	It remains to estimate the tangential term \(\left\|p-s\right\|\). By Lemma~\ref{lem1} and the decomposition \(y-s=(y-p)+(p-s)\),
	\begin{equation}\label{eq:prop4-upper-gap}
		F(y)-F(s)
		\le
		(1+L_s)\frac{R_\sigma(u)}{\sigma}\left\|y-s\right\|
		\le
		(1+L_s)\frac{R_\sigma(u)}{\sigma}
		\bigl(\left\|p-s\right\|+h_y\bigr).
	\end{equation}
	
	Assumption~\ref{asup_quadratic_growth} at \(p\in X_\varphi^*\cap\mathbb B_{\epsilon_q}(\bar x)\), followed by the split \(F(p)-F(s)=F(y)-F(s)+F(p)-F(y)\), gives
	\begin{equation}\label{eq:prop4-tangential}
		\kappa_q\left\|p-s\right\|^2\le F(y)-F(s)+F(p)-F(y)\le
		(1+L_s)\frac{R_\sigma(u)}{\sigma}
		\bigl(\left\|p-s\right\|+h_y\bigr)
		+M_Ph_y,
	\end{equation}
	where the second line uses \eqref{eq:prop4-upper-gap} and the convexity of \(F\),
	\[
	F(p)-F(y)
	\le
	\left\langle \nabla F(p),p-y\right\rangle
	\le
	M_P\left\|p-y\right\|
	=
	M_Ph_y.
	\]
	
	Applying Young's inequality to the two products in \eqref{eq:prop4-tangential} gives
	\[
	(1+L_s)\frac{R_\sigma(u)}{\sigma}\left\|p-s\right\|
	\le
	\frac{\kappa_q}{2}\left\|p-s\right\|^2
	+
	\frac{(1+L_s)^2R_\sigma(u)^2}
	{2\kappa_q\sigma^2},\ (1+L_s)\frac{R_\sigma(u)}{\sigma}h_y
	\le
	\frac{(1+L_s)^2R_\sigma(u)^2}
	{2\kappa_q\sigma^2}
	+
	\frac{\kappa_q}{2}h_y^2.
	\]
	Substituting these bounds into \eqref{eq:prop4-tangential} and rearranging yields
	\[
	\left\|p-s\right\|
	\le
	\sqrt{
		\frac{2(1+L_s)^2R_\sigma(u)^2}
		{\kappa_q^2\sigma^2}
		+h_y^2
		+\frac{2M_P}{\kappa_q}h_y
	}\le
	\frac{\sqrt{2}(1+L_s)}{\kappa_q}
	\frac{R_\sigma(u)}{\sigma}
	+h_y
	+\sqrt{\frac{2M_P}{\kappa_q}}\sqrt{h_y}.
	\]
	Returning to \eqref{eq:prop4-distance-split}, and using \(h_y\le\sqrt{h_y}\) since \(h_y\le1\), gives
	\[
	\begin{aligned}
		\dist(y,S^*)
		&\le
		\frac{\sqrt{2}(1+L_s)}{\kappa_q}
		\frac{R_\sigma(u)}{\sigma}
		+2h_y
		+\sqrt{\frac{2M_P}{\kappa_q}}\sqrt{h_y}\\
		&\le
		\frac{\sqrt{2}(1+L_s)}{\kappa_q}
		\frac{R_\sigma(u)}{\sigma}
		+
		\left(
		2+\sqrt{\frac{2M_P}{\kappa_q}}
		\right)\sqrt{h_y}\\
		&\le
		\frac{\sqrt{2}(1+L_s)}{\kappa_q}
		\frac{R_\sigma(u)}{\sigma}
		+
		\left(
		2+\sqrt{\frac{2M_P}{\kappa_q}}
		\right)
		\sqrt{c_m}\sqrt{R_L(y)},
	\end{aligned}
	\]
	where the last inequality uses Proposition~\ref{prop:lower-error-bound}(i), namely, \(h_y\le c_mR_L(y)\). Thus \eqref{eq:line-search-distance} holds with $c_P:= \max\{\sqrt2(1+L_s)/\kappa_q,(2+\sqrt{2M_P/\kappa_q})\sqrt{c_m}\}$.
\end{proof}

\subsection{Global Convergence}

We now establish finite termination of the fallback procedure, convergence of the global measures to zero, and accumulation-point optimality. We first prove finite termination under the standing assumptions, and then establish the subsequent convergence results under the boundedness of the generated sequence. 

These results do not require the upper-level quadratic-growth condition in Assumption~\ref{asup_quadratic_growth} or the Hessian Lipschitz continuity in Assumption~\ref{asup_hessian_lipschitz}. The additional local assumptions will be imposed in the final theorem to establish eventual acceptance of Newton unit steps.

We begin with the finite termination of the line-search procedure.

\begin{lemma}\label{lem:line-search-termination}
	Suppose that the standing assumptions hold, \(0<c_A<1/2\), and the parameter condition \eqref{eq:line_search_forcing} holds. Then every invoked line-search procedure is well defined and terminates after finitely many inner iterations.
	For each fixed invoked outer iteration \(k\),
	\begin{equation}\label{eq:uniform-line-search-step}
		t_{k,j}\ge\underline t_k:=
		\min\left\{\bar t,\beta,
		\frac{\beta(1-2c_A)\widehat{\sigma}_k\alpha_k}{L_s}
		\right\}>0,
	\end{equation}
	where the lower bound is independent of \(j\).
\end{lemma}
\begin{proof}
	Fix an invoked outer iteration \(k\). The quantities \(\widehat{\sigma}_k,\alpha_k,\tau_k\) are positive and fixed throughout the inner loop.
	
	From the model-decrease condition \eqref{eq:inner_model_decrease} for the inner subproblem, we obtain
	\[
	0
	\ge
	q_k(\widehat x^{k,j})-q_k(x^{k,j})=
	\left\langle \nabla f(x^{k,j})+\widehat{\sigma}_k\nabla F(x^{k,j}),d^{k,j}\right\rangle
	+g(x^{k,j}+d^{k,j})-g(x^{k,j})
	+\frac12\left\langle d^{k,j},H_kd^{k,j}\right\rangle.
	\]
	Hence,
	\begin{equation}\label{eq:inner_model_descent}
		\begin{aligned}
			-\left\langle \nabla f(x^{k,j})+\widehat{\sigma}_k\nabla F(x^{k,j}),d^{k,j}\right\rangle-g(x^{k,j}+d^{k,j})+g(x^{k,j})\ge
			\frac12\left\langle d^{k,j},H_kd^{k,j}\right\rangle\ge
			\frac12\widehat{\sigma}_k\alpha_k\left\|d^{k,j}\right\|^2,
		\end{aligned}
	\end{equation}
	where the last inequality follows from \(H_k\succeq\widehat{\sigma}_k\alpha_k I\).
	
	We next relate the optimality residual \(R_{\widehat{\sigma}_k}(x^{k,j})\) to the step \(\left\|d^{k,j}\right\|\). Since \(\widetilde x^{k,j}=x^{k,j}+d^{k,j}+r_k(\widetilde x^{k,j})\), the correction formula \(\widehat x^{k,j}=\widetilde x^{k,j}-r_k(\widetilde x^{k,j})\) and the definition of the inner model residual give
	\[
	x^{k,j}+d^{k,j}=\operatorname{Prox}_g\!
	\bigl(x^{k,j}-\nabla f(x^{k,j})-\widehat{\sigma}_k\nabla F(x^{k,j})
	+(I-H_k)\bigl(d^{k,j}+r_k(\widetilde x^{k,j})\bigr)\bigr).
	\]
	Since \(R_{\widehat{\sigma}_k}(x^{k,j})=\left\|x^{k,j}-\operatorname{Prox}_g(x^{k,j}-\nabla f(x^{k,j})-\widehat{\sigma}_k\nabla F(x^{k,j}))\right\|\), the triangle inequality and nonexpansiveness of \(\operatorname{Prox}_g\) therefore give
	\begin{equation}\label{eq:inner_residual_bound}
		\begin{aligned}
			R_{\widehat{\sigma}_k}(x^{k,j})
			\le\left\|d^{k,j}\right\|
			+\left\|(I-H_k)\bigl(d^{k,j}+r_k(\widetilde x^{k,j})\bigr)\right\|\le(1+\left\|I-H_k\right\|)\left\|d^{k,j}\right\|
			+\left\|I-H_k\right\|\,\left\|r_k(\widetilde x^{k,j})\right\|,
		\end{aligned}
	\end{equation}
	where the spectral bound $\left\|I-H_k\right\|\le\max\{|1-\widehat{\sigma}_k\alpha_k|,|1-M_B-\widehat{\sigma}_k\alpha_k|\}\le 1+M_B+\bar\sigma\bar\alpha$ follows from \(\widehat{\sigma}_k\alpha_k I\preceq H_k\preceq(M_B+\widehat{\sigma}_k\alpha_k)I\); for the fixed outer iteration \(k\), this bound is independent of the inner index \(j\). Whenever the inner loop performs a Newton/line-search step at \(x^{k,j}\), the stopping test \eqref{eq:line_search_stop} fails, so \(R_{\widehat{\sigma}_k}(x^{k,j})>\widehat{\sigma}_k\tau_k\). The residual inexactness condition therefore gives
	\[
	\left\|r_k(\widetilde x^{k,j})\right\|
	\le c_r\widehat{\sigma}_k\alpha_k\tau_k
	<c_r\alpha_kR_{\widehat{\sigma}_k}(x^{k,j}).
	\]
	Substituting this bound into \eqref{eq:inner_residual_bound} and rearranging, which is possible since the spectral bound above and \eqref{eq:line_search_forcing} give \(c_r\alpha_k\left\|I-H_k\right\|<1\), yields
	\begin{equation}\label{eq:residual-direction-bound}
		R_{\widehat{\sigma}_k}(x^{k,j})
		\le\frac{1+\left\|I-H_k\right\|}{1-c_r\alpha_k\left\|I-H_k\right\|}\left\|d^{k,j}\right\|.
	\end{equation}
	
	Next, we show that the backtracking rule yields the step-size bound \eqref{eq:uniform-line-search-step}. By the standing assumptions, the gradient \(\nabla(f+\widehat{\sigma}_kF)\) is Lipschitz continuous on \(\operatorname{dom}g\) with the constant \(L_s\) of \eqref{eq:global-gradient-bound}. Both \(x^{k,j}\) and \(x^{k,j}+d^{k,j}\) belong to \(\operatorname{dom}g\), so their joining segment lies in \(\operatorname{dom}g\). Combining the \(L_s\)-Lipschitz continuity of \(\nabla(f+\widehat{\sigma}_kF)\), the convexity of \(g\), and the decrease estimate \eqref{eq:inner_model_descent}, we obtain for any \(t\in(0,1]\),
	\[
	\begin{aligned}
		&\quad \psi_{\widehat{\sigma}_k}(x^{k,j})
		-\psi_{\widehat{\sigma}_k}(x^{k,j}+td^{k,j})
		-c_A\widehat{\sigma}_k\alpha_k t\left\|d^{k,j}\right\|^2\\
		&\ge -t\bigl[
		\left\langle \nabla f(x^{k,j})+\widehat{\sigma}_k\nabla F(x^{k,j}),d^{k,j}\right\rangle
		+g(x^{k,j}+d^{k,j})-g(x^{k,j})\bigr]
		-\frac{L_s}{2}t^2\left\|d^{k,j}\right\|^2
		-c_A\widehat{\sigma}_k\alpha_k t\left\|d^{k,j}\right\|^2\\
		&\ge\frac12\widehat{\sigma}_k\alpha_k t\left\|d^{k,j}\right\|^2
		-\frac{L_s}{2}t^2\left\|d^{k,j}\right\|^2
		-c_A\widehat{\sigma}_k\alpha_k t\left\|d^{k,j}\right\|^2\\
		& = \frac t2\left\|d^{k,j}\right\|^2
		\bigl((1-2c_A)\widehat{\sigma}_k\alpha_k-L_s t\bigr).
	\end{aligned}
	\]
	Hence the Armijo condition \eqref{eq:line_search_test} holds whenever
	\begin{equation}\label{eq:armijo_threshold}
		0<t\le\min\{1,(1-2c_A)\widehat{\sigma}_k\alpha_k/L_s\}.
	\end{equation}
	Thus backtracking terminates. If its first trial is accepted, \(t_{k,j}=\bar t\). Otherwise the preceding rejected trial \(t_{k,j}/\beta\) exceeds the threshold in \eqref{eq:armijo_threshold}, proving \eqref{eq:uniform-line-search-step}. Substituting the accepted step into \eqref{eq:line_search_test} shows that each accepted inner update satisfies
	\begin{equation}\label{eq:line-search-descent}
		\psi_{\widehat{\sigma}_k}(x^{k,j+1})-\psi_{\widehat{\sigma}_k}(x^{k,j})
		\le-c_A\widehat{\sigma}_k\alpha_k t_{k,j}\left\|d^{k,j}\right\|^2.
	\end{equation}
	The exact model minimizer satisfies both inexact conditions, so an admissible corrected trial exists at every inner iteration. Accepted updates have finite objective value and hence remain in \(\operatorname{dom}g\).
	
	We now convert the direction decrease \eqref{eq:line-search-descent} into a residual decrease. Inverting \eqref{eq:residual-direction-bound} and using \(t_{k,j}\ge\underline t_k\) from \eqref{eq:uniform-line-search-step} gives
	\begin{equation}\label{eq:residual-sufficient-decrease}
		\psi_{\widehat{\sigma}_k}(x^{k,j+1})-\psi_{\widehat{\sigma}_k}(x^{k,j})
		\le
		-c_A\widehat{\sigma}_k\alpha_k\underline t_k
		\left(\frac{1-c_r\alpha_k\left\|I-H_k\right\|}{1+\left\|I-H_k\right\|}\right)^2
		R_{\widehat{\sigma}_k}(x^{k,j})^2.
	\end{equation}
	
	It remains to prove finite termination of the inner loop. While the loop continues, \(R_{\widehat{\sigma}_k}(x^{k,j})>\widehat{\sigma}_k\tau_k\), so \eqref{eq:residual-sufficient-decrease} gives
	\[
	\psi_{\widehat{\sigma}_k}(x^{k,j+1})-\psi_{\widehat{\sigma}_k}(x^{k,j})
	<
	-c_A(\widehat{\sigma}_k)^3\alpha_k\underline t_k\tau_k^2
	\left(\frac{1-c_r\alpha_k\left\|I-H_k\right\|}{1+\left\|I-H_k\right\|}\right)^2.
	\]
	The right-hand side is bounded above by a strictly negative constant independent of \(j\): for the fixed outer iteration \(k\), the quantities \(\widehat{\sigma}_k,\alpha_k,\tau_k>0\) are fixed, the spectral bound on \(\left\|I-H_k\right\|\) above is independent of \(j\), \(\underline t_k>0\) by \eqref{eq:uniform-line-search-step}, and \eqref{eq:line_search_forcing} bounds \(1-c_r\alpha_k\left\|I-H_k\right\|\) away from zero uniformly in \(j\). If the inner loop were infinite, summing this bound over \(j\) would give \(\psi_{\widehat{\sigma}_k}(x^{k,j})\to-\infty\), contradicting
	 \(\psi_{\widehat{\sigma}_k}\ge\varphi^*+\widehat{\sigma}_k\inf_{z\in\mathbb R^n}F(z)>-\infty\). Hence the inner loop terminates finitely.
\end{proof}

We next establish the convergence of \(\omega_k\) and \(G_k\) to zero.

\begin{proposition}\label{prop:global-measure}
	Suppose that the standing assumptions hold. Let \(\{x^k\}\) be an infinite sequence generated by the SiB-Newton method (Algorithm~\ref{alg:global-newton}). If \(\{x^k\}\) is bounded, then \(\omega_k\to0\) and \(G_k\to0\). 
\end{proposition}
\begin{proof}
	By the parameter rule \eqref{eq:parameter-rule}, \(0<\tau_k\le\gamma^k\tau_0\to0\), and hence \(\widehat{\sigma}_k\le c_\sigma\tau_k^2\to0\) by \eqref{eq:line_search_penalty}.
	
	We first analyze the line-search branch. Suppose that it is invoked at outer iteration \(k\), and let \(x^{k,j}\) be the final inner iterate. By the line-search measure in \eqref{eq:G_update} and the returned point \eqref{eq:fallback_output},
	\[
	G_{k+1}
	=\frac{R_{\widehat{\sigma}_k}(x^{k,j})}{\widehat{\sigma}_k}
	+\sqrt{R_L(x^{k+1})},
	\qquad
	x^{k+1}
	=\operatorname{Prox}_g\!
	\bigl(x^{k,j}-\nabla f(x^{k,j})
	-\widehat{\sigma}_k\nabla F(x^{k,j})\bigr).
	\]
	The first term is at most \(\tau_k\) by the stopping test \eqref{eq:line_search_stop}. It therefore remains only to control \(R_L(x^{k+1})\).
	
	Applying Lemma~\ref{lem1} with \(u=x^{k,j}\), \(y=x^{k+1}\), and \(\sigma=\widehat{\sigma}_k\), we obtain
	\[
	w^{k+1}
	:=
	x^{k,j}-x^{k+1}
	+\nabla f(x^{k+1})-\nabla f(x^{k,j})
	+\widehat{\sigma}_k
	\bigl[\nabla F(x^{k+1})-\nabla F(x^{k,j})\bigr]
	\in\partial\psi_{\widehat{\sigma}_k}(x^{k+1}),
	\]
	and
	\[
	\left\|w^{k+1}\right\|
	\le
	(1+L_s)R_{\widehat{\sigma}_k}(x^{k,j})
	\le
	(1+L_s)\widehat{\sigma}_k\tau_k,
	\]
	where the last inequality follows from the stopping test \eqref{eq:line_search_stop}.
	
	Since the sequence \(\{x^k\}\) is bounded and \(\nabla F\) is Lipschitz continuous on \({\rm dom}\,g\), there exists a constant \(\widehat M>0\) such that
	\[
	\|\nabla F(x^k)\|\le \widehat M
	\qquad \forall\,k\ge0.
	\]
	
	Moreover, $w^{k+1}
	-\widehat{\sigma}_k\nabla F(x^{k+1})
	\in\partial\varphi(x^{k+1})$, so \eqref{eq:residual_by_subgradient} gives
	\[
		R_L(x^{k+1})\le
		\left\|w^{k+1}\right\|
		+\widehat{\sigma}_k\left\|\nabla F(x^{k+1})\right\|\le
		\widehat{\sigma}_k
		\bigl[(1+L_s)\tau_k+\widehat M\bigr]\le
		c_\sigma
		\bigl[(1+L_s)\tau_0+\widehat M\bigr]\tau_k^2.
	\]
	Taking square roots and using the line-search measure, we obtain
	\[
	G_{k+1}=\frac{R_{\widehat{\sigma}_k}(x^{k,j})}{\widehat{\sigma}_k}+\sqrt{R_L(x^{k+1})}
	\le(1+\widehat c)\tau_k,\qquad\widehat c:=\sqrt{c_\sigma\bigl[(1+L_s)\tau_0+\widehat M\bigr]},
	\]
	for every outer iteration \(k\) at which the line-search branch is invoked.

	If instead the unit step is accepted at iteration \(k\), the acceptance rule \eqref{eq:unit_acceptance} and \(G_{k+1}=G_k^{\mathrm t}\) give \(\omega_{k+1}=G_k^{\mathrm t}+\tau_{k+1}\le\theta\omega_k\). Combining the two branches, we obtain, for all \(k\ge0\),
	\begin{equation}\label{eq:omega_recursion}
		\omega_{k+1}
		\le
		\max\bigl\{\theta\omega_k,\,
		(1+\widehat c)\tau_k+\tau_{k+1}\bigr\}.
	\end{equation}
	Since \((1+\widehat c)\tau_k+\tau_{k+1}\to0\), this sequence is bounded, and induction in \eqref{eq:omega_recursion} shows that \(\{\omega_k\}\) is bounded, so
	\[
	L:=\limsup_{k\to\infty}\omega_k<\infty.
	\]
	Taking the limit superior in \eqref{eq:omega_recursion} then yields $L\le\max\{\theta L,0\}=\theta L$, which, since \(0<\theta<1\) and \(L\ge0\), forces \(L=0\). Therefore \(\omega_k\to0\), and \(0\le G_k\le\omega_k\) gives \(G_k\to0\).
\end{proof}

The boundedness condition in Proposition~\ref{prop:global-measure} is guaranteed when the upper-level objective is coercive, as shown in the following lemma.

\begin{lemma}\label{lem:bounded-iterates}
	Suppose that the standing assumptions hold and that the upper-level objective \(F\) is coercive. Let \(\{x^k\}\) be an infinite sequence generated by the SiB-Newton method (Algorithm~\ref{alg:global-newton}). Then \(\{x^k\}\) is bounded.
\end{lemma}

\begin{proof}
	We first show that the line-search outputs are uniformly bounded. Fix any \(s\in X_\varphi^*\), which is nonempty by the standing assumptions. Suppose that the line-search branch is invoked at iteration \(k\), and let \(x^{k,j}\) be the final inner iterate. By Lemma~\ref{lem1} and the stopping condition~\eqref{eq:line_search_stop},
	\[
	F(x^{k+1})-F(s)\le(1+L_s)\tau_k\|x^{k+1}-s\|.
	\]
	By the coercivity of \(F\), choose \(R>0\) such that
	\[
	F(x)\ge F(s)+1\qquad\forall\,\|x-s\|\ge R.
	\]
	We claim that all sufficiently late line-search outputs satisfy \(\|x^{k+1}-s\|\le R\). Otherwise, for some sufficiently large \(k\) with \(\|x^{k+1}-s\|>R\), define \(z:=s+R(x^{k+1}-s)/\|x^{k+1}-s\|\). Then \(\|z-s\|=R\). By the convexity of \(F\) and the preceding estimate,
	\[
	F(z)-F(s)\le\frac{R}{\|x^{k+1}-s\|}\bigl[F(x^{k+1})-F(s)\bigr]\le R(1+L_s)\tau_k.
	\]
	Since \(\tau_k\to0\), the right-hand side is smaller than \(1\) for all sufficiently large \(k\), contradicting the choice of \(R\). Since there are only finitely many earlier outputs, there exists \(R_0>0\) such that $\|x^{k+1}-s\|\le R_0$ whenever the line-search branch is invoked.
	
	Next, we show that the progress measures associated with the line-search outputs are uniformly bounded. By the continuity of \(R_L\), there exists \(M_L\ge0\) such that \(R_L(x)\le M_L\) for all \(x\in\overline B_{R_0}(s)\). Using \eqref{eq:G_update} and the stopping condition~\eqref{eq:line_search_stop}, we obtain
	\[
	\omega_{k+1}=G_{k+1}+\tau_{k+1}\le\tau_k+\sqrt{R_L(x^{k+1})}+\tau_{k+1}\le2\tau_0+\sqrt{M_L}
	\]
	whenever the line-search branch is invoked. Consequently, setting $M_\omega:=\max\{\omega_0,\,2\tau_0+\sqrt{M_L}\}$, we have \(\omega_r\le M_\omega\) whenever \(r=0\) or \(x^r\) is a line-search output.
	
	It remains to control the iterates generated by consecutive accepted Newton unit steps. Fix any \(\ell\ge0\), and let \(r\le\ell\) be the largest index such that either \(r=0\) or \(x^r\) is a line-search output. Then all iterations \(r,\ldots,\ell-1\) accept Newton unit steps. By the acceptance condition~\eqref{eq:unit_acceptance},
	\[
	\omega_{r+j}\le\theta^j\omega_r,\qquad 0\le j\le\ell-r.
	\]
	Moreover, by the definition of the unit-step progress measure,
	\[
	\|x^{r+j+1}-x^{r+j}\|^2\le G_{r+j+1}\le\omega_{r+j+1}\le\theta^{j+1}\omega_r.
	\]
	Therefore,
	\[
		\|x^\ell-x^r\|\le\sum_{j=0}^{\ell-r-1}\|x^{r+j+1}-x^{r+j}\|\le\sqrt{\omega_r}\sum_{j=0}^{\ell-r-1}\theta^{(j+1)/2}
		\le\frac{\sqrt{\theta}}{1-\sqrt{\theta}}\sqrt{M_\omega}.
	\]
	Since \(x^r=x^0\) or \(x^r\) is a line-search output, we also have \(\|x^r-s\|\le\max\{\|x^0-s\|,R_0\}\). Thus
	\[
	\|x^\ell-s\|\le\max\{\|x^0-s\|,R_0\}+\frac{\sqrt{\theta}}{1-\sqrt{\theta}}\sqrt{M_\omega}
	\]
	for every \(\ell\ge0\). Hence \(\{x^k\}\) is bounded.
\end{proof}

We now use the vanishing global measures to characterize the accumulation points of the iterates. Lower metric subregularity is needed only locally for an accumulation point arising from unit-step outputs.

\begin{theorem}\label{thm:global-subseq}
	Suppose that the standing assumptions and Assumption~\ref{asup:common-image} hold.
	Let \(\{x^k\}\) be an infinite sequence generated by the SiB-Newton method (Algorithm~\ref{alg:global-newton}), and assume that \(\{x^k\}\) is bounded. Then every accumulation point \(\bar x\) of \(\{x^k\}\) belongs to \(X_\varphi^*\). Moreover, the following assertions hold:
	\begin{enumerate}
		\item[(i)] If \(\bar x\) is an accumulation point of the line-search outputs, then \(\bar x\in S^*\).
		\item[(ii)] If Assumption~\ref{asup:lower_regular} holds at \(\bar x\), then \(\bar x\in S^*\).
	\end{enumerate}
\end{theorem}
\begin{proof}
	By Proposition~\ref{prop:global-measure}, \(\omega_k\to0\) and \(G_k\to0\). We first establish lower feasibility and then upper optimality. The penalty cap is therefore eventually inactive. The outer measure \eqref{eq:G_update} satisfies \(G_k\ge\sqrt{R_L(x^k)}\) for every \(k\ge1\), so \(R_L(x^k)\to0\). Continuity of \(R_L\) and Proposition~\ref{prop:simple-blo-geometry} give \(R_L(\bar x)=0\) and \(\bar x\in X_\varphi^*\).
	
	If some subsequence converging to \(\bar x\) contains infinitely many line-search outputs, let \(x^{k_i}\to\bar x\) be such a subsequence and restrict to those indices, retaining the notation \(k_i\). Fix \(s\in S^*\), which is nonempty by the standing assumptions. Let \(j_i\) denote the final inner index of iteration \(k_i-1\), and set $\sigma_i:=\widehat{\sigma}_{k_i-1}$. By \eqref{eq:line_search_penalty}, \(0<\sigma_i\le\bar\sigma\). Moreover, \eqref{eq:G_update} and \eqref{eq:fallback_output} give
	\[
	x^{k_i}
	=
	\operatorname{Prox}_g\!\left(
	x^{k_i-1,j_i}
	-\nabla f(x^{k_i-1,j_i})
	-\sigma_i\nabla F(x^{k_i-1,j_i})
	\right),
	\qquad
	\frac{R_{\sigma_i}(x^{k_i-1,j_i})}{\sigma_i}
	\le G_{k_i}\to0.
	\]
	Applying Lemma~\ref{lem1} with \(u=x^{k_i-1,j_i}\), \(y=x^{k_i}\), and \(\sigma=\sigma_i\), and using the fixed \(s\in S^*\subseteq X_\varphi^*\), gives
	\[
	F(x^{k_i})-F(s)
	\le
	(1+L_s)
	\frac{R_{\sigma_i}(x^{k_i-1,j_i})}{\sigma_i}
	\left\|x^{k_i}-s\right\|
	\longrightarrow0,
	\]
	where the convergence follows from $R_{\sigma_i}(x^{k_i-1,j_i})/\sigma_i\le G_{k_i}\to0$ and the boundedness of \(\{x^{k_i}\}\). Since \(x^{k_i}\to\bar x\), continuity of \(F\) gives $F(\bar x)\le F(s)=F^*$. The lower-level feasibility established above gives \(\bar x\in X_\varphi^*\), and hence \(F(\bar x)\ge F^*\) by the definition of \(F^*\). Therefore $F(\bar x)=F^*$ and $\bar x\in S^*$.

	It remains to consider the case in which Assumption~\ref{asup:lower_regular} holds at \(\bar x\). Let \(x^{k_i}\to\bar x\) be a convergent subsequence; it contains only finitely many line-search outputs, since otherwise the preceding case applies. Discarding these finitely many indices, we may assume that every \(x^{k_i}\) is a unit-step output. We prove \(-\nabla F(\bar x)\in N_{X_\varphi^*}(\bar x)\) by constructing \(p_i\in X_\varphi^*\) and \(v_i\) such that \(p_i\to\bar x\), \(v_i\to-\nabla F(\bar x)\), and
	\[
	\limsup_i\langle v_i,s-p_i\rangle\le0
	\qquad\forall s\in X_\varphi^*.
	\]
	
	On this unit-step tail, \eqref{eq:local_trial_measure} and \eqref{eq:unit_acceptance} give \(G_{k_i}=G_{k_i-1}^{\mathrm t}\le\theta\omega_{k_i-1}\). Since each of the first three terms in \(G_{k_i-1}^{\rm t}\) is bounded by \(G_{k_i}\le\theta\omega_{k_i-1}\), we have
	\[
	\frac{\left\|\widetilde r^{k_i-1}\right\|}{\sigma_{k_i-1}}\to0,\qquad
	\left\|x^{k_i}-x^{k_i-1}\right\|\to0,\qquad
	\alpha_{k_i-1}\left\|x^{k_i}-x^{k_i-1}\right\|\to0.
	\]
	Since \(\sigma_{k_i-1}\le\bar\sigma\), also \(\widetilde r^{k_i-1}\to0\). Hence \(x^{k_i-1}\to\bar x\) and \(z_i:=x^{k_i}-\widetilde r^{k_i-1}\to\bar x\). Choose \(i_0\) so that, for \(i\ge i_0\), these points lie in \(\mathbb B_{\epsilon_L}(\bar x)\), the fixed ball on which Proposition~\ref{prop:lower-error-bound}(i) applies at \(\bar x\). Since \(\omega_{k_i-1}\to0\), choose the same \(i_0\) so that \(\sigma_{k_i-1}=c_\sigma\omega_{k_i-1}^{2(1+\rho)}\) throughout this tail.
	
	Define the lower Taylor remainder
	\[
	\begin{aligned}
		\rho_i&:=\nabla f(x^{k_i})-\nabla f(x^{k_i-1})
		-\nabla^2 f(x^{k_i-1})(x^{k_i}-x^{k_i-1})\\
		&=A^\top\bigl[\nabla h(Ax^{k_i})-\nabla h(Ax^{k_i-1})
		-\nabla^2h(Ax^{k_i-1})A(x^{k_i}-x^{k_i-1})\bigr]
		\in\range A^\top.
	\end{aligned}
	\]
	By \eqref{eq:exact_hessian_residual}, proximal optimality at \(z_i\) gives
	\[
	\widetilde r^{k_i-1}-\nabla f(x^{k_i-1})
	-\sigma_{k_i-1}\nabla F(x^{k_i-1})-\bigl[\nabla^2 f(x^{k_i-1})+\sigma_{k_i-1}\nabla^2F(x^{k_i-1})
	+\sigma_{k_i-1}\alpha_{k_i-1}I\bigr](x^{k_i}-x^{k_i-1})
	\in\partial g(z_i).
	\]
	Define the approximate normal vector directly by
	\[
	v_i:=\frac{\widetilde r^{k_i-1}}{\sigma_{k_i-1}}
	+\frac{\nabla f(z_i)-\nabla f(x^{k_i})}{\sigma_{k_i-1}}-\nabla F(x^{k_i-1})
	-\nabla^2F(x^{k_i-1})(x^{k_i}-x^{k_i-1})
	-\alpha_{k_i-1}(x^{k_i}-x^{k_i-1}).
	\]
	Adding \(\nabla f(z_i)\) to the proximal inclusion and using the definition of \(\rho_i\) gives the exact relation $\sigma_{k_i-1}v_i+\rho_i\in\partial\varphi(z_i)$. To prove convergence of \(v_i\), add \(\nabla F(x^{k_i})\) and apply the triangle inequality:
	\begin{equation*}
		\begin{aligned}
			\left\|v_i+\nabla F(x^{k_i})\right\|&\le\left\|\frac{\widetilde r^{k_i-1}}{\sigma_{k_i-1}}
				+\frac{\nabla f(z_i)-\nabla f(x^{k_i})}{\sigma_{k_i-1}}\right\|+\left\|\nabla F(x^{k_i})-\nabla F(x^{k_i-1})
				-\nabla^2F(x^{k_i-1})(x^{k_i}-x^{k_i-1})\right\|\\
			&\quad +\alpha_{k_i-1}\left\|x^{k_i}-x^{k_i-1}\right\|.
		\end{aligned}
	\end{equation*}
	The first term is at most \((1+L_{hA})\left\|\widetilde r^{k_i-1}\right\|/\sigma_{k_i-1}\to0\) by the standing Lipschitz bound for \(\nabla f\), the identity \(z_i-x^{k_i}=-\widetilde r^{k_i-1}\), and the unit-step smallness above; the last term tends to zero by the same smallness. Since \(F\in C^2\) and \(x^{k_i-1},x^{k_i}\to\bar x\), the second term, the gradient Taylor remainder of \(F\) with base point \(x^{k_i-1}\), also tends to zero. Hence $v_i\to-\nabla F(\bar x)$.
	
	It remains to show that \(v_i\) asymptotically satisfies the normal inequality for \(X_\varphi^*\). Set $p_i:=\operatorname{Proj}_{X_\varphi^*}(z_i)$ and fix arbitrary \(s\in X_\varphi^*\). We start from the target quantity \(\left\langle v_i,s-p_i\right\rangle\). Since $\sigma_{k_i-1}v_i+\rho_i\in\partial\varphi(z_i)$, the subgradient inequality gives
	\[
	\left\langle \sigma_{k_i-1}v_i+\rho_i,s-z_i\right\rangle
	\le \varphi^*-\varphi(z_i)\le0.
	\]
	Moreover, Remark~\ref{rem:fixed-kernel} and \(\rho_i\in\operatorname{range}A^\top\) give \(\left\langle \rho_i,s-p_i\right\rangle=0\). Splitting \(s-p_i=(s-z_i)+(z_i-p_i)\) and using the Cauchy--Schwarz inequality therefore yields
	\begin{equation}\label{eq:unit-step-normal-inequality}
		\left\langle v_i,s-p_i\right\rangle=\frac{\left\langle \sigma_{k_i-1}v_i+\rho_i,s-p_i\right\rangle}{\sigma_{k_i-1}}\le\frac{\left\langle \sigma_{k_i-1}v_i+\rho_i,z_i-p_i\right\rangle}{\sigma_{k_i-1}}\le\left(\left\|v_i\right\|+\frac{\left\|\rho_i\right\|}{\sigma_{k_i-1}}\right)
		\dist(z_i,X_\varphi^*).
	\end{equation}
	Since \(v_i\to-\nabla F(\bar x)\), the sequence \(\{v_i\}\) is bounded. Also, \(z_i\to\bar x\in X_\varphi^*\), so $\left\|v_i\right\|\dist(z_i,X_\varphi^*)\to0$. It therefore remains only to prove \(\left\|\rho_i\right\|\dist(z_i,X_\varphi^*)/\sigma_{k_i-1}\to0\).
	
	On indices with \(A(x^{k_i}-x^{k_i-1})\ne0\), we write
	\[
	\frac{\left\|\rho_i\right\|}{\sigma_{k_i-1}}
	\dist(z_i,X_\varphi^*)
	=
	\frac{\left\|\rho_i\right\|}
	{\left\|A(x^{k_i}-x^{k_i-1})\right\|}
	\frac{\left\|A(x^{k_i}-x^{k_i-1})\right\|}
	{\sigma_{k_i-1}}
	\dist(z_i,X_\varphi^*).
	\]
	Since \(h\in C^2\) and \(Ax^{k_i-1},Ax^{k_i}\to A\bar x\), the definition of \(\rho_i\) gives
	\[
	\frac{\left\|\rho_i\right\|}
	{\left\|A(x^{k_i}-x^{k_i-1})\right\|}
	\le
	\left\| A\right\|\,
	\frac{
		\left\|
			\nabla h(Ax^{k_i})-\nabla h(Ax^{k_i-1})
			-\nabla^2h(Ax^{k_i-1})A(x^{k_i}-x^{k_i-1})
		\right\|
	}{
		\left\|A(x^{k_i}-x^{k_i-1})\right\|
	}
	\longrightarrow0.
	\]
	
	With the fixed modulus \(c_m\) of Proposition~\ref{prop:lower-error-bound}(i), the residual estimate, the measure bounds, and unit-step acceptance give
	\[
	\dist(x^{k_i-1},X_\varphi^*)\le c_mR_L(x^{k_i-1})
	\le c_mG_{k_i-1}^2\le c_m\omega_{k_i-1}^2,\ \dist(x^{k_i},X_\varphi^*)\le c_m\theta^2\omega_{k_i-1}^2.
	\]
	Hence Proposition~\ref{prop:simple-blo-geometry} yields
	\begin{equation*}
		\left\|A(x^{k_i}-x^{k_i-1})\right\|
		\le\left\| A\right\|\bigl[\dist(x^{k_i},X_\varphi^*)
		+\dist(x^{k_i-1},X_\varphi^*)\bigr]\le\left\| A\right\|\,c_m(1+\theta^2)\omega_{k_i-1}^2.
	\end{equation*}
	Moreover, \(\left\|\widetilde r^{k_i-1}\right\|/\sigma_{k_i-1} \le G_{k_i}\le\theta\omega_{k_i-1}\), and, for \(i\ge i_0\), $\sigma_{k_i-1} =c_\sigma\omega_{k_i-1}^{2(1+\rho)}$, and $\omega_{k_i-1}\le1$. Since \(z_i=x^{k_i}-\widetilde r^{k_i-1}\), we therefore have
	\[
	\dist(z_i,X_\varphi^*)\le\dist(x^{k_i},X_\varphi^*)
	+\left\|\widetilde r^{k_i-1}\right\|\le c_m\theta^2\omega_{k_i-1}^2
	+c_\sigma\theta\omega_{k_i-1}^{3+2\rho}\le (c_m\theta^2+c_\sigma\theta)
	\omega_{k_i-1}^2.
	\]
	Combining these estimates gives
	\[
	\frac{\left\|A(x^{k_i}-x^{k_i-1})\right\|}
	{\sigma_{k_i-1}}
	\dist(z_i,X_\varphi^*)
	\le
	\frac{
		\left\| A\right\|\,c_m(1+\theta^2)
		(c_m\theta^2+c_\sigma\theta)
	}{c_\sigma}
	\omega_{k_i-1}^{2-2\rho}\le
	\frac{
		\left\| A\right\|\,c_m(1+\theta^2)
		(c_m\theta^2+c_\sigma\theta)
	}{c_\sigma},
	\]
	where the last inequality follows from \(0<\rho\le1\) and \(\omega_{k_i-1}\le1\). Thus the second factor is bounded.
	
	Consequently, the product tends to zero on indices with \(A(x^{k_i}-x^{k_i-1})\ne0\). On the remaining indices the two image points coincide, so \(\rho_i=0\) and the product vanishes. Hence
	\[
	\frac{\left\|\rho_i\right\|}{\sigma_{k_i-1}}\dist(z_i,X_\varphi^*)\to0.
	\]
	
	The right-hand side of \eqref{eq:unit-step-normal-inequality} tends to zero by the limit just proved and the boundedness of \(\{v_i\}\), so \(\limsup_i\left\langle v_i,s-p_i\right\rangle\le0\). Since \(\bar x\in X_\varphi^*\), nonexpansiveness of the projection gives \(\left\|p_i-\bar x\right\|\le\left\|z_i-\bar x\right\|\to0\). Together with \(v_i\to-\nabla F(\bar x)\), this shows that \(\left\langle v_i,s-p_i\right\rangle\) converges to \(\left\langle -\nabla F(\bar x),s-\bar x\right\rangle\), and hence
	\[
	\left\langle -\nabla F(\bar x),s-\bar x\right\rangle\le0
	\qquad\forall\ s\in X_\varphi^*.
	\]
	Thus \(0\in\nabla F(\bar x)+N_{X_\varphi^*}(\bar x)\). Convexity of \(F\) and closed convexity of \(X_\varphi^*\) imply \(\bar x\in S^*\).
\end{proof}

We next show that the globalization mechanism eventually accepts every Newton unit step, so that Theorem~\ref{thm:local-recursion} applies to the tail sequence.
\begin{theorem}\label{thm:local-superlinear}
	Suppose that the standing assumptions and Assumption~\ref{asup:common-image} hold. Let \(\{x^k\}\) be an infinite sequence generated by the SiB-Newton method (Algorithm~\ref{alg:global-newton}), and let \(\bar x\) be an accumulation point of \(\{x^k\}\). If \(\{x^k\}\) is bounded, Assumptions~\ref{asup:lower_regular}--\ref{asup_hessian_lipschitz} hold at \(\bar x\), and the  Hessian approximations satisfy \eqref{eq:hessian-approx}, then \(\{\omega_k\}\) converges \(Q\)-superlinearly to zero, while \(\{\dist(x^k,S^*)\}\) and \(\{\left\|x^k-\bar x\right\|\}\) converge \(R\)-superlinearly to zero. If \(\rho=1\), the corresponding rates are \(Q\)-quadratic and \(R\)-quadratic, respectively.
\end{theorem}

\begin{proof}
	Proposition~\ref{prop:global-measure} gives \(G_k,\omega_k\to0\), and Theorem~\ref{thm:global-subseq} yields \(\bar x\in S^*\). Let \(\epsilon_S,\delta_S,c_S\) be the constants of Lemma~\ref{lem:measure-distance-iteration} and \(\epsilon_P,c_P\) those of Proposition~\ref{prop:line-search-distance}, all at \(\bar x\), and set \(c_b:=\max\{c_S,c_P\}\). By Theorem \ref{thm:local-recursion}, let \(\epsilon_+,\delta_+,c_\omega\) be the local constants corresponding to this choice of \(c_b\). Choose \(\epsilon,\delta>0\) such that
	\[
	\epsilon\le\min\{\epsilon_+,\epsilon_S/2,\epsilon_P\},
	\qquad
	\delta\le\min\{\delta_+,\delta_S,\epsilon^2\},
	\qquad
	c_\omega\delta^\rho\le\theta/2.
	\]
	Since \(\omega_k\to0\), choose \(K_0\) such that \(\omega_k\le\delta\) for all \(k\ge K_0\). Take a subsequence \(x^{k_i}\to\bar x\), and choose \(K:=k_i\ge K_0+1\) sufficiently large that \(x^K\in\mathbb B_\epsilon(\bar x)\). In particular, \(\omega_K\le\delta\).
	
	To verify the remaining compatibility condition, consider the branch producing \(x^K\). If it was an accepted unit step, then \eqref{eq:trial_measure} and \eqref{eq:G_update} give \(\left\|x^K-x^{K-1}\right\|^2\le G_K\le\omega_K\), so
	\[
	\left\|x^{K-1}-\bar x\right\|
	\le\left\|x^K-\bar x\right\|+\sqrt{\omega_K}
	<\epsilon+\sqrt\delta
	\le2\epsilon\le\epsilon_S.
	\]
	Both model points therefore lie in \(\mathbb B_{\epsilon_S}(\bar x)\), and \(\omega_{K-1}\le\delta_S\). Lemma~\ref{lem:measure-distance-iteration} gives \(\dist(x^K,S^*)\le c_S\omega_K\). If \(x^K\) was returned by the line-search branch, Proposition~\ref{prop:line-search-distance}, applied to \eqref{eq:fallback_output}, and the measure \eqref{eq:G_update} give \(\dist(x^K,S^*)\le c_PG_K\le c_P\omega_K\). Thus in either case,
	\begin{equation}\label{eq:entry-distance}
	\dist(x^K,S^*)\le c_b\omega_K.
	\end{equation}
	Moreover, \eqref{eq:G_update} gives \(G_k\ge\sqrt{R_L(x^k)}\) for either branch, and \(\tau_k>0\). Hence the entry state \((x^K,G_K,\tau_K)\) satisfies all initialization conditions of Theorem~\ref{thm:local-recursion}; the same residual bound supplies the input of Lemma~\ref{lem:lower-distance} whenever a local estimate is applied.
	
	Suppose that \(m\ge K\) is the first iteration at which the unit step is rejected. The accepted updates at indices \(K,\ldots,m-1\) follow the local iteration. The finite induction established in the proof of Theorem~\ref{thm:local-recursion}, initialized at \(K\), therefore gives
	\[
	x^m\in\mathbb B_{2\epsilon_+}(\bar x),
	\qquad
	\omega_m\le\omega_K\le\delta,
	\qquad
	\dist(x^m,S^*)\le c_b\omega_m.
	\]
	The local radii and thresholds in that construction ensure that Lemmas~\ref{lem:exact-residual-bridge}, \ref{lem:step_estimate}, and \ref{lem:lower-stationarity} apply to the Newton trial at \(m\). Their one-step bound \eqref{eq:residual-estimate}, established in the proof of Theorem~\ref{thm:local-recursion}, together with \eqref{eq:parameter-rule}, yields
	\[
	G_m^{\mathrm t}+\tau_{m+1}
	\le c_\omega\omega_m^{1+\rho}
	\le c_\omega\delta^\rho\omega_m
	\le\frac\theta2\omega_m
	<\theta\omega_m.
	\]
	This contradicts rejection under \eqref{eq:unit_acceptance}. Thus every unit step from iteration \(K\) onward is accepted.
	
	The tail now satisfies the local iteration and all hypotheses of Theorem~\ref{thm:local-recursion}. Applying that theorem with initial state \((x^K,G_K,\tau_K)\) gives
	\[
	\omega_{k+1}\le c_\omega\omega_k^{1+\rho},
	\qquad
	\dist(x^k,S^*)\le c_b\omega_k
	\qquad \forall\ k\ge K.
	\]
	Moreover, the step estimate and geometric decay yield
	\[
	\sum_{j=k}^{\infty}\|d^j\|
	\le c_d\sum_{j=k}^{\infty}\omega_j
	\le 2c_d\omega_k
	\qquad \forall\,k\ge K.
	\]
	Hence \(\{x^k\}\) converges to some \(x^\infty\), with $\|x^k-x^\infty\|\le 2c_d\omega_k$. Since the tail contains arbitrarily late elements of the subsequence \(x^{k_i}\to\bar x\), we have \(x^\infty=\bar x\). The asserted convergence rates follow.
	
\end{proof}

    \section{Numerical Experiments}\label{sec5}
In this section, we provide a numerical illustration of the convergence behavior of SiB-Newton, and compare it with several first-order baselines on a linear regression problem. All experiments were implemented in MATLAB and conducted on a laptop with an Intel  i7-1260P CPU (2.10 GHz) and 32 GB RAM.

We consider the following simple bilevel problem:
\begin{equation}\label{LP}
	\min_{x\in\mathbb{R}^n}\ \frac{1}{2}\|A_1 x - b_1\|^2
	\quad\text{s.t.}\quad
	x \in \operatorname*{argmin}_{x\in\mathbb{R}^n}\ \frac{1}{2}\|A_2 x - b_2\|^2,
\end{equation}
where $A_1\in\mathbb{R}^{m_1\times n}$, $A_2\in\mathbb{R}^{m_2\times n}$, $b_1\in\mathbb{R}^{m_1}$ and $b_2\in\mathbb{R}^{m_2}$. The lower-level problem is an unconstrained least-squares problem, whose solution set is
\[
X_\varphi^* = \{x \in \mathbb{R}^n : A_2 x = b_2\}.
\]

The data are generated as follows. The design matrices \(A_1\in\mathbb{R}^{m_1\times n}\) and \(A_2\in\mathbb{R}^{m_2\times n}\) are generated independently from the standard normal distribution \(\mathcal{N}(0,1)\). The coefficient vector \(x^*\in\mathbb{R}^n\) is constructed to be \(5\%\)-sparse: we randomly select \(k=\max\{1,\lfloor 0.05\,n \rfloor\}\) indices, assign those entries with i.i.d.\ samples from \(\mathcal{N}(0,1)\), and set all other components to zero. The lower-level response is then set as \(b_2 := A_2 x^*\), ensuring that $x^*\in X_{\varphi}^*$. To guarantee that $x^*$ is also the unique upper-level solution, we proceed as follows. We take a QR decomposition of \(A_1^\top\), and denote \(P_1=QQ^\top\), \(M=(I-P_1)A_2^\top\), and pick a nonzero vector \(\lambda\in\ker M\). Solve \(A_1^\top y=A_2^\top\lambda\) and define
\[
b_1 := A_1 x^* + y,
\qquad
A_1^\top\bigl(b_1 - A_1 x^*\bigr)=A_2^\top\lambda.
\]
This construction ensures that the bilevel structure in~\eqref{LP} is consistent and \(x^*\) is well defined. 

We compare the proposed method with several representative algorithms, including BiG-SAM~\cite{sabach2017first}, DBGD~\cite{gong2021bi}, a-IRG~\cite{kaushik2021method}, BiSG~\cite{merchav2023convex}, SEA~\cite{shen2023online}, AGM-BiO~\cite{cao2024accelerated}, R-APM~\cite{samadi2024achieving}, PB-APG\cite{chen2024penalty}, Bisec\cite{wang2024near}, R-VFISTA~\cite{samadi2025iteratively}, IRE-APG\cite{shtern2026convergence}, FBi-PG\cite{merchav2026dynamic}, and Algorithm~2 of Bo\c{t} et al.~\cite{boct2026accelerating}. 

We evaluate their performance using two criteria: the solution error \(\|x^k-x^*\|\) and the KKT residual \(R_{\mathrm{KKT}}(x^k)\). Since the lower-level least-squares problem is consistent, problem~\eqref{LP} is equivalent to minimizing the upper-level objective subject to the linear constraint \(A_2x=b_2\). The corresponding KKT conditions are
\[
A_2x=b_2,\ A_1^\top(A_1x-b_1)+A_2^\top\lambda=0
\]
for some multiplier \(\lambda\in\mathbb R^{m_2}\). Eliminating the multiplier through orthogonal projection, we define the KKT residual by
\[
	R_{\mathrm{KKT}}(x):=\sqrt{\|A_2x-b_2\|^2+\|(I-P_2)A_1^\top(A_1x-b_1)\|^2},
\]
where $P_2$ denotes the orthogonal projector onto $\operatorname{Range}(A_2^\top)$. For plotting purposes, each error sequence is replaced by its running minimum. All methods start from $x_0=(0,0,\cdots,0)$, and are terminated upon meeting the tolerance $\text{tol}=10^{-8}$, the maximum iteration limit $\text{maxiter}=10^6$, or the time limit of $\text{maxtime}=1000$s.

We tune all baselines by grid search on predefined ranges and select the hyperparameters that minimize the distance to the optimal solution based on the independently generated instance with random seed 1. Concretely, for BiG\text{-}SAM we search the scaling factors $(\eta,\gamma)$ over $\{1/2,1/4,1/8\}^2$ (applied to the stepsizes $s=\eta\cdot(2/\lambda_{\max}(A_2^\top A_2))$ and $t=\gamma\cdot(1/\lambda_{\max}(A_1^\top A_1))$) and the averaging parameter $b$ over $\{10^{-3},10^{-2},10^{-1},1,10\}$ and select $(\eta,\gamma,b)=(1/4,1/2,10)$; for DBGD we search the parameters of the dynamic barrier function and the stepsize $(\alpha,\beta,\epsilon)$ over $\{10^{-3},10^{-2},10^{-1},1,10\}^3$ and select $(\alpha,\beta,\epsilon)=(1,10^{-2},10^{-3})$; for a\text{-}IRG we search the stepsize-related coefficients $(\eta,\gamma)$ over $\{10^{-3},10^{-2},10^{-1},1,10\}^2$ and $b$ over $\{0.05,0.1,0.2,0.3,0.4\}$ and select $(\eta,\gamma,b)=(0.1,10^{-2},0.2)$; for BiSG we search the stepsize-related coefficient $c$ over $\{10^{-3},10^{-2},10^{-1},1,10\}$ and select $10^{-1}$; for SEA we search the primal and dual step sizes and the Slater constant $(\gamma,\beta,r)$ over $\{10^{-5},10^{-4},10^{-3},10^{-2},10^{-1},1\}^3$ and select $(\gamma,\beta,r)=(10^{-5},10^{-5},10^{-2})$; for AGM-BiO we search the stepsize-related coefficient $\gamma$ over $\{10^{-3},10^{-2},10^{-1},1,10\}$ and select $10^{-3}$; for R\text{-}APM we search the stepsizes $(\eta,\gamma)$ over $\{10^{-6},10^{-5},10^{-4},10^{-3}\}^2$ and select $(\eta,\gamma)=(10^{-5},10^{-3})$; for PB\text{-}APG we search the stepsize-related coefficient $\eta$ over $\{10^{-6},2\times10^{-6},5\times10^{-6},10^{-5},2\times10^{-5},5\times10^{-5},10^{-4},10^{-3}\}$ and select $\eta=10^{-5}$; for R\text{-}VFISTA we search the stepsize and the regularization coefficient $(\eta,\mu_f)$ over $\{10^{-5},10^{-4},10^{-3},10^{-2},10^{-1},1\}^2$ and select $(\eta,\mu_f)=(1,10^{-5})$; for IRE\text{-}APG we search the coefficient $\beta$ over $\{0.1,0.25,0.5,0.75,1,1.25,1.5,1.75,2\}$ and select $\beta=1.25$; for FBi\text{-}PG we fix $a=2$, search the coefficient $\gamma$ over $\{0.5,1,1.3,1.5,1.9,2,3\}$, and select $\gamma=1$; for Bisec we set the accuracy parameters $(\epsilon_f,\epsilon_g)=(10^{-8},10^{-8})$; for Algorithm~2 of Bo\c{t} et al.~\cite{boct2026accelerating} we search the acceleration-related coefficients $(\alpha,\gamma)\in\{3.01,3.1,4,13\}\times\{10^{-2},10^{-1},1,10\}$, the stepsize-related coefficient $s_{\rm frac}\in\{10^{-2},10^{-1},1\}$, and the regularization-related coefficients $(c,\beta,\delta)\in\{10^{-2},10^{-1},1,10\}^3$, and select $(\alpha,\gamma,s_{\rm frac},c,\beta,\delta)=(13,10,1,10^{-1},10,1)$. Here $s=s_{\rm frac}/\lambda_{\max}(A_2^\top A_2)$, and hence the selected value $s_{\rm frac}=1$ corresponds to $s=1/\lambda_{\max}(A_2^\top A_2)$.

SiB-Newton uses the following hyperparameters: $\tau_0=1,\ \rho=0.5,\ \gamma_\tau=0.95,\ \theta=0.85,\ \beta=0.5,\ \bar{\sigma}=1,\ \bar{\alpha}=1,\ \bar{t}=1,\ c_\sigma=0.15,\ c_\alpha=1,\ c_\tau=1,\ c_A=10^{-4}$. The parameter $c_r$ is automatically chosen as $c_r=0.25/[\bar{\alpha}(1+M_B+\bar{\sigma}\bar{\alpha})]$, where $M_B=\|A_2^\top A_2\|+\bar{\sigma}\|A_1^\top A_1\|$. At each iteration, solving the quadratic subproblems in \eqref{eq:q_model_unit} reduces to solving the symmetric positive definite linear system
\[
(\sigma_k A_1^\top A_1+A_2^\top A_2+\sigma_k\alpha_k I)x=\sigma_k A_1^\top b_1+A_2^\top b_2+\sigma_k\alpha_k x^k.
\]
The linear systems are solved by MATLAB's PCG routine, where the relative tolerance is adaptively determined according to the inexact model condition and restricted to the interval $[10^{-14},10^{-1}]$. If the prescribed residual accuracy cannot be achieved due to numerical errors, a direct linear solve is used as a safeguard. The maximum number of backtracking reductions is set to $80$, and the maximum number of inner line-search iterations is set to $10^4$.

\begin{figure}[!htb]
	\vspace{-5pt}
	\centering
	
	\includegraphics[width=0.95\linewidth]{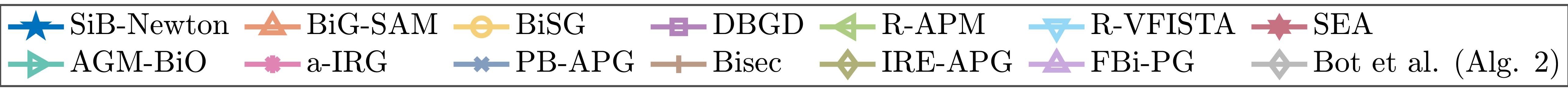}
	
	\vspace{2pt}
	\hspace{-10pt}
	\begin{minipage}[t]{.45\linewidth}
		\includegraphics[width=\textwidth]{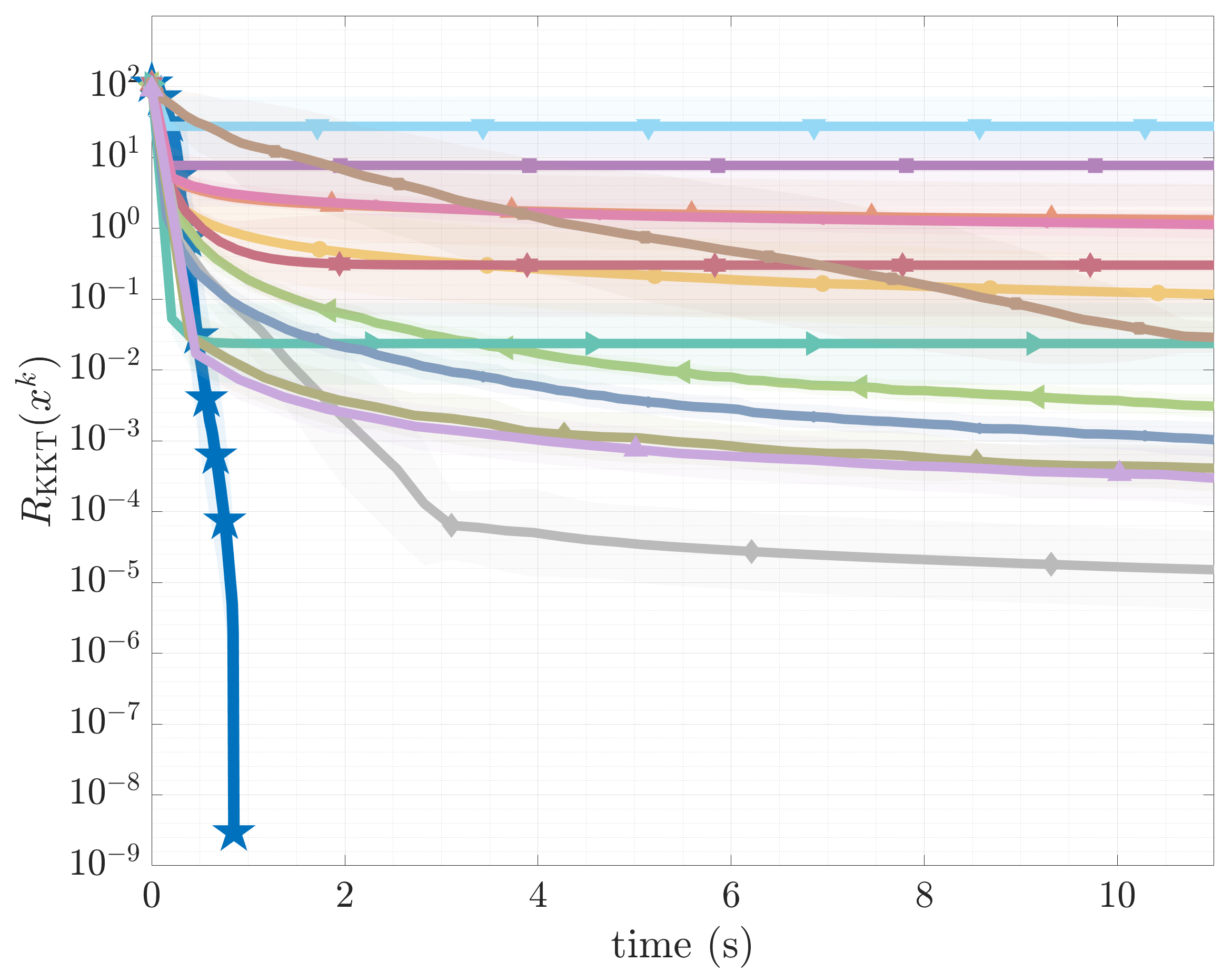}
	\end{minipage}
	\hspace{20pt}
	\begin{minipage}[t]{.45\linewidth}
		\includegraphics[width=\textwidth]{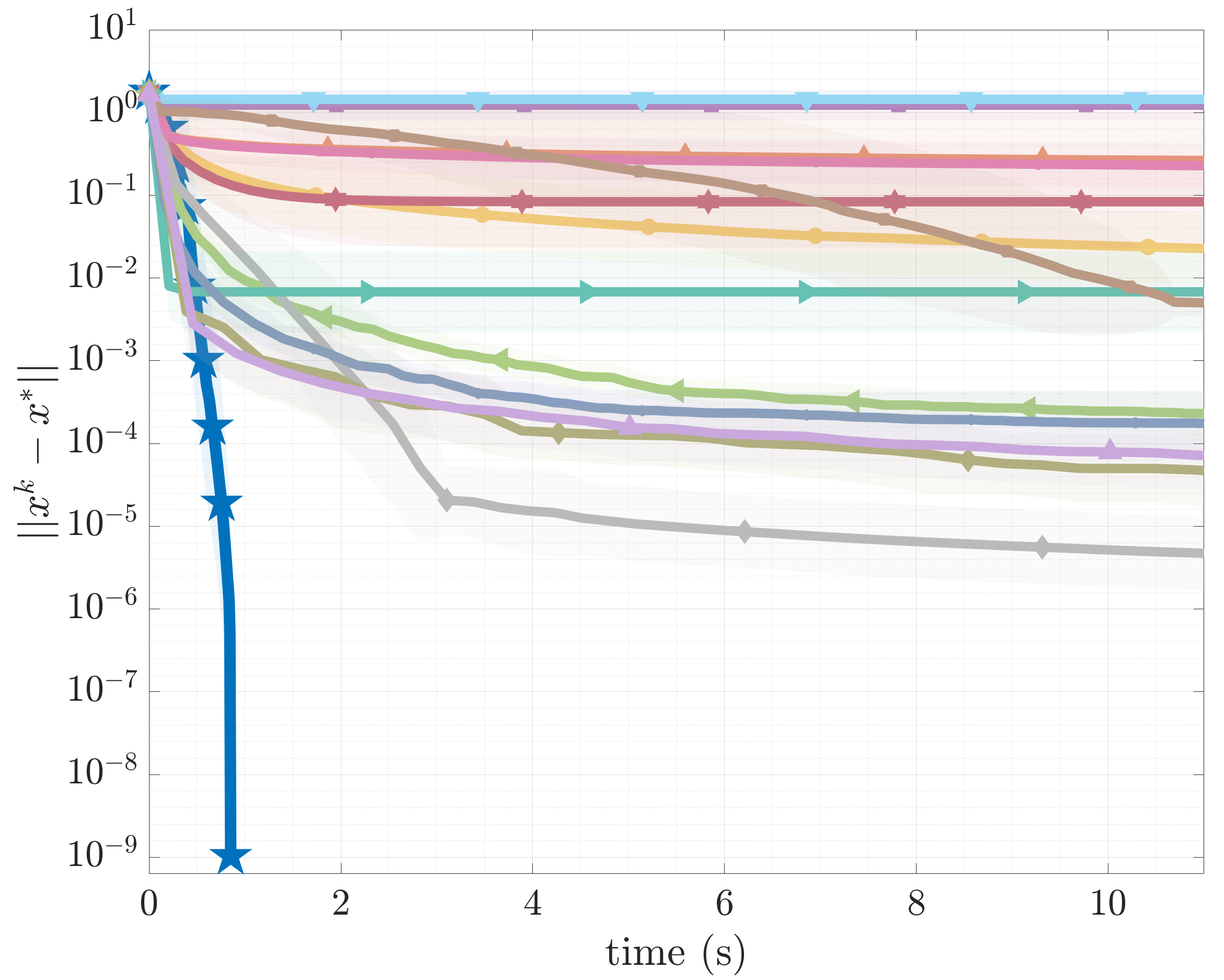}
	\end{minipage}
	
	\vspace{-10pt}
	\caption{Comparison of different methods on the linear regression problem with $m_1=m_2=60$ and $n=100$.}
	\label{fig-1}
	\vspace{-10pt}
\end{figure}

The problem dimensions are set as $n=100$, $m_1=m_2=60$, and $\operatorname{rank}(A_2)=60$. All experiments are repeated over 10 independent trials with random seeds 1\text{-}10. In Figure~\ref{fig-1}, solid lines denote the mean of the running-minimum error sequences over the 10 runs, interpolated onto a common time grid, and shaded areas represent the corresponding one-standard-deviation bands, approximated on the logarithmic scale. On this instance, SiB-Newton achieves substantially smaller KKT residuals and solution errors than the compared methods within a shorter computational time. After approximately 0.5 seconds, both error curves of SiB-Newton decrease rapidly, reaching the order of $10^{-8}$ for both the KKT residual and the solution error within approximately one second, which is consistent with the local superlinear convergence behavior predicted by the analysis. By contrast, the baselines exhibit a markedly slower decay over the same time horizon. In particular, AGM-BiO and R-VFISTA exhibit pronounced error plateaus after an initial decrease, while R-APM continues to improve at a much slower rate. Among the other gradient-type methods, BiSG shows gradual error reduction, whereas a-IRG and SEA exhibit relatively limited improvement after their initial decrease. Although several other methods achieve progressively smaller errors over time, none reaches the accuracy attained by SiB-Newton within the given computational budget. Taken together, Figure~\ref{fig-1} shows that SiB-Newton is both fastest and most accurate, consistent with the predicted local superlinear convergence.

	\bibliographystyle{plain}
	\bibliography{references}
\end{document}